\documentclass[11pt]{amsart}
\usepackage[T1]{fontenc}
\usepackage{lmodern,microtype}
\usepackage{amsmath,amssymb,amsthm,mathtools}
\usepackage[margin=28mm]{geometry}
\usepackage{enumitem}
\usepackage[colorlinks=true,linkcolor=blue,citecolor=blue,urlcolor=blue,hypertexnames=false]{hyperref}
\allowdisplaybreaks
\numberwithin{equation}{section}
\newtheorem{theorem}{Theorem}[section]
\newtheorem{proposition}[theorem]{Proposition}
\newtheorem{lemma}[theorem]{Lemma}
\newtheorem{corollary}[theorem]{Corollary}
\theoremstyle{definition}

\theoremstyle{remark}
\newtheorem{remark}[theorem]{Remark}

\newcommand{\avg}[1]{\langle #1\rangle}

\providecommand{\physnorm}[1]{\lVert #1\rVert}
\providecommand{\physabs}[1]{\lvert #1\rvert}
\title[Infinite-dimensional tori for the perturbed mKdV equation]{Infinite-dimensional invariant tori for the perturbed modified KdV equation}
\author[Y. Wu]{Yuchen Wu}
\address{School of Mathematical Sciences, Fudan University, Shanghai 200433, China}
\email{22110180044@m.fudan.edu.cn}
\author[X. Yuan]{Xiaoping Yuan}
\address{School of Mathematical Sciences, Fudan University, Shanghai 200433, China}
\email{xpyuan@fudan.edu.cn}
\thanks{Xiaoping Yuan is the corresponding author.}
\date{}
\hypersetup{pdftitle={Infinite-dimensional invariant tori for the perturbed modified KdV equation},pdfauthor={Yuchen Wu and Xiaoping Yuan},pdfsubject={Infinite-dimensional invariant tori for the perturbed mKdV equation without external parameters},pdfkeywords={Hamiltonian PDE, KAM theory, Gardner bracket, almost periodic solutions}}
\begin{document}
\begin{abstract}
In this paper, we consider the perturbed modified KdV equation
\(u_t=-u_{xxx}+6u^2u_x+\partial_x(u^5)\) on the circle, without external parameters.
We prove that it admits plenty of almost periodic solutions.
\end{abstract}
\maketitle
\tableofcontents
\section{Introduction}\label{sec:introduction}

The periodic KdV and mKdV equations are completely integrable,
and their Birkhoff coordinates yield infinite-dimensional invariant tori,
including full-dimensional ones~\cite{KappelerPoschel2003,KST2008}.
In this paper, we consider the perturbed modified Korteweg--de Vries
equation
\begin{equation}\label{eq:pde}
 u_t=-u_{xxx}+6u^2u_x+\partial_x(u^5),\qquad
 x\in\mathbb T=\mathbb R/(2\pi\mathbb Z),\qquad \langle u\rangle=0,
\end{equation}
for real-valued functions. We prove the existence of small
infinite-dimensional invariant tori carrying almost periodic solutions.
Equation~\eqref{eq:pde} has no external parameters; the squared Fourier
amplitudes serve as internal parameters.

The equation is Hamiltonian for the Gardner Poisson structure, with
Hamiltonian
\[
 H(u)=\left\langle\frac12u_x^2+\frac12u^4+\frac16u^6\right\rangle.
\]
Thus \eqref{eq:pde} is a Hamiltonian perturbation of mKdV by a quintic
derivative nonlinearity.

\subsection{Context and background}

For nonlinear evolution equations viewed as infinite-dimensional
dynamical systems, a fundamental question is to understand their
dynamical properties, especially over long times. A natural first step
is therefore to seek solutions with recurrent properties, of which
time-periodic solutions are the simplest. Early results on periodic
solutions of nonlinear wave equations go back to
Rabinowitz~\cite{Rabinowitz1967}; see also Br\'ezis, Coron and
Nirenberg~\cite{BCN1980} for a variational approach.

A more general notion is that of a quasi-periodic solution, depending
periodically on finitely many phases. From the point of view of
dynamical systems, finding such solutions amounts to constructing
embedded finite-dimensional tori on which the flow is conjugate to a
linear rotation. The construction of such tori leads to small-divisor
problems. Kolmogorov--Arnold--Moser (KAM) theory provides a powerful
tool for handling these problems under suitable nonresonance and
nondegeneracy assumptions.

Since the properties of dispersive equations vary widely, it is useful
first to distinguish the orders of the linear and nonlinear vector fields.
Consider a Hamiltonian PDE
\[
 \partial_t w=Aw+F(w),\qquad x\in\mathbb T^d,\qquad
 p:=\operatorname{ord}A>0,\quad \widetilde p:=\operatorname{ord}F.
\]
Here $Aw$ is a linear Hamiltonian vector field and $F(w)$ is a nonlinear
Hamiltonian vector field. The order describes the loss of spatial
regularity: for Sobolev indices $v$ in the relevant range, $A$ maps $H^v$
into $H^{v-p}$, while $F$ maps a neighborhood of the origin in $H^v$
into $H^{v-\widetilde p}$. Accordingly, the perturbation is called bounded
when $\widetilde p\leq0$ and unbounded when $\widetilde p>0$; the case
$\widetilde p<0$ corresponds to a regularizing nonlinear vector field.

For bounded perturbations, there is a well-established theory concerning
the existence of quasi-periodic solutions. Typical examples are nonlinear
Schr\"odinger equations with derivative-free nonlinearities and nonlinear
wave equations written as first-order Hamiltonian systems. The pioneering
works of Kuksin~\cite{Kuksin1987} and Wayne~\cite{Wayne1990}, followed by
the KAM theorems of Kuksin and P\"oschel~\cite{KuksinPoschel1996} and
P\"oschel~\cite{Poschel1996}, established the existence of finite-dimensional
invariant tori for such equations. An alternative is the Newton method of Craig and
Wayne~\cite{CraigWayne1993}, developed further by
Bourgain~\cite{Bourgain1994,Bourgain1998,Bourgain2005Book} for the
construction of periodic and quasi-periodic solutions. This is known as
the Craig--Wayne--Bourgain (CWB) method.

For PDEs with unbounded Hamiltonian perturbations, the small-divisor
problem is accompanied by a loss of spatial derivatives. Consequently,
estimating the inverse of the homological operator is not sufficient by
itself: one must also ensure that the resulting changes of variables are
well defined in the phase space. For one-dimensional dispersive models,
the comparison between the nonlinear order $\widetilde p$ and the growth
of the linear frequency gaps, $(n+1)^p-n^p\sim pn^{p-1}$, leads to a
further distinction between the range $0<\widetilde p<p-1$ and the
critical case $0<\widetilde p=p-1$.

In the range $0<\widetilde p<p-1$, the first KAM results for unbounded
perturbations are due to Kuksin. His estimates for small-divisor equations
with large variable coefficients~\cite{Kuksin1997} provide the main tool
for a KAM theorem for equations of KdV type~\cite{Kuksin1998}. In
particular, they can be used to study the persistence of finite-gap
solutions of the periodic KdV equation under Hamiltonian perturbations;
see also Kappeler and P\"oschel~\cite{KappelerPoschel2003}. The perturbed
mKdV equation~\eqref{eq:pde} belongs to this class: its linear part has
order $p=3$, whereas the nonlinear vector field
$6u^2u_x+\partial_x(u^5)$ has order $\widetilde p=1$. Thus the nonlinearity
is unbounded, but it is not critical in the above sense.
We also mention the work of Xu and Yan~\cite{XuYan2016} on smooth
quasi-periodic solutions of perturbed mKdV equations.

The critical case $0<\widetilde p=p-1$ is more delicate. A typical example
is the derivative nonlinear Schr\"odinger equation, for which $p=2$ and
$\widetilde p=1$. In this case, the growth associated with the derivative
loss is of the same order as the linear frequency gaps, and Kuksin's
lemma is no longer applicable in its original form. Liu and
Yuan~\cite{LiuYuan2010,LiuYuan2011} generalized the small-divisor estimates
and obtained quasi-periodic solutions for a class of derivative nonlinear
Schr\"odinger equations with Dirichlet boundary conditions and perturbed
Benjamin--Ono equations with periodic boundary conditions. Their later
work~\cite{LiuYuan2014} treats derivative nonlinear Schr\"odinger equations
with periodic boundary conditions.

A further class consists of quasi-linear and fully nonlinear
perturbations, for which the order of the nonlinear term can equal or
exceed that of the linear part. Of particular relevance to
the present problem are the results of Baldi, Berti and
Montalto~\cite{BBM2016KdV,BBM16}, who constructed quasi-periodic solutions
for autonomous quasi-linear perturbations of KdV and mKdV.
For other representative results, see Baldi, Berti, Haus and
Montalto~\cite{BBHM2018} on gravity water waves in finite depth and
Berti, Hassainia and Masmoudi~\cite{BHM2023} on vortex patches for the
planar Euler equation.

It should be pointed out that the preceding results concern
finite-dimensional invariant tori. Even when tori of every prescribed
finite dimension can be constructed, the smallness conditions generally
depend on that dimension. More precisely, increasing the number of basic
frequencies changes the nonresonance conditions and the constants in the
iteration, and generally requires a smaller perturbation. Therefore, the
existence of infinite-dimensional invariant tori for a fixed equation
does not follow from these results by merely letting the dimension tend
to infinity.

This leads to the more general question of almost periodic solutions,
which can be approximated uniformly in time by quasi-periodic functions
in the phase-space norm. Nonresonant infinite-dimensional invariant
tori provide such solutions with infinitely many independent
frequencies. Their construction requires control of the frequency
corrections, the excluded parameters and the changes of variables as
the number of tangential sites increases. The relevant results are
recalled in Section~\ref{subsec:historical-remarks}.

\subsection{Statement of the main result}

Let $H^v_0=H^v_0(\mathbb T;\mathbb R)$ denote the real Sobolev space of
functions of zero mean, equipped with the norm
\[
 \|u\|_v^2=\sum_{n>0}\langle n\rangle^{2v}|u_n|^2,
 \qquad u_{-n}=\overline{u_n},\qquad
 \langle n\rangle=(1+n^2)^{1/2}.
\]
Every countable product of circles is equipped with the product topology.
Following \cite[Definition~1.1]{BGR25}, a
nonresonant infinite-dimensional Kronecker torus is an invariant set
homeomorphic to such a product, with dynamics conjugate to a rotation
whose frequency sequence has no nonzero finitely supported integer
relation. We use Bohr almost periodicity in the indicated Sobolev norm;
its definition is recalled in Theorem~\ref{thm:main}.

Our main result is stated as follows.

\begin{theorem}
\label{thm:main}
Let \(s\geq2\), let \(\sigma=s+4\), and let
\(S_\infty\subset\mathbb N\) be any infinite set satisfying
\[
 \sum_{j\in S_\infty}\frac1j<\infty.
\]
For every sufficiently small \(\varepsilon>0\) and every
\(0<\vartheta<1\), there exist positive numbers
\(b_j=b_j(\varepsilon,\vartheta)\), \(j\in S_\infty\), with
\[
 \sum_{j\in S_\infty}\langle j\rangle^{2\sigma}b_j<\infty,
\]
and a Borel set
\[
 \mathcal G\subset\mathcal B=
       \prod_{j\in S_\infty}[b_j/2,b_j],\qquad
 \mathbb P(\mathcal G)\geq1-\vartheta,\qquad
 \mathbb P=\bigotimes_{j\in S_\infty}
       \operatorname{Unif}([b_j/2,b_j]).
\]
For each \(\xi\in\mathcal G\), define
\(z_\xi(\theta)\), \(\theta\in\mathbb T^{S_\infty}\), by
\[
 (z_\xi(\theta))_j=
 \begin{cases}
   \sqrt{\xi_j}e^{i\theta_j},&j\in S_\infty,\\
   0,&j>0,\ j\notin S_\infty,
 \end{cases}
 \qquad T_\xi=z_\xi(\mathbb T^{S_\infty}),
\]
with the negative Fourier coefficients determined by the reality condition.
There exist embeddings \(\Psi_\xi:T_\xi\to H^s_0\), real
frequencies \((\omega_j(\xi))_{j\in S_\infty}\), and a constant
\(1\leq C\leq e^{1/64}\), independent of \(\xi\), such that the following assertions hold.
\begin{enumerate}[label=\textup{(\roman*)}]
\item For every \(\xi\in\mathcal G\),
\[
 C^{-1}\|u-v\|_s\leq\|\Psi_\xi(u)-\Psi_\xi(v)\|_s
              \leq C\|u-v\|_s\qquad(u,v\in T_\xi).
\]
The image \(\mathcal T_\xi=\Psi_\xi(T_\xi)\) is compact and
homeomorphic to \(\mathbb T^{\mathbb N}\), both in \(H^s_0\)
and in \(H^{s+3}_0\). Its points belong to \(H^{s+4}_0\), with
\[
 \sup_{w\in\mathcal T_\xi}\|w\|_{s+4}\leq\varepsilon,
\]
and \(\theta\mapsto\Psi_\xi z_\xi(\theta)\) is continuous
from the product topology into \(H^{s+3}_0\).

\item For every initial angle \(\theta\), the function
\[
 u(t)=\Psi_\xi z_\xi(\theta+t\omega(\xi)),\qquad t\in\mathbb R,
\]
is a global solution of \eqref{eq:pde} and belongs to
\(C(\mathbb R;H^{s+3}_0)\cap C^1(\mathbb R;H^s_0)\).
It is Bohr almost periodic in \(H^{s+3}_0\): for each \(a>0\),
the set of \(\tau\in\mathbb R\) satisfying
\(\sup_t\|u(t+\tau)-u(t)\|_{s+3}<a\) is relatively dense.
\item The frequencies are rationally independent:
\[
 \sum_{j\in S_\infty}k_j\omega_j(\xi)\ne0
\]
for every nonzero finitely supported integer sequence \(k\).
Every solution in \textup{(ii)} has a dense forward orbit in \(\mathcal T_\xi\)
with respect to the \(H^{s+3}_0\) topology.
\end{enumerate}
\end{theorem}

The action intervals, and hence the product probability, may depend on
both $\varepsilon$ and $\vartheta$. In particular, for any fixed $a>0$,
the choice $\vartheta=\varepsilon^a$ gives a family of parameter boxes
on which the relative measure of the corresponding good sets tends to
one as $\varepsilon\to0$.

\begin{remark}
We expect our method to apply to the more general family
\[
 u_t=-u_{xxx}+6u^2u_x+\partial_x f(u),
\]
where $f$ is a fixed real-analytic function near the origin with
$f(u)=O(u^4)$. The proposed extension relies on a selective normalization of the
quintic Hamiltonian terms and is expected to involve only additional technical
details, without introducing essential new difficulties.
\end{remark}

\begin{corollary}\label{cor:intro-not-quasiperiodic}
The solutions in Theorem~\ref{thm:main}, regarded as functions with
values in $H^s_0$, are not quasi-periodic with finitely many basic
frequencies.
\end{corollary}
\begin{proof}
Apply \cite[Lemma~10.2]{BGR25} in $H^s_0$ to the homeomorphism
$\varphi\mapsto\Psi_\xi z_\xi(\theta+\varphi)$; rational independence
of the frequencies gives an infinite-dimensional $\mathbb Q$-span.
\end{proof}

\subsection{Historical remarks on almost periodic solutions}\label{subsec:historical-remarks}

The first result on almost periodic solutions for nonintegrable Hamiltonian
PDEs is due to
Bourgain~\cite{Bourgain1996}, who considered a wave equation with a
random potential. A few years later, the same author considered NLS with
a random convolution potential and constructed full-dimensional
invariant tori with a slower decrease of the Fourier
amplitudes~\cite{Bourgain2005}. This result has been revisited and
improved in several directions. Cong, Liu, Shi and
Yuan~\cite{CLSY2018} extended the admissible Gevrey decay and proved
long-time stability of the resulting tori; Cong~\cite{Cong2024}
obtained a substantially slower decay, beyond the Gevrey class. For NLS
on the circle with a convolution potential, Biasco, Massetti and
Procesi~\cite{BMP2021} proved the persistence of both finite- and
infinite-dimensional invariant tori within a single KAM scheme. See
also Cong, Mi, Shi and Wu~\cite{CMSW2019} for spatially dependent
nonlinearities and Cong and Yuan~\cite{CongYuan2021} for nonlinear
wave equations. Further,
using sparsity arguments, Biasco, Massetti and
Procesi~\cite{BMP2023} constructed almost periodic solutions of finite
Sobolev regularity, including nonclassical solutions. The corresponding
tori have sparse reference support and are not full-dimensional.

In these constructions, the external parameters play a crucial role in
controlling the small divisors. In particular, the Fourier coefficients
of the convolution potential provide independent parameters for the
linear frequencies. Improving the regularity of the potential is a
natural step towards removing these parameters. In this direction,
Corsi, Gentile and Procesi~\cite{CGP2023} obtained almost periodic
solutions for convolution potentials of arbitrarily high finite
regularity, while the solutions remain Gevrey in the space variable.
Biasco, Corsi, Gentile and Procesi~\cite{BCGP2024} extended this result
by proving that their set of almost periodic solutions is asymptotically
of full measure.

For a fixed PDE, the parameters must instead be extracted from the
initial data. As pointed out by Bourgain~\cite{Bourgain2005}, fast
decay of the action variables favors convergence but allows less
frequency modulation and leads to worse small divisors. A different
scheme was proposed by P\"oschel~\cite{Poschel2002}: rather than
constructing an infinite-dimensional torus directly, one constructs
iteratively a convergent sequence of invariant tori whose dimensions
tend to infinity. His result concerns NLS with a potential and a
regularizing nonlinearity; see Geng and Xu~\cite{GengXu2013} for a
development using Toeplitz--Lipschitz estimates. With internal
parameters, Berti and Biasco~\cite{BB2011} obtained finite sequences
of bifurcations into higher-dimensional elliptic tori.

Following P\"oschel's scheme, Bernier, Gr\'ebert and Robert~\cite{BGR25}
constructed nonresonant infinite-dimensional invariant tori supporting
almost periodic but not quasi-periodic solutions of NLS without external
parameters. To control the twist, they use a normal form in which
the nonlinearity is regularizing up to a gauge transformation, together
with a sparsity condition on the tangential sites. Bernier and
Procesi~\cite{BP26} removed this sparsity condition by using an
action--Toeplitz decomposition of the Hamiltonian to obtain uniform
estimates on the amplitude--frequency map. In particular, they obtained
full-dimensional tori whose tangent cone at each point is a Lagrangian
subspace.

The preceding results on almost periodic solutions concern bounded or
regularizing nonlinear vector fields. For unbounded nonlinearities,
Geng and Hong~\cite{GengHong2017} obtained full-dimensional invariant
tori for a fifth-order equation of KdV type with external parameters,
and Wu and Yuan~\cite{WuYuan2026} treated a derivative NLS equation
with a convolution potential. We also mention the almost periodic
response solutions of Corsi, Montalto and Procesi~\cite{CMP2021} for
a forced quasilinear Airy equation. In contrast to this response
problem, the equation considered here is autonomous and has no external
parameters.

We summarize the results for one-dimensional Hamiltonian PDEs in
Table~\ref{tab:almost-periodic-results}, following the classification
in~\cite{WuYuan2026}. A dash means that no result is listed in the
corresponding class.
\begin{table}[htbp]
\centering
\small
\renewcommand{\arraystretch}{1.15}
\begin{tabular}{|c|c|c|c|c|}
\hline
 & \multicolumn{2}{c|}{external parameters}
 & \multicolumn{2}{c|}{internal parameters} \\
\cline{2-5}
 & $\tilde p\leq0$ & $\tilde p>0$
 & $\tilde p\leq0$ & $\tilde p>0$ \\
\hline
infinite but non-full dimension
 & \cite{BMP2023} & --- & \cite{BGR25} & this paper \\
\cline{1-3}\cline{4-5}
full dimension
 & \cite{Bourgain1996,Bourgain2005,Poschel2002}$\cdots$
 & \cite{CMP2021,GengHong2017,WuYuan2026} & \cite{BP26} & --- \\
\hline
\end{tabular}
\caption{Infinite-dimensional invariant tori for one-dimensional PDEs.}
\label{tab:almost-periodic-results}
\end{table}

In this paper, we consider the perturbed mKdV equation~\eqref{eq:pde}
and construct sparse infinite-dimensional invariant tori using only
internal parameters. Our proof follows P\"oschel's construction by successive openings~\cite{Poschel2002} and uses the action--Toeplitz frequency analysis
of~\cite{BP26}, adapted to the Gardner Poisson structure and the
derivative nonlinearity. For related normal-form results on the
long-time dynamics of generalized KdV and Benjamin--Ono equations,
we refer to Bernier and Gr\'ebert~\cite{BG21}.

\subsection{Main difficulties and proof strategy}

The finite-stage frequencies have the form
\begin{align*}
 \omega_j&=j^3+j\mathfrak a-6j\xi_j+jf_j,
       &&j\in S,\\
 \Omega_n&=n^3+n\mathfrak a+g_n,
       &&n\in\mathbb N\setminus S.
\end{align*}
Momentum conservation cancels the common term $\mathfrak a$ in the
homological divisors. Following \cite{BP26}, two successive action
limits control the common frequency correction and its variation under
later openings. The derivatives of the decaying remainder with respect to the amplitude
parameters are summable when $\sum_{j\in S_\infty}j^{-1}<\infty$.
Keeping each integrable correction as a function of the parameters present
when it is constructed yields Lipschitz estimates for the frequencies
that are uniform under subsequent openings.

The Gardner bracket contributes a Fourier multiplier to each normal
contraction. We introduce a class of generating Hamiltonians whose vector fields decompose into
a bounded vector field and a real scalar multiple of the generator of
spatial translations in the normal variables.
Translation invariance allows us to integrate the latter by a
state-dependent spatial translation. The resulting Sobolev flows satisfy the pullback identities,
conjugacy relations and estimates required for infinite compositions.

The paper is organized as follows.
Section~\ref{sec:hamiltonian} gives the Hamiltonian formulation and
identifies the quartic resonances. In Section~\ref{sec:verified-coefficients}
we introduce the coefficient norms and estimate the Poisson brackets.
We solve the homological equations in Section~\ref{sec:max4-homological}
and construct the corresponding Hamiltonian flows in Section~\ref{sec:verified-flows}.
Section~\ref{sec:verified-newton} gives the KAM iteration and the cubic
normalization for a fixed tangential set.
Sections~\ref{sec:verified-opening} and~\ref{sec:preparatory-normal-forms}
give the estimates and normal forms required to open new modes.
In Section~\ref{sec:loop} we construct a sequence of finite-dimensional
invariant tori by successive openings. In Section~\ref{sec:xhigh3-outer}
we pass to the limit and prove Theorem~\ref{thm:main}.

\section{Hamiltonian formulation}\label{sec:hamiltonian}

Let $\mathbb T=\mathbb R/(2\pi\mathbb Z)$ and
$\langle f\rangle=(2\pi)^{-1}\int_0^{2\pi}f(x)\,dx$.
For $v\ge0$, we work in the real mean-zero Sobolev space $H^v_0$,
with Fourier coordinates
\[
 u(x)=\sum_{n\ne0}u_ne^{inx},\qquad
 u_{-n}=\overline{u_n},\qquad
 \|u\|_v^2=\sum_{n>0}\langle n\rangle^{2v}|u_n|^2,
 \qquad \langle n\rangle=(1+n^2)^{1/2}.
\]
Equation~\eqref{eq:pde} can be written in Hamiltonian form
$u_t=\partial_x\nabla H(u)$, with
\begin{equation}\label{ham:hamiltonian}
 H(u)=\left\langle\frac12u_x^2+\frac12u^4+\frac16u^6\right\rangle
 =H_2+H_4+H_6,
\end{equation}
where $\nabla H$ denotes the variational gradient and $H_j$ is
homogeneous of degree $j$.

Treating $u_n$ and $u_{-n}$ as independent complex variables, the
Gardner Poisson bracket and Hamiltonian vector field are
\begin{equation}\label{ham:bracket}
 \{F,G\}=i\sum_{n\ne0}n F_{u_n}G_{u_{-n}},\qquad
 (X_G)_n=\{u_n,G\}=inG_{u_{-n}}.
\end{equation}
The corresponding weak symplectic form is
\begin{equation}\label{ham:symplectic}
 \Omega(h,k)=\langle(\partial_x^{-1}h)k\rangle
 =\sum_{n>0}\frac{h_nk_{-n}-h_{-n}k_n}{in},
\end{equation}
where $\partial_x^{-1}$ has Fourier multiplier $1/(in)$ for $n\ne0$.
It is constant and continuous, with
$|\Omega(h,k)|\le2\|h\|_v\|k\|_v$, and is nondegenerate since
$\Omega(h,\partial_x^{-1}h)=\langle(\partial_x^{-1}h)^2\rangle>0$
for $h\ne0$. Our convention is $\Omega(X_G,h)=dG[h]$, so
$X_G=\partial_x\nabla G$. These identities are first understood on
smooth functions; their Sobolev realization is established in
Section~\ref{sec:verified-flows}.

Set $\rho_n=u_nu_{-n}$ for $n>0$. Then
\[
 H_2=\sum_{n>0}n^2\rho_n,\qquad
 M=\sum_{n>0}\rho_n=\frac12\langle u^2\rangle.
\]
Since $(X_M)_n=inu_n$, the Hamiltonian $M$ generates the spatial
translations $\tau_bu(x)=u(x+b)$.

\begin{lemma}[Quartic resonant projection]
The part of $H_4$ resonant with $H_2$ is
\begin{equation}\label{ham:quartic}
 Z_4=6M^2-3\sum_{n>0}\rho_n^2.
\end{equation}
The angular frequencies of $H_2+Z_4$ are
\begin{equation}\label{ham:frequencies}
 \Omega_n=n\partial_{\rho_n}(H_2+Z_4)
 =n^3+12nM-6n\rho_n\qquad(n>0),
\end{equation}
with $\Omega_{-n}=-\Omega_n$. The normalizing transformation with a
finite cutoff is constructed in Theorem~\ref{thm:max-initial-cutoff}.
\end{lemma}
\begin{proof}
The ordered Fourier coefficient of $H_4$ is $1/2$ on tuples
$(a,b,c,d)$ with nonzero entries and $a+b+c+d=0$. The linear flow
multiplies the corresponding monomial by $e^{it(a^3+b^3+c^3+d^3)}$.
Since
\[
 a^3+b^3+c^3+d^3=-3(a+b)(a+c)(b+c),
\]
a quartic monomial is resonant precisely when its indices form two
opposite pairs. There are $24$ permutations of $(n,-n,m,-m)$ for
$0<m<n$ and $6$ permutations of $(n,n,-n,-n)$. Thus
\[
 Z_4=3\sum_{n>0}\rho_n^2+12\sum_{0<m<n}\rho_m\rho_n
 =6M^2-3\sum_{n>0}\rho_n^2.
\]
For an action Hamiltonian $K$, $(X_K)_n=in(\partial_{\rho_n}K)u_n$,
which gives \eqref{ham:frequencies}.
\end{proof}

Fix a finite tangential set $S\subset\mathbb N$ and positive parameters
$\xi=(\xi_j)_{j\in S}$. Introduce partial action--angle variables by
\begin{equation}\label{ham:polar}
 u_j=\sqrt{\xi_j+y_j}\,e^{i\theta_j},\qquad
 u_{-j}=\sqrt{\xi_j+y_j}\,e^{-i\theta_j}\qquad(j\in S),
\end{equation}
where $\xi_j+y_j>0$, and retain $u_n$ as normal variables for
$n\in\mathcal N_S:=\mathbb Z\setminus(S\cup(-S)\cup\{0\})$.
Since $\{\theta_j,y_k\}=j\delta_{jk}$, the bracket becomes
\begin{equation}\label{ham:polar-bracket}
 \{F,G\}=\sum_{j\in S}j(F_{\theta_j}G_{y_j}-F_{y_j}G_{\theta_j})
 +i\sum_{n\in\mathcal N_S}nF_{u_n}G_{u_{-n}}.
\end{equation}
Here $\rho_j=\xi_j+y_j$ is the squared Fourier amplitude; the canonical
action is $I_j=\rho_j/j$, for which $\{\theta_j,I_k\}=\delta_{jk}$.

A Hamiltonian invariant under spatial translations and $u\mapsto-u$
contains only monomials $e^{ik\cdot\theta}y^m u_{n_1}\cdots u_{n_q}$
satisfying
\begin{equation}\label{ham:symmetries}
 \pi_S(k)+\sum_{a=1}^q n_a=0,\qquad
 \sum_{j\in S}k_j+q\equiv0\pmod2,
 \qquad \pi_S(k)=\sum_{j\in S}jk_j.
\end{equation}
Indeed, translation by $b$ multiplies such a monomial by
$e^{ib(\pi_S(k)+\sum_a n_a)}$, whereas sign reversal multiplies it by
$(-1)^{\sum_jk_j+q}$. Both transformations preserve the Poisson bracket,
so this class is closed under the bracket. The weighted Taylor degree
of the monomial is $2|m|_1+q$; the angular Fourier order is not included.

\section{Spaces of Hamiltonians}
\label{sec:verified-coefficients}
We introduce spaces of Hamiltonians with two successive action limits,
following the action--Toeplitz construction of Bernier and
Procesi~\cite[Sections~5.1--5.2]{BP26}. The coefficient weights below are
adapted to the Fourier multiplier in the Gardner bracket.

\subsection{Norms and action limits}
\label{sec:coefficient-definitions}

\medskip\noindent\textit{Fourier--Taylor expansions.}
Fix a finite tangential set $S\subset\mathbb N$, and write
\[
 \mathcal N=\mathbb Z\setminus(S\cup(-S)\cup\{0\}),\qquad
 \mathcal N_+=\mathcal N\cap\mathbb N,\qquad J=\max(\{1\}\cup S).
\]
The variables are $\theta=(\theta_j)_{j\in S}$, $y=(y_j)_{j\in S}$,
and $u=(u_n)_{n\in\mathcal N}$; $u_n$ and $u_{-n}$ are treated as
independent complex variables. For $k\in\mathbb Z^S$ and
$m\in\mathbb N_0^S$, where $\mathbb N_0=\{0,1,2,\ldots\}$, set
\[
 \begin{gathered}
 |k|_1=\sum_{j\in S}|k_j|,\qquad
 |m|_1=\sum_{j\in S}m_j,\qquad |m|_\infty=\max_{j\in S}m_j,\\
 k\cdot\theta=\sum_{j\in S}k_j\theta_j,\qquad
 y^m=\prod_{j\in S}y_j^{m_j},\qquad
 \pi(k)=\sum_{j\in S}jk_j.
 \end{gathered}
\]
For a normal tuple $\mathbf j=(j_1,\ldots,j_q)\in\mathcal N^q$, put
$u_{\mathbf j}=\prod_{a=1}^q u_{j_a}$ and
$|\mathbf j|_\infty=\max_{1\le a\le q}|j_a|$.
Empty sums and maxima are zero, empty products are one, and
$\mathcal N^0$ consists of the empty tuple.

A formal Hamiltonian is a Fourier--Taylor series
\begin{equation}\label{eq:coefficient-expansion}
 P=\sum_{\substack{k\in\mathbb Z^S,\ m\in\mathbb N_0^S\\q\ge0}}
       \sum_{\mathbf j\in\mathcal N^q}
       p_{k,m,q}(\mathbf j)e^{ik\cdot\theta}y^m u_{\mathbf j}.
\end{equation}
For each fixed $(k,m,q)$, the coefficient array
$p_{k,m,q}:\mathcal N^q\to\mathbb C$ is symmetric:
\[
 p_{k,m,q}(j_{\sigma(1)},\ldots,j_{\sigma(q)})
 =p_{k,m,q}(j_1,\ldots,j_q)
 \qquad(\sigma\in\mathfrak S_q),
\]
where $\mathfrak S_q$ is the permutation group on $q$ positions.
The inner sum in \eqref{eq:coefficient-expansion} runs over ordered tuples,
including repeated indices, and contains no factor $1/q!$.
Thus, if
$\nu_n(\mathbf j)=\#\{a:j_a=n\}$, the coefficient of the monomial
$e^{ik\cdot\theta}y^m\prod_{n\in\mathcal N}u_n^{\nu_n(\mathbf j)}$
after collecting equal terms is
\[
 \frac{q!}{\prod_{n\in\mathcal N}\nu_n(\mathbf j)!}\,
 p_{k,m,q}(\mathbf j).
\]
We call the array with fixed $(k,m,q)$ a coefficient block and often
denote it by $p$. Throughout this section we impose momentum conservation:
\[
 p_{k,m,q}(\mathbf j)=0
 \quad\hbox{unless}\quad \pi(k)+\sum_{a=1}^q j_a=0.
\]
The coefficient spaces are complex. The reality condition is
$p_{-k,m,q}(-\mathbf j)=\overline{p_{k,m,q}(\mathbf j)}$, where
$-\mathbf j=(-j_1,\ldots,-j_q)$. It ensures real values on
$\theta,y\in\mathbb R^S$, $u_{-n}=\overline{u_n}$ wherever the series
converges. When evenness is imposed, the
coefficients also vanish unless $q+\sum_{j\in S}k_j$ is even, as in
\eqref{ham:symmetries}.

\medskip\noindent\textit{Projections.}
The normal degree of $e^{ik\cdot\theta}y^m u_{\mathbf j}$ is $q$, and
its weighted Taylor degree is $2|m|_1+q$. In particular, angular Fourier
indices do not contribute to the weighted Taylor degree.
For a set $A\subset\mathbb N_0$, define the degree projection by
\[
 (\Pi_{\deg\in A} P)_{k,m,q}(\mathbf j)
 =\mathbf1_{\{2|m|_1+q\in A\}}p_{k,m,q}(\mathbf j),
\]
where $\mathbf1_E$ denotes the indicator of a condition $E$.
For $d\in\mathbb N_0$, we use the abbreviations
\[
 \begin{gathered}
 \Pi_d=\Pi_{=d}=\Pi_{\deg\in\{d\}},\qquad
 \Pi_{\le d}=\Pi_{\deg\in\{0,\ldots,d\}},\qquad
 \Pi_{\ge d}=\Pi_{\deg\in\{d,d+1,\ldots\}},\\
 \Pi_{4,5}=\Pi_4+\Pi_5.
 \end{gathered}
\]
A tuple is called balanced if
$\nu_n(\mathbf j)=\nu_{-n}(\mathbf j)$ for every $n\in\mathcal N_+$;
the empty tuple is balanced. The integrable projection is defined by
\begin{equation}\label{eq:coefficient-integrable-projection}
 (\Pi_{\mathrm{int}}P)_{k,m,q}(\mathbf j)
 =\mathbf1_{\{k=0\}}
  \mathbf1_{\{\nu_n(\mathbf j)=\nu_{-n}(\mathbf j)\ \forall n\in\mathcal N_+\}}
  p_{k,m,q}(\mathbf j).
\end{equation}
Thus $\Pi_{\mathrm{int}}P$ depends only on $y$ and the products
$\rho_n=u_nu_{-n}$, $n\in\mathcal N_+$. We call $P$ integrable when
$\Pi_{\mathrm{int}}P=P$, and write $1-\Pi_{\mathrm{int}}$ for the
complementary projection. Both projections commute with the degree
projections.

\medskip\noindent\textit{Action limits.}
For $q\ge2$, let $\mathcal C_{q,1}$ be the set of all permutations of
\[
 (n,-n,\mathbf v),\qquad
 n\in\mathcal N_+,\quad \mathbf v\in\mathcal N^{q-2},\quad
 n>|\mathbf v|_\infty.
\]
For $q\ge4$, let $\mathcal C_{q,2}$ be the set of all permutations of
\[
 (n,-n,h,-h,\mathbf v),\qquad
 n,h\in\mathcal N_+,\quad \mathbf v\in\mathcal N^{q-4},\quad
 n>h>|\mathbf v|_\infty.
\]
Set $\mathcal C_{q,1}=\varnothing$ for $q<2$ and
$\mathcal C_{q,2}=\varnothing$ for $q<4$.
These sets consist of tuples whose two or four largest entries in
absolute value form, respectively, one or two opposite pairs of distinct
magnitudes. The strict inequalities uniquely determine the marked
pairs, independently of the ordering of the tuple.

Fix a coefficient block $p=p_{k,m,q}$. Its first and second action limits
are the arrays
\[
 \begin{aligned}
 U(\mathbf v)=(L_1p)(\mathbf v)
   &=\lim_{\substack{n\to\infty\\n\in\mathcal N_+}}
          p(n,-n,\mathbf v),
       &&\mathbf v\in\mathcal N^{q-2},\\
 V(\mathbf v)=(L_2p)(\mathbf v)
   &=\lim_{\substack{h\to\infty\\h\in\mathcal N_+}}
          U(h,-h,\mathbf v),
       &&\mathbf v\in\mathcal N^{q-4}.
 \end{aligned}
\]
We require these limits to exist pointwise, with $k,m$ and the remaining
tuple fixed. In the second limit, $n$ tends to infinity first and then
$h$ tends to infinity; equivalently, $L_2p=L_1(L_1p)$.
We set $L_1p=0$ if $q<2$ and $L_2p=0$ if $q<4$.

The operators $T_1,T_2$ extend these lower-degree arrays to arrays of
normal degree $q$:
\[
 \begin{aligned}
 (T_1p)(\mathbf j)&=
 \begin{cases}
 U(\mathbf v),&\mathbf j\text{ is a permutation of }(n,-n,\mathbf v)
                         \in\mathcal C_{q,1},\\
 0,&\mathbf j\notin\mathcal C_{q,1},
 \end{cases}\\[2pt]
 (T_2p)(\mathbf j)&=
 \begin{cases}
 V(\mathbf v),&\mathbf j\text{ is a permutation of }(n,-n,h,-h,\mathbf v)
                         \in\mathcal C_{q,2},\\
 0,&\mathbf j\notin\mathcal C_{q,2}.
 \end{cases}
 \end{aligned}
\]
Symmetry of $U,V$ and uniqueness of the marked pairs make these
definitions unambiguous. Set
\[
 p_0=p-T_1p,\qquad p_1=T_1p-T_2p,\qquad p_2=T_2p,
 \qquad p=p_0+p_1+p_2.
\]
The component $p_0$ has zero first action limit, $p_1$ is independent of
the highest marked pair and has zero second action limit, and $p_2$ is
independent of both marked pairs on its support.

\medskip\noindent\textit{The Hamiltonian norm.}
Define weights on $\mathcal N^q$ by
\[
 \begin{aligned}
 w_0(\mathbf j)&=
 \begin{cases}
 1+n,&\mathbf j\sim(n,-n,\mathbf v)\in\mathcal C_{q,1},\\
 1,&\mathbf j\notin\mathcal C_{q,1},
 \end{cases}
 &w_1(\mathbf j)&=
 \begin{cases}
 1+h,&\mathbf j\sim(n,-n,h,-h,\mathbf v)\in\mathcal C_{q,2},\\
 1,&\mathbf j\notin\mathcal C_{q,2},
 \end{cases}\\
 w_2(\mathbf j)&=1,
 \end{aligned}
\]
where $\sim$ means equality up to a permutation. The norm of a
coefficient block is
\[
 a(p)=\max_{0\le b\le2}\sup_{\mathbf j\in\mathcal N^q}
                w_b(\mathbf j)|p_b(\mathbf j)|.
\]
We also write $a_S(p)=a(p)$ when the dependence on $S$ is relevant.
In particular, it controls the remainders through
\[
 \begin{aligned}
 |p(n,-n,\mathbf v)-U(\mathbf v)|
      &\le\frac{a(p)}{1+n} &&(n>|\mathbf v|_\infty),\\
 |U(h,-h,\mathbf v)-V(\mathbf v)|
      &\le\frac{a(p)}{1+h} &&(h>|\mathbf v|_\infty).
 \end{aligned}
\]
For each block, put
\[
 D=|k|_1+2|m|_1+q,\qquad W=\max(1,|k|_1,|m|_\infty,q).
\]
For $\eta\ge0$ and $r>0$, the space $\mathcal H_{\eta,r}(S)$ consists
of the symmetric, momentum-preserving series
\eqref{eq:coefficient-expansion} with the above action limits and finite norm
\begin{equation}\label{eq:coefficient-hamiltonian-norm}
 \|P\|_{\mathcal H_{\eta,r}(S)}
 =\sum_{k,m,q}W e^{\eta D}r^{2|m|_1+q}a(p_{k,m,q}).
\end{equation}
We omit $S$ from the norm when it is fixed.

\medskip\noindent\textit{Generating Hamiltonians.}
For $\mathbf j\in\mathcal N^q$, let
$j_1^*\ge\cdots\ge j_q^*\ge0$ be the decreasing rearrangement of
$|j_1|,\ldots,|j_q|$, counted with multiplicity, and set
\[
 \mu_4(\mathbf j)=
 \begin{cases}j_4^*,&q\ge4,\\0,&q<4.\end{cases}
 \qquad
 \mu_4(P)=\sup_{p_{k,m,q}(\mathbf j)\ne0}\mu_4(\mathbf j).
\]
Here and below, a supremum of nonnegative quantities over an empty set
is zero; in particular, $\mu_4(0)=0$.

Let $\mu\ge0$. A generating Hamiltonian
$\chi\in\mathcal C_{G,\eta,r}(S;\mu)$ is an element of
$\mathcal H_{\eta,r}(S)$ with coefficient arrays $c_{k,m,q}$ such that
\[
 c_{k,m,q}(\mathbf j)=0
 \quad\text{if }q>5\text{ or }\mu_4(\mathbf j)>\mu,
\]
and with finite norm specified below. Decompose each block
$c=c_{k,m,q}$ as $c=c_0+c_1+c_2$ using the same operators $T_1,T_2$.
The fourth-index cutoff implies $L_2c=0$, hence $c_2=0$: for fixed
$\mathbf v$, the coefficient $c(n,-n,h,-h,\mathbf v)$ vanishes whenever
$n>h>\max(\mu,|\mathbf v|_\infty)$.
Define
\[
 b(c)=\max\left\{
 \begin{aligned}
 &\sup_{\mathbf j\notin\mathcal C_{q,1}}
        (1+|\mathbf j|_\infty)^2|c_0(\mathbf j)|,\\
 &\sup_{\mathbf j\in\mathcal C_{q,1}}
        (1+|\mathbf j|_\infty)|c_0(\mathbf j)|,\qquad
   \sup_{\mathbf j\in\mathcal N^q}w_1(\mathbf j)|c_1(\mathbf j)|
 \end{aligned}\right\},
\]
and
\begin{equation}\label{eq:coefficient-generating-norm}
 \|\chi\|_{\mathcal C_{G,\eta,r}(S;\mu)}
 =\sum_{k,m,q}W e^{\eta D}r^{2|m|_1+q}b(c_{k,m,q}).
\end{equation}
The general generating class places no restriction on the degree in
$y$. Its affine subclass $\mathcal C^{\mathrm{aff}}_{G,\eta,r}(S;\mu)$
consists of the elements satisfying
\[
 c_{k,m,q}=0\quad\text{unless}\quad
 |m|_1\le1\quad\text{and}\quad 2|m|_1+q\le5.
\]
Thus its elements are affine in $y$ and have weighted Taylor degree at
most five. This subclass carries the restriction of
\eqref{eq:coefficient-generating-norm}; we suppress $\eta,r,S,\mu$ when
these parameters are fixed.

\medskip\noindent\textit{Poisson bracket and parameter norms.}
We use the Gardner bracket in partial action--angle variables,
\[
 \{P,Q\}_G=i\sum_{n\in\mathcal N}nP_{u_n}Q_{u_{-n}}
 +\sum_{j\in S}j(P_{\theta_j}Q_{y_j}-P_{y_j}Q_{\theta_j}).
\]
For a family $P:O\to\mathcal X$ on a metric parameter set
$(O,d_O)$, where $\mathcal X$ is any of the preceding normed spaces
with fixed $S,\eta,r,\mu$ as applicable, define
\[
 \begin{aligned}
 \|P\|^{\sup}_{\mathcal X;O}
     &=\sup_{\xi\in O}\|P(\xi)\|_{\mathcal X},\\
 \|P\|^{\mathrm{lip}}_{\mathcal X;O}
     &=\sup_{\substack{\xi,\zeta\in O\\\xi\ne\zeta}}
       \frac{\|P(\xi)-P(\zeta)\|_{\mathcal X}}{d_O(\xi,\zeta)},\\
 \|P\|^{\mathrm{tot}}_{\mathcal X;O}
     &=\|P\|^{\sup}_{\mathcal X;O}
        +r^2\|P\|^{\mathrm{lip}}_{\mathcal X;O}.
 \end{aligned}
\]
The Lipschitz seminorm is zero when $O$ has at most one point, and both
terms are zero when $O=\varnothing$. Subscripts are omitted when the
space and parameter set are clear.

\begin{theorem}[The initial perturbation]
\label{thm:coefficient-initial}
The degree projections, $\Pi_{\mathrm{int}}$, $1-\Pi_{\mathrm{int}}$
and their compositions are contractions on $\mathcal H_{\eta,r}(S)$,
both in the Hamiltonian norm and in the total parameter norm.

For \(S=\varnothing\), put
\[
 H_4=\tfrac12\langle u^4\rangle,\quad
 H_6=\tfrac16\langle u^6\rangle,\quad
 \rho_n=u_nu_{-n},\quad M=\sum_{n>0}\rho_n,\quad
 Z_4=6M^2-3\sum_{n>0}\rho_n^2.
\]
Here \(u(x)=\sum_{n\ne0}u_ne^{inx}\) formally and
\(\langle\cdot\rangle\) extracts its zero Fourier coefficient.
Then $R=H_4-Z_4+H_6$ is real and even, belongs to
$\mathcal H_{\eta,r}(\varnothing)$, and satisfies
\[
 \Pi_{\le3}R=\Pi_{\rm int}\Pi_{=4}R=0,\qquad
 \|R\|_{\mathcal H_{\eta,r}}=2e^{4\eta}r^4+e^{6\eta}r^6.
\]
In particular, $\|R\|_{\mathcal H_{\eta,r}}\le r^{4-2\delta}$ for
$\delta>0$ and
\[
 0<r\leq\min\left(1,(2e^{4\eta}+e^{6\eta})^{-1/(2\delta)}\right).
\]
\end{theorem}
\begin{proof}
The degree and integrable projections commute with $T_1,T_2$: the former
are constant on each block, and the latter is unchanged by inserting an
opposite pair. Their contractivity follows from the definition of the norms.

By \eqref{ham:quartic}, the ordered coefficient of $Q=H_4-Z_4$ is
\[
 q(\mathbf j)=\frac12\mathbf1_{\{\sum_a j_a=0\}}
                   \mathbf1_{\{\mathbf j\text{ is not balanced}\}}.
\]
A zero-momentum quartet in $\mathcal C_{4,1}$ is balanced, so
$T_1q=T_2q=0$ and $a(Q)=1/2$.

Writing $p$ for the ordered coefficient of $H_6$, we have
\[
 (p_0,p_1,p_2)
 =\frac16\mathbf1_{\{\sum_a j_a=0\}}
 \bigl(\mathbf1_{\mathcal N^6\setminus\mathcal C_{6,1}},
       \mathbf1_{\mathcal C_{6,1}\setminus\mathcal C_{6,2}},
       \mathbf1_{\mathcal C_{6,2}}\bigr).
\]
Thus $a(H_6)=1/6$, and
\[
 \|R\|_{\mathcal H_{\eta,r}}
 =4e^{4\eta}r^4\frac12+6e^{6\eta}r^6\frac16
 =2e^{4\eta}r^4+e^{6\eta}r^6.
\]
The symmetries follow from the coefficient formulas, and, for $r\le1$,
\[
 \|R\|_{\mathcal H_{\eta,r}}
 \le(2e^{4\eta}+e^{6\eta})r^4\le r^{4-2\delta}
\]
under the stated smallness condition.
\end{proof}

\subsection{Poisson brackets and Lie series}

\begin{theorem}[Poisson brackets and Lie series]
\label{thm:max-coefficient-algebra}
The space $\mathcal H_{\eta,r}(S)$ is Banach. Let $\eta\ge0$, $r>0$,
$\mu\ge0$ and $0<\delta\le\alpha\le1$. For
$P\in\mathcal H_{\eta+\delta,r}(S)$ and
$\chi\in\mathcal C_{G,\eta+\alpha,r}(S;\mu)$, we have
\[
 \|\{P,\chi\}_G\|_{\mathcal H_{\eta,r}}
 \leq\frac{C_S(1+\mu)^2}{r^2\alpha^2\delta}
 \|P\|_{\mathcal H_{\eta+\delta,r}}
 \|\chi\|_{\mathcal C_{G,\eta+\alpha,r}(S;\mu)}.
\]
Here $C_S$ depends only on $S$. The estimate also holds in total norm
on any fixed metric parameter set, with $S,\mu,r$ fixed.

Consequently, if \(0<a\leq1/2\), \(P\in\mathcal H_{\eta+a,r}\), \(\chi\in\mathcal C_{G,\eta+2a,r}(S;\mu)\), and
\[
 \operatorname{ad}_\chi P=\{P,\chi\}_G,\qquad
 K=\frac{C_S(1+\mu)^2}{a^3r^2}\|\chi\|_{\mathcal C_{G,\eta+2a,r}},
\]
then, for every integer \(n\geq1\),
\[
\|\operatorname{ad}_\chi^nP\|_{\mathcal H_{\eta,r}}
\leq(Kn)^n\|P\|_{\mathcal H_{\eta+a,r}}.
\]
If \(eK\leq1/2\), the series
\[
\sum_{n\geq0}\frac{t^n}{n!}\operatorname{ad}_\chi^nP
\]
converges absolutely in \(\mathcal H_{\eta,r}\), uniformly for complex \(|t|\leq1\), with norm at most \(2\|P\|_{\mathcal H_{\eta+a,r}}\). These assertions also hold in total norm.

\end{theorem}
\begin{proof}
We first establish the estimate for a pair of coefficient blocks and
then sum over their Fourier--Taylor labels. All brackets in this proof
are defined coefficientwise. Their identification with the Poisson
bracket of the corresponding Hamiltonian functions is proved in
Proposition~\ref{prop:flows-evaluation}, under the hypotheses stated there.
We write
\[
 \{F,G\}_{\mathrm{nor}}=i\sum_{n\in\mathcal N}nF_{u_n}G_{u_{-n}},
 \qquad
 \{F,G\}_{\mathrm{tan}}=
 \sum_{j\in S}j(F_{\theta_j}G_{y_j}-F_{y_j}G_{\theta_j}),
\]
so that $\{F,G\}_G=\{F,G\}_{\mathrm{nor}}+\{F,G\}_{\mathrm{tan}}$.
We denote the $a$th summand of the tangential bracket by
$\{F,G\}_{\mathrm{tan},a}$.
The constants $C$ and $C_S$ may change from line to line; the former is
numerical, and the latter depends only on $S$.

\medskip\noindent\textit{1. Completeness and cutoffs.}
Let $p$ be a coefficient array, with $U=L_1p$ and $V=L_2p$.
The component decomposition gives
\[
 \sup|p|\le3a(p),\qquad \sup|U|\le2a(p),\qquad \sup|V|\le a(p),
\]
and
\[
 \begin{aligned}
 |p(n,-n,\mathbf v)-U(\mathbf v)|
   &\le\frac{a(p)}{1+n} &&(n>|\mathbf v|_\infty),\\
 |U(h,-h,\mathbf v)-V(\mathbf v)|
   &\le\frac{a(p)}{1+h} &&(h>|\mathbf v|_\infty).
 \end{aligned}
\]
Indeed, the two differences are respectively $p_0(n,-n,\mathbf v)$
and $p_1(n,-n,h,-h,\mathbf v)$, with $n>h$ in the latter identity.
For a Cauchy sequence in $a$, each component converges in its weighted
supremum norm. Passing to the limit in these inequalities identifies
$U,V$ as the successive action limits of the limiting array. Symmetry
and momentum conservation are closed conditions. Each coefficient
space is therefore complete, and so is their weighted $\ell^1$ sum
$\mathcal H_{\eta,r}(S)$.

Consider now a generator block $c=c_{k_\chi,m_\chi,s}$, where $s\le5$.
The fourth-index cutoff implies
\[
 c(n,-n,h,-h,\mathbf v)=0
 \quad\text{if }n>h>\max(\mu,|\mathbf v|_\infty).
\]
Taking the limits in this order gives $L_2c=0$. Likewise, $L_1c(\mathbf v)=0$
if the second-largest absolute entry of $\mathbf v$ exceeds $\mu$.
Consequently, $c_1=T_1c$ and $c_0=c-c_1$ both satisfy the same cutoff.
If either component is nonzero at $(n,-n,\mathbf v)\in\mathcal C_{s,1}$,
then $\mathbf v$ has at most three entries, all but its largest entry
have modulus at most $\mu$, and
\[
 \sum_a v_a=-\pi(k_\chi),\qquad
 |\mathbf v|_\infty\le2\mu+|\pi(k_\chi)|.
\]
We shall use the block-dependent threshold
\[
 N_\chi:=2+J+2\mu+|\pi(k_\chi)|.
\]
In particular, all entries remaining after removal of the largest
opposite pair are strictly smaller than $N_\chi$ in absolute value.

\medskip\noindent\textit{2. Bracket coefficients.}
Fix a Fourier--Taylor block of $P$ and one of $\chi$, with indices
$(k_P,m_P,q)$ and $(k_\chi,m_\chi,s)$ and coefficient arrays $p$ and $c$.
Thus $q,s$ denote normal degrees, whereas $p,c$ always denote arrays.
For a set of positions $A\subset\{1,\ldots,d\}$, let $\mathbf j_A$
be the corresponding subtuple of $\mathbf j$, and let $A^c$ be its
complement. Put
\[
 k=k_P+k_\chi,\qquad A_0=a(p),\qquad B_0=b(c).
\]
For $q,s\ge1$, the normal bracket has normal degree $d=q+s-2$,
Taylor label $m_P+m_\chi$, and ordered coefficient
\begin{equation}\label{eq:coefficient-normal-contraction}
 R(\mathbf j)=\frac{iqs}{\binom d{q-1}}
 \sum_{\substack{A\subset\{1,\ldots,d\}\\|A|=q-1}}
 \ell_A p(\mathbf j_A,\ell_A)c(\mathbf j_{A^c},-\ell_A),
 \qquad
 \ell_A=-\pi(k_P)-\sum_{a\in A}j_a.
\end{equation}
The coefficient is zero unless $\pi(k)+\sum_a j_a=0$, and terms with
$\ell_A\notin\mathcal N$ are omitted. Differentiation gives the factor
$qs$, momentum conservation determines the contracted index $\ell_A$,
and symmetrization gives the average over subsets of positions.
In particular, no unrestricted sum over contracted indices remains.

For $a\in S$, the tangential contraction has normal degree $d=q+s$,
Taylor label $m_P+m_\chi-e_a$, and coefficient
\begin{equation}\label{eq:coefficient-tangential-contraction}
 \frac{ia(k_{P,a}m_{\chi,a}-m_{P,a}k_{\chi,a})}{\binom d q}
 \sum_{\substack{A\subset\{1,\ldots,d\}\\|A|=q}}
 p(\mathbf j_A)c(\mathbf j_{A^c}),
\end{equation}
where $e_a$ is the unit multi-index in the $a$th Taylor coordinate.
A negative Taylor exponent gives a zero contribution. The sum of the
absolute derivative factors is bounded by
\begin{equation}\label{eq:coefficient-tangential-factor}
 T:=J\bigl(|k_P|_1|m_\chi|_\infty+
          |k_\chi|_1|m_P|_\infty\bigr).
\end{equation}
When several output labels occur, a block estimate below means the sum
of the corresponding coefficient seminorms. We also write $a(F)=a(f)$
for a single block $F$ with coefficient array $f$.

For brevity, brackets involving arrays such as $\{U,c_0\}_{\mathrm{nor}}$
denote the ordered coefficients obtained by the same formulas: $U$
retains the labels $(k_P,m_P)$ and has normal degree $q-2$; $V$ retains
these labels and has degree $q-4$. An array of negative degree, or a
contraction requiring more positions than are available, is understood
to be zero. All ratios below are used only when their required positions
exist; otherwise the corresponding contribution is zero.

\medskip\noindent\textit{3. First action limit.}
In this step, let $R$ be the normal bracket coefficient
\eqref{eq:coefficient-normal-contraction} with $c$ replaced by $c_0$.
The generator weights give $|\ell c_0|\le B_0$, and hence
\[
 \sup|R|\le3qsA_0B_0.
\]
For an output tuple $(n,-n,\mathbf v)$ with $n>|\mathbf v|_\infty$,
we separate the summands according to the positions of the marked
entries $n,-n$.

\emph{(i) One marked index lies in each factor.}
For a fixed summand, write $\ell=\ell_A$ and define
\[
 N:=|(\mathbf j_{A^c},-\ell)|_\infty
   =\max(n,|\ell|).
\]
Thus $N$ is precisely the largest absolute index of the generator tuple
in this summand. Unless its marked entry and contracted index form an
opposite pair of magnitude $n$, this tuple lies outside $\mathcal C_{s,1}$:
a larger pair would require another output entry larger than $n$, and
a smaller pair could not be the largest pair. Therefore
\[
 (1+n)|\ell c_0|
 \le B_0\frac{(1+n)|\ell|}{(1+N)^2}\le B_0.
\]
In the exceptional case, exchanging the two marked positions between
the inputs changes $\ell=n$ to $\ell=-n$, or conversely. Both input
arrays have the same values by symmetry, and the two subsets have the
same averaging weight. Their Gardner multipliers have opposite signs,
so the terms cancel. This is a pairing of positions and remains valid
when other indices repeat.

\emph{(ii) Both marked indices lie in the generator.}
The contracted index $\ell$ is independent of $n$. If $|\ell|\ge n$,
the generator tuple is outside $\mathcal C_{s,1}$, and the preceding
estimate applies. If $|\ell|<n$, it lies in $\mathcal C_{s,1}$; for a
nonzero summand, step~1 gives
\[
 |c_0|\le\frac{B_0}{1+n},\qquad |\ell|<N_\chi,
 \qquad (1+n)|\ell c_0|\le N_\chi B_0.
\]
These terms have zero action limit.

\emph{(iii) Both marked indices lie in $P$.}
Again $\ell$ is independent of $n$. Let $\mathbf w$ be the tuple of
unmarked output entries assigned to $P$. If $n>|\ell|$, subtraction of
the first action limit replaces the coefficient of $P$ by
$p_0(n,-n,\mathbf w,\ell)$, whose modulus is at most $A_0/(1+n)$.
If $n\le|\ell|$, then $-\ell$ is the unique largest generator index,
and
\[
 |\ell c_0|\le\frac{B_0}{1+|\ell|}\le\frac{B_0}{1+n}.
\]
Using $|p|\le3A_0$ and $|U|\le2A_0$ gives the same remainder estimate
in this case.

Only case~(iii) contributes to the action limit. For $q\ge3$, its
proportion in the position average is
\[
 \frac{\binom{d-2}{q-3}}{\binom d{q-1}}
 =\frac{(q-1)(q-2)}{d(d-1)},\qquad d=q+s-2.
\]
Comparing the derivative factors $qs$ and $(q-2)s$ yields
\begin{equation}\label{eq:coefficient-first-limit}
 L_1R=\frac{q(q-1)}{d(d-1)}\{U,c_0\}_{\mathrm{nor}}.
\end{equation}
For a nonzero contribution, the prefactor is at most
$q/(q-2)\le3$. Each of the three classes has total averaging weight
at most one. We have therefore proved
\begin{equation}\label{eq:coefficient-first-remainder}
 \sup_{n>|\mathbf v|_\infty}(1+n)
 |R(n,-n,\mathbf v)-L_1R(\mathbf v)|
 \le C N_\chi qs A_0B_0.
\end{equation}

\medskip\noindent\textit{4. Second limit and tangential bracket.}
Repeat the first-action-limit calculation for
$\{U,c_0\}_{\mathrm{nor}}$, replacing $q,p,U$ by $q-2,U,V$.
It uses only
\[
 \sup|U|\le2A_0,\qquad \sup|V|\le A_0,\qquad
 |U(h,-h,\mathbf v)-V(\mathbf v)|\le\frac{A_0}{1+h}
 \quad(h>|\mathbf v|_\infty).
\]
Thus no further action limit is required. Symmetry and momentum
conservation pass to $U$, so the cancellation in step~3(i) still applies.
With $(x)_t=x(x-1)\cdots(x-t+1)$, we obtain
\begin{equation}\label{eq:coefficient-second-limit}
 L_2R=\frac{(q)_4}{(d)_4}\{V,c_0\}_{\mathrm{nor}},
\end{equation}
and
\[
 \sup_{h>|\mathbf v|_\infty}(1+h)
 |L_1R(h,-h,\mathbf v)-L_2R(\mathbf v)|
 \le C N_\chi qs A_0B_0.
\]
For a nonzero contribution, $q\ge5$ and
\[
 \frac{(q)_4}{(d)_4}\le\frac q{q-4}\le5,
 \qquad (q-4)s\le qs.
\]
These estimates, together with the uniform bounds for $R,L_1R,L_2R$,
control all three components of the output and give
\begin{equation}\label{eq:coefficient-decaying-bracket}
 a(\{P,c_0\}_{\mathrm{nor}})
 \le C N_\chi qs\,a(p)b(c).
\end{equation}

For the tangential bracket, there is no contracted normal index.
The same three cases give, respectively, the bounds
$B_0/(1+n)^2$, $B_0/(1+n)$, and the first remainder of $p$.
Repeating the calculation with $U,V$ gives the second remainder.
If $R_{\mathrm{tan}}$ denotes the coefficient of the $a$th tangential
contraction with $c_0$, the position averages give
\begin{equation}\label{eq:coefficient-tangential-limits}
 L_1R_{\mathrm{tan}}=
 \frac{(q)_2}{(q+s)_2}\{U,c_0\}_{\mathrm{tan},a},\qquad
 L_2R_{\mathrm{tan}}=
 \frac{(q)_4}{(q+s)_4}\{V,c_0\}_{\mathrm{tan},a}.
\end{equation}
Both factors are at most one. Removal of normal indices does not change
its derivative factor, so summing over $a\in S$ yields
\[
 a(\{P,c_0\}_{\mathrm{tan}})\le CTa(p)b(c).
\]
The argument also covers normal degree zero.

\medskip\noindent\textit{5. High-frequency mass.}
Fix $K\ge J$. Define $a_K$ exactly as $a$, but replace a weight $1+n$
or $1+h$ by $1$ when its marked index is at most $K$. The components
and action limits are unchanged, and
\begin{equation}\label{eq:coefficient-relaxed-weights}
 a_K(p)\le a(p)\le(1+K)a_K(p),\qquad
 \sup|p|\le3a_K(p),\quad \sup|U|\le2a_K(p),\quad
 \sup|V|\le a_K(p).
\end{equation}
Set
\[
 M_{>K}:=\sum_{\substack{n\in\mathcal N_+\\n>K}}u_nu_{-n}.
\]
For a tuple $\mathbf v$ of length $l$, let
$\mathbf v^{\widehat a\widehat b}$ be the tuple obtained by deleting
positions $a,b$. Define the pair-deletion operator by
\[
 (\mathcal D_K f)(\mathbf v):=
 \sum_{\substack{1\le a,b\le l,\ a\ne b\\v_a>K,\ v_b=-v_a}}
 f(\mathbf v^{\widehat a\widehat b}).
\]
Only the positive entry orients each pair. The sum is over positions,
so repeated indices are counted with their multiplicities.

If $P$ has normal degree $q$, write $d=q+2$. Its product with $M_{>K}$
has ordered coefficient
\begin{equation}\label{eq:coefficient-mass-product}
 (M_{>K}P)(\mathbf j)=\frac{(\mathcal D_Kp)(\mathbf j)}{d(d-1)}.
\end{equation}
Indeed, the ordered coefficient of $u_nu_{-n}$ is $1/2$ at each of
$(n,-n)$ and $(-n,n)$; symmetrizing the product gives precisely this
positive-oriented sum. Deleting a newly inserted largest pair, or a
pair in the remaining tuple, gives
\begin{equation}\label{eq:coefficient-mass-limits}
 L_1(M_{>K}P)=\frac{p+\mathcal D_KU}{d(d-1)},\qquad
 L_2(M_{>K}P)=\frac{2U+\mathcal D_KV}{d(d-1)}.
\end{equation}
At a largest pair $n>K$, deletion of that pair cancels when subtracting
the first limit. Every other deletion leaves
$p(n,-n,\mathbf v')-U(\mathbf v')$, of modulus at most
$a_K(p)/(1+n)$. Similarly, at the second pair $h>K$, the numerator of
the second remainder is
\begin{equation}\label{eq:coefficient-mass-remainder}
 \begin{split}
 p(h,-h,\mathbf w)-U(\mathbf w)
 +\sum_{\substack{a\ne b\\w_a>K,\ w_b=-w_a}}
 \bigl[U(h,-h,\mathbf w^{\widehat a\widehat b})
       -V(\mathbf w^{\widehat a\widehat b})\bigr].
 \end{split}
\end{equation}
Each term is bounded by $a_K(p)/(1+h)$, and the number of terms divided
by $d(d-1)$ is uniformly bounded. Deletion of smaller pairs leaves the
larger marked pairs unchanged, also when indices repeat. For marked
indices at most $K$, use the uniform bounds in
\eqref{eq:coefficient-relaxed-weights}. Consequently,
\begin{equation}\label{eq:coefficient-mass-bound}
 a_K(M_{>K}P)\le C a_K(P),
\end{equation}
uniformly in $q\ge0$, with the negative-degree convention of step~2.

We also need estimates for a momentum-preserving block $F$ whose normal indices have
modulus at most $K$. Suppose its normal degree is $s\le5$ and the
supremum of its coefficients is $f$. In the symmetrized product $PF$,
every marked pair above $K$ lies entirely in $P$. Its action limits
therefore have the factors $(q)_2/(q+s)_2$ and $(q)_4/(q+s)_4$, both
at most one, and its remainders are those of $p$ and $U$.
For marked indices at most $K$, the uniform bounds suffice. Thus
\begin{equation}\label{eq:coefficient-finite-core-product}
 a_K(PF)\le C a_K(P)f.
\end{equation}
In $\{P,F\}_{\mathrm{nor}}$, the contracted index has modulus at most
$K$, so the same argument applies, now with the factors in
\eqref{eq:coefficient-first-limit} and
\eqref{eq:coefficient-second-limit}, bounded by $3$ and $5$.
Together with the tangential formula and
\[
 qs\le W_PW_F,\qquad T\le2JW_PW_F,
\]
this gives
\begin{equation}\label{eq:coefficient-finite-core-bracket}
 a_K(\{P,F\}_G)\le C_S(1+K)W_PW_F a_K(P)f.
\end{equation}
As in step~2, the left side is summed over the normal and tangential
output labels. The constants do not depend on the number of normal
tuples in the support of $F$.

\medskip\noindent\textit{6. The generator's mass term.}
Let $\chi_1$ be the block with coefficient $c_1$, and put $K=N_\chi$.
Step~1 shows that $c_1$ is independent of its largest marked pair and
that all remaining indices have modulus less than $K$. Define a block
$A$ of normal degree $s-2$, with the same angle and Taylor labels as
$\chi$, by its coefficient array
\[
 a_A(\mathbf v):=s(s-1)(L_1c)(\mathbf v).
\]
Set $A=0$ if $s<2$. Then $A$ preserves momentum, is supported on normal
indices of modulus at most $K$, and $\sup|a_A|\le20b(c)$.
For a tuple with exactly one opposite pair above $K$ and all remaining
entries at most $K$, only the deletion of that pair can contribute in
\eqref{eq:coefficient-mass-product}. The resulting coefficient is
$a_A/[s(s-1)]=c_1$. If further entries exceed $K$, all remaining
coefficients of $A$ vanish. Hence
\begin{equation}\label{eq:coefficient-generator-mass-decomposition}
 \chi_1=M_{>K}A+B,
\end{equation}
where $B$ has normal degree $s$, is supported on indices of modulus at
most $K$, and has coefficient supremum at most $b(c)$.
For $s<2$, all three blocks are zero.

Define the total mass in these coordinates and its low-frequency part by
\[
 \begin{aligned}
 M&:=\sum_{j\in S}y_j+\sum_{n\in\mathcal N_+}u_nu_{-n},\\
 M_{\le K}&:=\sum_{j\in S}y_j+
       \sum_{\substack{n\in\mathcal N_+\\n\le K}}u_nu_{-n}
       =M-M_{>K}.
 \end{aligned}
\]
For every momentum-preserving monomial,
\[
 \{e^{ik\cdot\theta}y^m u_{j_1}\cdots u_{j_q},M\}_G
 =i\Bigl(\pi(k)+\sum_a j_a\Bigr)
       e^{ik\cdot\theta}y^m u_{j_1}\cdots u_{j_q}=0.
\]
Thus $\{P,M_{>K}\}_G=-\{P,M_{\le K}\}_G$, and the product rule gives
\begin{equation}\label{eq:coefficient-one-action-bracket}
 \{P,\chi_1\}_G=M_{>K}\{P,A\}_G
       -A\{P,M_{\le K}\}_G+\{P,B\}_G.
\end{equation}
This is an identity of coefficients, each of which is a finite sum
over positions.

By \eqref{eq:coefficient-mass-bound} and
\eqref{eq:coefficient-finite-core-bracket}, the first and last terms
have $a_K$-bounds at most
$C_S(1+K)W_PW_\chi a(P)b(c)$, since $W_A\le W_\chi$ and $W_B=W_\chi$.
The middle bracket multiplies the coefficient indexed by $\mathbf j$ by
\begin{equation}\label{eq:coefficient-finite-momentum}
 i\left(\pi(k_P)+\sum_{a:\,|j_a|\le K}j_a\right).
\end{equation}
Its modulus is at most
$J|k_P|_1+Kq\le C_S(1+K)W_P$.
Insertion or removal of an opposite pair leaves this multiplier
unchanged, so multiplication by it commutes with both action limits
and their extensions. Apply \eqref{eq:coefficient-finite-core-product}
to the subsequent product with $A$, and then restore the original
weights by \eqref{eq:coefficient-relaxed-weights}. Since $K=N_\chi\ge1$,
we obtain
\begin{equation}\label{eq:coefficient-one-action-bound}
 a(\{P,\chi_1\}_G)
 \le C_S K^2W_PW_\chi a(P)b(c).
\end{equation}
Every output has weighted Taylor degree $d_P+d_\chi-2$, where
$d_P=2|m_P|_1+q$ and $d_\chi=2|m_\chi|_1+s$.
In the middle term of \eqref{eq:coefficient-one-action-bracket}, this
follows because $A$ has normal degree $s-2$.

\medskip\noindent\textit{7. Summation.}
Combining steps~3--6 gives, for each pair of input blocks,
\begin{equation}\label{eq:coefficient-block-bound}
 \sum_\lambda a(r_\lambda)
 \le C_S N_\chi^2 W_PW_\chi a(p)b(c).
\end{equation}
Here $r_\lambda$ is the coefficient array of the output block with label
$\lambda$, and $P,\chi$ still denote the fixed input blocks. Recall that
$D_X=|k_X|_1+d_X$ and
$W_X=\max(1,|k_X|_1,|m_X|_\infty,\text{normal degree of }X)$.
Each nonzero output satisfies
\[
 \begin{gathered}
 d_{\mathrm{out}}=d_P+d_\chi-2,\qquad
 D_{\mathrm{out}}\le D_P+D_\chi-2,\qquad
 W_{\mathrm{out}}\le W_P+W_\chi,\\
 N_\chi^2\le C_S(1+\mu)^2(1+D_\chi)^2.
 \end{gathered}
\]
Since $W_X\le1+D_X$ and
$\sup_{x\ge0}x^je^{-\beta x}\le C_j\beta^{-j}$, the assumptions
$0<\delta\le\alpha\le1$ imply
\begin{equation}\label{eq:coefficient-width-loss}
 \begin{split}
 &(W_P+W_\chi)(1+D_\chi)^2
       e^{-\delta D_P-\alpha D_\chi}\\
 &\hspace{1cm}\le C(\delta^{-1}\alpha^{-2}+\alpha^{-3})
       \le C\delta^{-1}\alpha^{-2}.
 \end{split}
\end{equation}
The factors $W_PW_\chi$ are absorbed by the input norms. Moreover,
\[
 r^{d_{\mathrm{out}}}=r^{-2}r^{d_P}r^{d_\chi},\qquad
 e^{\eta D_{\mathrm{out}}}\le e^{\eta D_P}e^{\eta D_\chi}.
\]
Summing \eqref{eq:coefficient-block-bound} with these weights proves
\[
 \|\{P,\chi\}_G\|_{\mathcal H_{\eta,r}}
 \le\frac{C_S(1+\mu)^2}{r^2\alpha^2\delta}
       \|P\|_{\mathcal H_{\eta+\delta,r}}
       \|\chi\|_{\mathcal C_{G,\eta+\alpha,r}}.
\]
To justify the summation, first restrict to finitely many input blocks,
without truncating their normal arrays. Each coefficient contraction is
the finite sum \eqref{eq:coefficient-normal-contraction}. The estimate
then gives absolute convergence of the full block sum in
$\mathcal H_{\eta,r}$. Continuity of coefficient evaluation and of both
action limits, established in step~1, identifies the sum with the
formal bracket and permits passage of both limits through the block sum.

\medskip\noindent\textit{8. Parameter norms and Lie series.}
For two parameters $\xi,\zeta$, bilinearity gives
\[
 \begin{split}
 \{P(\xi),\chi(\xi)\}_G-\{P(\zeta),\chi(\zeta)\}_G
 ={}&\{P(\xi)-P(\zeta),\chi(\xi)\}_G\\
   &+\{P(\zeta),\chi(\xi)-\chi(\zeta)\}_G.
 \end{split}
\]
The preceding estimate bounds its Lipschitz seminorm by
\[
 \frac{C_S(1+\mu)^2}{r^2\alpha^2\delta}
 \left(\|P\|^{\mathrm{lip}}\|\chi\|^{\sup}
       +\|P\|^{\sup}\|\chi\|^{\mathrm{lip}}\right),
\]
with the input spaces as above. Multiplication by $r^2$ and addition
of the supremum estimate prove the asserted total-norm bound.
Bounded Lipschitz families are complete because the coefficient spaces
are Banach.

For the iterated bracket, fix $n\ge1$ and set
\[
 \eta_j=\eta+a-\frac{ja}{n},\qquad j=0,\ldots,n.
\]
At step $j$, apply the bracket estimate at width $\eta_j$, with
$\delta=a/n$ and $\alpha=a$. Indeed,
$\eta_{j-1}=\eta_j+a/n$ and $\eta_j+a\le\eta+2a$, so the generator
is defined on the required strip. Hence
\[
 \|\operatorname{ad}_\chi^jP\|_{\mathcal H_{\eta_j,r}}
 \le Kn\,\|\operatorname{ad}_\chi^{j-1}P\|_{\mathcal H_{\eta_{j-1},r}},
\]
which gives the stated $(Kn)^n$ bound. Finally, $n^n/n!\le e^n$ implies,
for $|t|\le1$ and $eK\le1/2$,
\[
 \sum_{n\ge0}\frac{|t|^n}{n!}
       \|\operatorname{ad}_\chi^nP\|_{\mathcal H_{\eta,r}}
 \le\sum_{n\ge0}(eK)^n\|P\|_{\mathcal H_{\eta+a,r}}
 \le2\|P\|_{\mathcal H_{\eta+a,r}}.
\]
The same argument in total norm completes the proof.
\end{proof}
The dispersive Hamiltonian $H_2=\sum_{n>0}n^2u_nu_{-n}$ lies outside
$\mathcal H_{\eta,r}$, since its ordered coefficient at $(n,-n)$ is
$n^2/2$. We use the homological equation to compute its first bracket;
the preceding estimates then apply to the resulting Hamiltonian.

\section{Homological equations}
\label{sec:max4-homological}

We prove the homological and measure estimates under explicit hypotheses
on the frequencies. These hypotheses are verified for the successive
normal forms in Sections~\ref{sec:max3-retained-bootstrap} and~\ref{sec:loop}.

Constants denoted by $C_{S,X_0,B_0,\tau}$ depend only on the indicated
quantities. In particular, they are uniform in the normal indices, the
parameter set, its interval lengths and the fourth-index cutoff.

Fix a nonempty finite set $S\subset\mathbb N$, let $d=|S|$ and $J=\max S$,
and put $\mathcal N_S=\mathbb Z\setminus(S\cup(-S)\cup\{0\})$. Let
\[
 B=\prod_{j\in S}[a_j,b_j],\qquad
 \ell_* = \min_{j\in S}(b_j-a_j)>0,\qquad E\subset B
\]
with $E$ Borel. The coordinates $\xi_j$ are squared Fourier amplitudes.
The total parameter norm of Section~\ref{sec:coefficient-definitions}
becomes
\begin{equation}\label{eq:max4-total}
 \|P\|^{\mathrm{tot}}_{\mathcal H,\eta,r;E}
 =\sup_{\xi\in E}\|P(\xi)\|_{\mathcal H,\eta,r}
 +r^2\sup_{\substack{\xi,\zeta\in E\\\xi\ne\zeta}}
 \frac{\|P(\xi)-P(\zeta)\|_{\mathcal H,\eta,r}}
      {\|\xi-\zeta\|_1}.
\end{equation}
We use the same convention for generating Hamiltonians.

\subsection{Transversality}

Let \(f:E\to\mathbb R^S\), \(g_n:E\to\mathbb R\), and
\(a:E\to\mathbb R\). Assume
\begin{align}
 \omega_j(\xi)&=j^3+ja(\xi)-6j\xi_j+jf_j(\xi), &&j\in S,\notag\\
 \Omega_n(\xi)&=n^3+na(\xi)+g_n(\xi), &&n>0,\ n\notin S.
 \label{eq:max4-frequencies}
\end{align}
Extend $\Omega_n$ and $g_n$ oddly to negative normal indices; $g_n$
includes the Gardner factor $n$. Assume
$0\le\epsilon\le1/12$, $0\le B_0\le1/5$, and
\begin{align}
 \|f(\xi)-f(\zeta)\|_1
 &\le\epsilon\|\xi-\zeta\|_1, &
 \sup_{n>0,\ n\notin S}|g_n(\xi)-g_n(\zeta)|
 &\le\epsilon\|\xi-\zeta\|_1,\label{eq:max4-differences}\\
 \sup_{\xi,n}|g_n(\xi)|&\le B_0, &
 \sup_{\xi\in E}\bigl(\|\xi\|_1+\|f(\xi)\|_1\bigr)&\le X_0.
 \label{eq:max4-supbounds}
\end{align}
Momentum conservation cancels the common term $a$, so only $f$ and $g$
enter the estimates below.

For $k\in\mathbb Z^S$ set
\[
 \pi(k)=\sum_{j\in S}jk_j,\qquad w_j=jk_j,\qquad K=1+|k|_1.
\]
An admissible normal tuple $\boldsymbol n=(n_1,\ldots,n_q)\in\mathcal N_S^q$
has length $0\le q\le5$ and satisfies
\begin{equation}\label{eq:max4-selection}
 \pi(k)+\sum_{a=1}^q n_a=0,\qquad
 q+\sum_{j\in S}k_j\equiv0\pmod2.
\end{equation}
Its divisor is
\begin{equation}\label{eq:max4-divisor}
 D_{k,\boldsymbol n}
 =\sum_{j\in S}j^3k_j+\sum_a n_a^3
   -6w\cdot\xi+w\cdot f(\xi)+\sum_a g_{n_a}(\xi).
\end{equation}

\begin{lemma}[Transversality on the parameter set]\label{lem:max4-transversality}
For $q\le5$,
\[
 \operatorname{Lip}_{E,\ell^1}D_{k,\boldsymbol n}\le7JK.
\]
If $k\ne0$, then for every $h>0$,
\begin{equation}\label{eq:max4-strip}
 \bigl|\{\xi\in E:|D_{k,\boldsymbol n}(\xi)|<h\}\bigr|
 \le\frac{2h}{5\ell_*}|B|.
\end{equation}
\end{lemma}
\begin{proof}
By \eqref{eq:max4-divisor},
\[
 \begin{aligned}
 |D(\xi)-D(\zeta)|
 &\le\bigl((6+\epsilon)\|w\|_\infty+q\epsilon\bigr)\|\xi-\zeta\|_1\\
 &\le7JK\|\xi-\zeta\|_1.
 \end{aligned}
\]
For $k\ne0$, choose $i$ with $|w_i|=\|w\|_\infty\ge1$. On each fiber
parallel to the $\xi_i$-axis,
\[
 \begin{aligned}
 |D(\xi)-D(\zeta)+6w_i(\xi_i-\zeta_i)|
 &\le(\epsilon\|w\|_\infty+q\epsilon)|\xi_i-\zeta_i|\\
 &\le6\epsilon\|w\|_\infty|\xi_i-\zeta_i|.
 \end{aligned}
\]
Hence every difference quotient has the sign of $-w_i$ and modulus at
least $5\|w\|_\infty$. The portion of $\{|D|<h\}$ on this Borel fiber
has diameter at most $2h/(5\|w\|_\infty)$. Fubini's theorem gives
\[
 \bigl|\{|D|<h\}\cap E\bigr|
 \le\frac{2h}{5\|w\|_\infty(b_i-a_i)}|B|.
\]
This proves \eqref{eq:max4-strip}, without any convexity assumption on $E$.
\end{proof}

\subsection{Measure estimates}

Delete all opposite pairs $(n,-n)$ from the tuple of normal indices.
The remaining tuple, called the \emph{reduced tuple}, is unique up to
permutation. This cancellation leaves the divisor and the selection rules
unchanged and does not increase $\mu_4$. A tuple is integrable if $k=0$
and its reduced tuple is empty.

\begin{lemma}[Dispersion and divisor counting]
\label{lem:max4-dispersion}
Let $\mu\ge0$. There exists a constant
\[
 C_0=C_0(S,X_0,B_0)\ge8
\]
such that every reduced admissible tuple of length at most five with $\mu_4\le\mu$ satisfies
\begin{equation}\label{eq:max4-large}
 N:=\max_a|n_a|>T(k,\mu):=C_0(1+\mu+|k|_1)^2
 \quad\Longrightarrow\quad |D_{k,\boldsymbol n}|\ge\tfrac14N^2
 \quad\hbox{on }E.
\end{equation}
For fixed $k\ne0$, at most
\begin{equation}\label{eq:max4-count}
 C_{S,X_0,B_0}(1+\mu)^8K^8
\end{equation}
distinct divisor functions among all admissible tuples with $q\le5$ and
$\mu_4\le\mu$ can have modulus less than one anywhere on $E$.
For $k=0$, every nonintegrable admissible tuple has $|D|\ge5$ on $E$.
\end{lemma}
\begin{proof}
Write $A=1+\mu+|k|_1$. For a reduced tuple of length at least three, choose its three
largest absolute entries $x,y,z$ and write
\[
 P=x+y+z=-\pi(k)-\sum_{\mathrm{low}}n_a,
 \qquad |P|\le J|k|_1+2\mu\le(J+2)A.
\]
There are at most two low entries, each bounded by $\mu$. The identity
\[
 x^3+y^3+z^3=P^3-3(P-x)(P-y)(P-z)
\]
has three nonzero integer factors: each is a sum of two entries of a reduced tuple. Relabel the entries so that $|x|=N$ and $|y|\ge|z|$. If
$N\ge8(1+|P|)$, then
\[
 |P-x|\ge\frac78N,\qquad
 |y|\ge\frac7{16}N,\qquad |P-y|\ge\frac5{16}N.
\]
The remaining factor has magnitude at least one, so three times the product
has magnitude at least $105N^2/128$.
All other terms of the divisor have total magnitude at most $HA^3$, where
\[
 H=(J+2)^3+2+J^3+6JX_0+5B_0+1.
\]
Indeed the affine correction is at most
$\|w\|_\infty(6\|\xi\|_1+\|f\|_1)\le6J|k|_1X_0$.
Choose
\[
C_0\ge\max\{8,8(J+3),2\sqrt H\}.
\]
Then $N>C_0A^2$ implies both
$N\ge8(1+|P|)$ and $HA^3\le N^2/4$, proving
\eqref{eq:max4-large} for these tuples.

For a reduced tuple of length two, $P=x+y=-\pi(k)$ is a nonzero integer. If
$N=|x|\ge8(1+|P|)$, then $|y|\ge7N/8$ and
\[
 x^3+y^3=P^3-3xyP,\qquad |3xyP|\ge21N^2/8.
\]
The same estimate for the remaining terms proves \eqref{eq:max4-large}.
For a reduced tuple with one entry, momentum conservation gives $N\le J|k|_1$,
so the hypothesis of \eqref{eq:max4-large} cannot hold. We set $N=0$ for
the empty tuple.

For $k=0$, parity restricts the original length to $0,2,4$. A zero-momentum
pair is integrable. A zero-momentum quartet satisfies
\[
 a^3+b^3+c^3+(-a-b-c)^3=-3(a+b)(a+c)(b+c).
\]
This vanishes exactly for a tuple consisting of opposite pairs. Otherwise
it is a nonzero multiple of six, because $n^3\equiv n\pmod6$ and the
momentum is zero. The normal correction has magnitude at most $4B_0$, so
$|D|\ge6-4B_0\ge26/5>5$.

A divisor whose modulus is less than one at some point of $E$ must have a reduced tuple with
all entries bounded by $T=T(k,\mu)$. The first $r-1$ entries of a reduced tuple of
length $r$ determine the last by momentum. Thus the number of possible
functions is at most
\[
 1+\sum_{r=1}^5(2\lceil T\rceil)^{r-1}
 \le C_{S,X_0,B_0}(1+\mu)^8K^8.
\]
Restoring opposite pairs leaves the divisor unchanged.
\end{proof}

\begin{proposition}[Melnikov estimates]\label{prop:max4-exclusion}
Let $\tau>d+8$, $0<\kappa\le1$, and $\mu\ge0$. Define $E_\kappa$ to be
the set of all points of $E$ at which
\begin{equation}\label{eq:max4-melnikov}
 |D_{k,\boldsymbol n}|\ge\kappa K^{-\tau}
\end{equation}
for every nonintegrable admissible tuple with $q\le5$ and $\mu_4\le\mu$.
Then $E_\kappa$ is Borel and
\begin{equation}\label{eq:max4-loss}
 |E\setminus E_\kappa|
 \le C_{S,X_0,B_0,\tau}(1+\mu)^8\frac\kappa{\ell_*}|B|.
\end{equation}
Here $|B|$ is the measure of the original parameter box.
\end{proposition}
\begin{proof}
The set $E_\kappa$ is Borel by continuity of the divisors. No exclusion
is needed for $k=0$. For $k\ne0$, combining
\eqref{eq:max4-strip} and \eqref{eq:max4-count} yields
\[
 |E\setminus E_\kappa|
 \le C_{S,X_0,B_0}(1+\mu)^8\frac\kappa{\ell_*}|B|
       \sum_{k\in\mathbb Z^S\setminus\{0\}}(1+|k|_1)^{8-\tau}.
\]
The series converges for $\tau>d+8$.
\end{proof}

\subsection{Homological estimates}

We use the coefficient and generating spaces of
Section~\ref{sec:verified-coefficients}.

\begin{proposition}[The homological equation]\label{prop:max4-inverse}
Suppose $P$ is real, momentum preserving and even, and its support satisfies
\[
 \Pi_{\mathrm{int}}P=0,\qquad
 2|m|_1+q\le5,\qquad |m|_1\le1,\qquad \mu_4(P)\le\mu.
\]
Assume $\|P\|^{\mathrm{tot}}_{\mathcal H,\eta+\sigma,r;E}<\infty$,
where $\eta\ge0$, $r>0$ and $0<\sigma\le1$. On $E_\kappa$ define
\[
 \chi_{k,m,q}(\boldsymbol n)
 =\frac{p_{k,m,q}(\boldsymbol n)}{iD_{k,\boldsymbol n}}
\]
on the nonintegrable coefficient indices in the support of $P$, and set
all remaining coefficients to zero.
Then $\chi\in\mathcal C_G^{\mathrm{aff}}(S;\mu)$ and
\begin{align}
 \|\chi\|^{\mathrm{tot}}_{\mathcal C_G^{\mathrm{aff}},\eta,r;E_\kappa}
 &\le C_{S,X_0,B_0,\tau}(1+\mu)^6
 \left(\kappa^{-1}\sigma^{-\tau-4}
       +r^2\kappa^{-2}\sigma^{-2\tau-5}\right)
 \|P\|^{\mathrm{tot}}_{\mathcal H,\eta+\sigma,r;E}.
 \label{eq:max4-inverse}
\end{align}
The estimate involves no inverse powers of the tangential amplitudes.
For
\[
 Q=\sum_{j\in S}\frac{\omega_j}{j}y_j
   +\sum_{n>0,\ n\notin S}\frac{\Omega_n}{n}u_nu_{-n},
\]
the generating Hamiltonian solves $\{\chi,Q\}_G=P$ coefficientwise.
The solution is unique with zero integrable projection in the stated
support class.
\end{proposition}
\begin{proof}
The cutoff implies $p_2=0$ and is preserved by the components $p_0,p_1$, as in the
proof of Theorem~\ref{thm:max-coefficient-algebra}.

For every admissible tuple, define $M_{k,\boldsymbol n}=1/(iD_{k,\boldsymbol n})$
when the tuple is nonintegrable and $D_{k,\boldsymbol n}\ne0$, and set
$M_{k,\boldsymbol n}=0$ otherwise. No cutoff is imposed on this definition.
Adding an opposite pair leaves both integrability and the divisor
unchanged. Thus for $N>|v|_\infty$,
\[
 M_{k,(N,-N,v)}=M_{k,v}.
\]
Thus $L_1,L_2,T_1,T_2$ commute with multiplication by $M$, and
\begin{equation}\label{eq:max4-components}
 c_0=p_0M,\qquad c_1=p_1M,\qquad c_2=0.
\end{equation}
These identities also hold for parameter differences. Support is
preserved, and $D_{-k,-\boldsymbol n}=-D_{k,\boldsymbol n}$ gives reality.

We estimate $M$ only on the supports of $p_0,p_1$ and their action
limits. These supports satisfy the cutoff, so
\eqref{eq:max4-melnikov} holds there at every parameter in $E_\kappa$.
Consider such a tuple outside $\mathcal C_{q,1}$, with maximum absolute index $N$. If all
occurrences of $N$ disappear in reduction, the positive and negative
multiplicities are equal. They cannot both equal one, which would place
the tuple in $\mathcal C_{q,1}$. They are therefore both at least two, and the fourth
index is $N\le\mu$. In the remaining case an occurrence of $N$ survives,
and $N$ is the largest absolute index in the reduced tuple. It is then either at most
$T(k,\mu)$ or satisfies the uniform large-index bound
\eqref{eq:max4-large}. In the bounded cases,
\[
 (1+N)^2\le C_{S,X_0,B_0}(1+\mu)^4K^4.
\]
In the large case, $(1+N)^2/|D|$ is bounded by an absolute constant. Thus
\begin{equation}\label{eq:max4-reciprocal}
 (1+N)^2|M|
 \le C_{S,X_0,B_0}(1+\mu)^4\kappa^{-1}K^{\tau+4}
 \qquad (\boldsymbol n\notin\mathcal C_{q,1}).
\end{equation}
The empty tuple, for which $N=0$, belongs to the bounded case.

On these supports, the reciprocal identity at two distinct parameters gives
\[
 \frac{|M(\xi)-M(\zeta)|}{\|\xi-\zeta\|_1}
 \le\frac{7JK}{|D(\xi)D(\zeta)|}.
\]
In the bounded cases use both Melnikov bounds and the bound on $(1+N)^2$.
In the large case both denominators are at least $N^2/4$, so the weighted
left side is at most $CJK/N^2$. In both cases,
\begin{equation}\label{eq:max4-reciprocal-difference}
 (1+N)^2\frac{|M(\xi)-M(\zeta)|}{\|\xi-\zeta\|_1}
 \le C_{S,X_0,B_0}(1+\mu)^4\kappa^{-2}K^{2\tau+5}.
\end{equation}
On $\mathcal C_{q,1}$, the weight of the original zeroth component provides the
required factor $1+N$. For the first component, the weight associated with
the second action limit is unchanged. In these two components the unweighted bounds
$|M|\le\kappa^{-1}K^\tau$ and
$\operatorname{Lip}M\le7J\kappa^{-2}K^{2\tau+1}$ suffice.

Set
\[
 A_k=C_{S,X_0,B_0}(1+\mu)^4\kappa^{-1}K^{\tau+4},\qquad
 B_k=C_{S,X_0,B_0}(1+\mu)^4\kappa^{-2}K^{2\tau+5}.
\]
Equations \eqref{eq:max4-components}--\eqref{eq:max4-reciprocal-difference}
give, block by block,
\begin{align*}
 b(c(\xi))&\le A_k a(p(\xi)),\\
 b(c(\xi)-c(\zeta))
 &\le A_k a(p(\xi)-p(\zeta))
      +B_k a(p(\zeta))\|\xi-\zeta\|_1.
\end{align*}
To sum the blocks, use $\Delta=|k|_1+2|m|_1+q\ge K-1$ and
\[
 K^t e^{-\sigma\Delta}
 \le e\,\sigma^{-t}\sup_{x\ge0}x^te^{-x}
 \le C_t\sigma^{-t}.
\]
Applying this inequality with $t=\tau+4$ and $t=2\tau+5$ and taking the
total parameter norm proves \eqref{eq:max4-inverse}, in fact with
$(1+\mu)^4$ in place of $(1+\mu)^6$.

Finally,
\[
 \{e^{ik\cdot\theta}y^m u_{\boldsymbol n},Q\}_G
 =iD_{k,\boldsymbol n}e^{ik\cdot\theta}y^m u_{\boldsymbol n}.
\]
Thus $\{\chi,Q\}_G=P$. Since the nonintegrable divisors do not vanish,
the solution with zero integrable projection is unique.
\end{proof}

By Theorem~\ref{thm:max-coefficient-algebra}, the core indices of the
first action-limit component are bounded by $|\pi(k)|+2\mu$.

\subsection{Cubic nonresonance}

At a fixed opening, suppose the intervals of the new squared amplitudes have length
comparable to $a^2=r^{2-2\delta}$ and are shorter than the already fixed old
intervals. Choosing
\[
 \kappa_0=c_p a^2r^{\delta/2}=c_pr^{2-1.5\delta}
\]
with a sufficiently small fixed $c_p$ makes the loss in
\eqref{eq:max4-loss} at most the assigned multiple of $r^{\delta/2}|B|$.
The fixed data in this choice include an upper bound for $X_0$ and the
fourth-index cutoff. The two radius factors in the inverse are
\[
 \kappa_0^{-1}\asymp_p r^{-2+1.5\delta},\qquad
 r^2\kappa_0^{-2}\asymp_p r^{-2+3\delta}.
\]
With $\kappa_\nu=2^{-\nu}\kappa_0$ and
$\sigma_\nu=\sigma_0 2^{-\nu-c}$, the factors involving the width loss and the nonresonance threshold grow
respectively as $2^{(\tau+5)\nu}$ and $2^{(2\tau+7)\nu}$.
After fixing all constants at the given stage, choose $r$ sufficiently
small that they are bounded by $r^{-0.1\delta}$. This gives
\[
 \|L_\nu^{-1}\|^{\mathrm{tot}}
 \le C_p2^{(2\tau+7)\nu}r^{-2+1.4\delta}.
\]
This order of choices is used in Proposition~\ref{prop:coherent-production}.
The measure losses are summable, since $\sum_{\nu\ge0}\kappa_\nu=2\kappa_0$.

\begin{corollary}[Cubic nonresonance]\label{cor:max4-fresh-cubic}
Let $E_\infty\subset B$ be Borel, and suppose the limiting frequencies on
this set satisfy \eqref{eq:max4-frequencies}--\eqref{eq:max4-supbounds}
with fixed constants. For $0<\kappa_0\le2$ there exists a Borel set
$E_*\subset E_\infty$ with
\[
 |E_\infty\setminus E_*|
 \le C_{S,X_0,B_0,\tau}\frac{\kappa_0}{\ell_*}|B|,
\]
on which all divisors of degree-three homological equations with real,
momentum-preserving even right-hand sides satisfy the nonresonance
condition with threshold $\kappa_0/2$ for every Fourier mode. For a
right-hand side of finite total norm, the inverse of the homological operator satisfies
\eqref{eq:max4-inverse} with $\mu=0$ and $\kappa=\kappa_0/2$.
\end{corollary}
\begin{proof}
For a cubic source,
\[
 (|m|_1,q)\in\{(0,3),(1,1)\},\qquad
 \mu_4(P)=0,\qquad\Pi_{\mathrm{int}}P=0.
\]
Apply Propositions~\ref{prop:max4-exclusion} and~\ref{prop:max4-inverse}
to the limiting frequencies with $\mu=0$ and $\kappa=\kappa_0/2$.
Restriction to $E_*$ decreases the total norms.
\end{proof}

\subsection{Initial normal form}

\begin{proposition}[The quartic homological equation]
\label{prop:max4-empty-inverse}
Let $S=\varnothing$, $\mu\ge0$, $\eta\ge0$, and $r>0$. Let
$\Omega_n=n^3+na+g_n$ be real and odd, where $a\in\mathbb R$ and
$\sup_{n>0}|g_n|\le B_0\le1/5$. Suppose that $P$ is real, even and
momentum preserving, has degree at most five, and satisfies
$\Pi_{\mathrm{int}}P=0$ and $\mu_4(P)\le\mu$.
Assume also $\|P\|_{\mathcal H,\eta,r}<\infty$. Then $P$ consists only of
nonintegrable quartic terms. Division by
$iD_{\boldsymbol n}=i\sum_{a=1}^4\Omega_{n_a}$ gives the unique
momentum-preserving even solution of degree at most five with zero
integrable projection, and
\begin{equation}\label{eq:max4-empty-bound}
 \|\chi\|_{\mathcal C_G,\eta,r}
 \le16(1+\mu)\|P\|_{\mathcal H,\eta,r}.
\end{equation}
No parameter exclusion is required.
\end{proposition}
\begin{proof}
Evenness permits degrees zero, two and four. Degree zero and
zero-momentum quadratic terms are integrable and hence absent. A
zero-momentum quartet containing one opposite pair contains a second,
so every nonzero quartic coefficient lies outside $\mathcal C_{4,1}$.
Both action limits vanish. Writing $p$ and $c$ for the ordered
quartic coefficients and $N(\boldsymbol n)=\max_a|n_a|$, the norms reduce to
\[
 \|P\|_{\mathcal H,\eta,r}=4e^{4\eta}r^4\sup_{\boldsymbol n}|p(\boldsymbol n)|,
 \qquad
 \|\chi\|_{\mathcal C_G,\eta,r}
 =4e^{4\eta}r^4\sup_{\boldsymbol n}(1+N(\boldsymbol n))^2|c(\boldsymbol n)|.
\]
Fix a nonintegrable quartet $(x,y,z,t)$ with $|x|=N$. Put
$A=x+y$, $B=x+z$, and $C=y+z$. Its cubic dispersion is
\[
 \Delta=-3ABC,
\]
and each factor is a nonzero integer. Since $2x=A+B-C$, one factor
has magnitude at least $2N/3$. The product has magnitude at least that
of its largest factor, so $|\Delta|\ge2N$. It is also a nonzero multiple
of six. The common term cancels and therefore
\begin{equation}\label{eq:max4-empty-relative}
 |D|\ge|\Delta|-4B_0\ge\frac{13}{15}|\Delta|,
 \qquad |D|\ge6-4B_0\ge\frac{26}{5}.
\end{equation}

Let $m$ be the smallest absolute entry, so $m\le\mu$ on the support.
If $m\le N/2$, choose $t$ with $|t|=m$ while keeping $|x|=N$.
Then
\[
 |C|=|x+t|\ge N/2,\qquad |A+B|=|x-t|\ge N/2.
\]
One of $A,B$ has magnitude at least $N/4$ and the other at least one;
hence $|\Delta|\ge3N^2/8$. Using $(1+N)^2\le4N^2$ gives
\[
 \frac{(1+N)^2}{|D|}\le\frac{160}{13}<13.
\]
If $m>N/2$, then $N<2m\le2\mu$, and
\[
 \frac{(1+N)^2}{|D|}
 \le\frac{15}{13}\left(\frac N2+1+\frac1{2N}\right)
 \le\frac{15}{13}\left(\mu+\frac32\right)
 \le2(1+\mu).
\]
These two estimates prove \eqref{eq:max4-empty-bound}. For $\mu<1$ the
source vanishes, since every normal index is a nonzero integer.
Reality, symmetry and uniqueness follow from the corresponding divisor
properties, as in Proposition~\ref{prop:max4-inverse}.

\end{proof}

We apply this inverse to the finite-cutoff quartic projection in
Theorem~\ref{thm:max-initial-cutoff}. The cutoff is fixed before the
normalization radius is chosen.

\section{Hamiltonian flows}
\label{sec:verified-flows}
We construct the generating flows by writing their vector fields as a
bounded part and a scalar transport term. This also identifies the
coefficient Lie series with the corresponding Hamiltonian pullbacks.

Fix \(s\geq2\), \(\sigma=s+4\), a finite \(S\subset\mathbb N\), and
\(H^v=H^v_0(\mathbb T;\mathbb R)\), with
\[
 \|u\|_v^2=\sum_{n>0}\langle n\rangle^{2v}|u_n|^2,
 \qquad u_{-n}=\overline{u_n}.
\]
Write \(\mathcal N=\mathbb Z\setminus(S\cup(-S)\cup\{0\})\),
\(\pi(k)=\sum_{j\in S}jk_j\), and let \(\Pi_S\) denote the Fourier projection onto \(S\cup(-S)\).
Set \(Q=I-\Pi_S\), \(D=\partial_xQ\), and \(D_S=\partial_x\Pi_S\).
Fix \(0<r\leq1\), \(A_0\geq1\), and
\(10r\leq\sqrt{\xi_j}\leq A_0\) for \(j\in S\). For
\(s\leq v\leq\sigma\) and \(0<R\leq3\), define
\[
 \mathcal U_v(R)=\left\{u\in H^v:
 \sum_{j\in S}\bigl||u_j|^2-\xi_j\bigr|<R^2r^2,\quad
 \|Qu\|_v<Rr\right\}.
\]
For $S=\varnothing$, tangential conditions are omitted. Maps on different
Sobolev spaces are compatible if they agree on their common domains.
On these annuli, set
\[
 y_j=|u_j|^2-\xi_j,\qquad e^{i\theta_j}=u_j/|u_j|.
\]
The characters $e^{ik\cdot\theta}$ are independent of the local angles.
Write
\[
 \tau_bu(x)=u(x+b),\qquad
 \Omega(f,g)=\langle\partial_x^{-1}f,g\rangle,\qquad
 \langle f,g\rangle=(2\pi)^{-1}\int_{\mathbb T}fg\,dx .
\]
Here \(\partial_x^{-1}\) has multiplier \(1/(in)\) on \(n\ne0\);
the pairing extends continuously to \(H^{s-1}\times H^s\).

\subsection{Analytic Hamiltonians}
\label{subsec:flows-evaluation}

For a real, momentum-preserving symmetric ordered array
\[
 P=\sum_{k,m,q}\sum_{\mathbf n\in\mathcal N^q}
 p_{k,m,q}(\mathbf n)e^{ik\cdot\theta}y^m
 u_{n_1}\cdots u_{n_q},
 \qquad \pi(k)+\sum_a n_a=0,
\]
set
\[
 \mathfrak M_{\eta,r}(P)=
 \sum_{k,m,q}W e^{\eta(|k|+2|m|+q)}r^{2|m|+q}
 \sup_{\mathbf n}|p_{k,m,q}(\mathbf n)|,
 \quad W=\max(1,|k|,|m|_\infty,q).
\]
We use the coefficient conventions of
Section~\ref{sec:coefficient-definitions}; in particular,
$\mathfrak M_{\eta,r}(P)\leq3\|P\|_{\mathcal H_{\eta,r}}$.
Let
\[
 C_s=4\left(1+\sum_{n\geq1}\langle n\rangle^{-2s}\right)^{1/2},
 \qquad \eta_* =\sigma+3+\log(8C_s).
\]

\begin{proposition}[Evaluation of Hamiltonians]\label{prop:flows-evaluation}
If $\eta\geq\eta_*$ and $\mathfrak M_{\eta,r}(P)<\infty$, the
array defines a real analytic Hamiltonian with a real analytic, mean-zero
$L^2$ gradient $\nabla P:\mathcal U_v(3)\to H^v$, for $s\leq v\leq\sigma$.
The gradient is taken with respect to the normalized spatial average, and
\begin{equation}\label{eq:flows-evaluation-bound}
 \sup_{\mathcal U_v(3)}
 \bigl(|P|+r\|Q\nabla P\|_v
             +r^2\|\Pi_S\nabla P\|_v\bigr)
 \leq C\mathfrak M_{\eta,r}(P),
\end{equation}
where $C=C(S,A_0,s,\sigma)$. The same estimate holds for absolute
majorants. Moreover,
\[
 \|P\|_{C^1(\mathcal U_s(3);\mathbb R)}
 \le Cr^{-2}\mathfrak M_{\eta,r}(P),
\]
where derivatives are measured in $(H^s)^*$.

If $P,R$ satisfy these hypotheses, the formal Poisson bracket
\begin{equation}\label{eq:flows-bracket}
 \{P,R\}_G=
 \sum_{j\in S}j(P_{\theta_j}R_{y_j}-P_{y_j}R_{\theta_j})
 +i\sum_{n\in\mathcal N}nP_{u_n}R_{u_{-n}}
\end{equation}
is absolutely summable at every point of $\mathcal U_s(3)$, with bound
$Cr^{-2}\mathfrak M_{\eta,r}(P)\mathfrak M_{\eta,r}(R)$.
Each coefficient series converges absolutely, and
\[
 \{P,R\}_G(u)=\langle\nabla P(u),\partial_x\nabla R(u)\rangle;
\]
for $u\in\mathcal U_\sigma(3)$ it equals $dP[X_R]$, with $X_R=\partial_x\nabla R$.
Whenever the collected bracket belongs to a coefficient class with width at
least $\eta_*$, its evaluation agrees with this Poisson bracket.

For any array \(\chi\) satisfying the same evaluation hypotheses and
\(H_2=\sum_{n>0}n^2|u_n|^2\), the bracket \(\{H_2,\chi\}_G\) has the
coefficient multiplier
\begin{equation}\label{eq:flows-dispersive-multiplier}
 -i\left(\sum_{j\in S}j^3 k_j+\sum_{a=1}^q n_a^3\right).
\end{equation}
Its series converges absolutely on $\mathcal U_\sigma(3)$ and equals
$dH_2[X_\chi]$.

\end{proposition}

\begin{proof}
Put $d=|S|$ and $c=[100(d+1)(A_0+4)]^{-1}$. At each real point
of $\mathcal U_v(3)$, $9r<|u_j|<A_0+3$. Complexifying the signed
Fourier coordinates independently, the ball of radius $cr^2$ satisfies
\[
 \sum_{j\in S}|y_j|<16r^2,\qquad
 \|u_{\mathcal N,+}\|_v,\ \|u_{\mathcal N,-}\|_v<4r,
 \qquad \max_j|\operatorname{Im}\theta_j|<1,
\]
and $8r<|u_{\pm j}|<A_0+4$. Indeed, the action increment is at most
$d[2(A_0+3)\sqrt2cr^2+2c^2r^4]<r^2$, while the local formula
$\theta_j=(\log u_j-\log u_{-j})/(2i)$ gives angular increments
$O(cr)$ and preserves the chosen logarithmic branches.

Let $f_n=|u_n|$ on the signed normal set and zero elsewhere.  The
Cauchy--Schwarz inequality gives
$\|f\|_{\ell^1},\|f\|_{\ell^2_v}\leq4C_s r=:F$.
Writing $c_{k,m,q}=\sup|p_{k,m,q}|$, the scalar block is bounded by
\begin{equation}\label{eq:flows-scalar-block}
 c_{k,m,q}e^{|k|}(16r^2)^{|m|}F^q.
\end{equation}
A normal derivative removes one of the $q$ ordered positions. By momentum
conservation, the remaining convolution is evaluated at $n-\pi(k)$,
so its majorant is
$q c_{k,m,q}e^{|k|}(16r^2)^{|m|}f^{*(q-1)}(n-\pi(k))$.
For $q\geq2$,
\[
 \left\langle\pi+\sum_{a=1}^{q-1}n_a\right\rangle^v
 \leq q^{v-1}\left(\langle\pi\rangle^v+
                  \sum_{a=1}^{q-1}\langle n_a\rangle^v\right).
\]
Young's inequality, with one weighted factor, gives
\[
 \|f^{*(q-1)}(\cdot-\pi)\|_{\ell^2_v}
 \leq q^v\langle\pi\rangle^vF^{q-1}.
\]
For $q=1$, use the point mass at $\pi$. Hence the normal gradient block
is bounded by
\begin{equation}\label{eq:flows-gradient-block}
 c_{k,m,q}e^{|k|}(16r^2)^{|m|}
 q^{v+1}\langle\pi(k)\rangle^vF^{q-1}.
\end{equation}
If $J=\max(S\cup\{0\})$, then
$\langle\pi(k)\rangle^v\leq(1+J)^\sigma e^{\sigma|k|}$ and
$q^v\leq e^{\sigma q}$.  Use $W\geq q$; the definition of $\eta_*$
absorbs these factors and the powers of $16$ and $4C_s$ in \eqref{eq:flows-scalar-block}--\eqref{eq:flows-gradient-block}.
This proves the scalar and normal gradient estimates of \eqref{eq:flows-evaluation-bound}.

Likewise $|m|\leq dW$ and $|k|\leq W$ give absolute bounds
\[
 \sum_j|P_{y_j}|_{\rm maj}\leq Cr^{-2}\mathfrak M_{\eta,r}(P),
 \qquad
 \sum_j|P_{\theta_j}|_{\rm maj}\leq C\mathfrak M_{\eta,r}(P).
\]
In physical coordinates
\[
 (\nabla P)_j=u_jP_{y_j}+\frac{iP_{\theta_j}}{2u_{-j}},
 \qquad
 (\nabla P)_{-j}=u_{-j}P_{y_j}-\frac{iP_{\theta_j}}{2u_j}.
\]
The amplitude bounds give the remaining part of \eqref{eq:flows-evaluation-bound}.
The resulting series and its gradient converge locally uniformly on the
complexified domain. Analyticity, the $C^1$ estimate and continuity of
evaluation follow.

For the contraction, let $A_{P,n}$ be the sum of the absolute terms of
$P_{u_{-n}}$.  The estimate for the normal gradient gives
$\|A_P\|_{\ell^2_s}\leq Cr^{-1}\mathfrak M_{\eta,r}(P)$.
The expanded normal bracket is bounded by
\[
 \sum_{n\in\mathcal N}|n|A_{P,-n}A_{R,n}
 \leq\|A_P\|_{\ell^2_1}\|A_R\|_{\ell^2_0}
 \leq Cr^{-2}\mathfrak M_{\eta,r}(P)\mathfrak M_{\eta,r}(R).
\]
The tangential contraction satisfies the same bound. Choose a point in
the domain with every $y_j>0$ and every normal coordinate nonzero. For
any fixed output monomial, its absolute majorant subseries is summable;
division by the nonzero monomial value proves absolute convergence of
its coefficient series. Symmetrization and rearrangement are therefore
legitimate. Passing from finite sums to the limit in the coordinate
identity \eqref{eq:flows-bracket} gives
\[
 \{P,R\}_G(u)=\langle\nabla P(u),\partial_x\nabla R(u)\rangle.
\]
For $u\in\mathcal U_\sigma(3)$ this is $dP(u)[X_R(u)]$.

Finally direct differentiation of $H_2$ gives \eqref{eq:flows-dispersive-multiplier}.  The absolute value of the tangential
contraction is bounded by $J^3\sum_j|\chi_{\theta_j}|_{\rm maj}$.
The normal cubic multiplier is controlled by
\[
 \sum_{n\in\mathcal N}|n|^3|u_{-n}|A_{\chi,n}
 \leq\sqrt2\|Qu\|_3\|A_\chi\|_{\ell^2_0}.
\]
This proves the identity on $\mathcal U_\sigma(3)$.
\end{proof}

\subsection{Generating vector fields}
\label{subsec:flows-affine}

\begin{proposition}[The generating vector field]
\label{prop:flows-affine-realization}
Fix \(\eta_0>0\), \(\mu\geq0\), and \(\eta\geq\eta_0\).
Let \(\chi\in\mathcal C_{G,\eta,r}(S;\mu)\) be a real,
momentum-preserving symmetric ordered generating Hamiltonian supported on
\[
 |m|_1\leq1,\qquad 2|m|_1+q\leq5.
\]
Put $\epsilon=\|\chi\|_{\mathcal C_G,\eta,r}<\infty$. Then $\chi$
defines a real analytic Hamiltonian on $\mathcal U_v(3)$, and
\[
 X_\chi=V+aD u,
\]
where \(a:\mathcal U_v(3)\to\mathbb R\) and \(V:\mathcal U_v(3)\to H^v_0\) are compatible real
analytic maps, \(a\) is invariant under full spatial translations, and \(V\)
is equivariant. There is \(C=C(S,A_0,s,\eta_0,\mu)\), independent of
\(r,\xi,\eta\), such that
\begin{align*}
 \|V(u)\|_s&\leq C\epsilon r^{-2},&
 \|QV(u)\|_v&\leq C\epsilon r^{-1},\\
 \|V(u)-V(w)\|_v&\leq C\epsilon r^{-3}\|u-w\|_v,&
 |a(u)|&\leq C\epsilon r^{-2},\\
 |a(u)-a(w)|&\leq C\epsilon r^{-4}\|u-w\|_s,&
 \sum_{j\in S}|2\operatorname{Re}(\overline{u_j}V_j(u))|
 &\leq C\epsilon.
\end{align*}
The two-point estimates hold for every pair \(u,w\in \mathcal U_v(3)\). For every integer
\(h\geq0\), with \(W=V-aD_Su\),
\[
 \sup_{u\in \mathcal U_v(3)}\bigl(
 \|D^hW(u)\|_{\mathcal L^h(H^v;H^v)}+
 \|D^ha(u)\|_{\mathcal L^h(H^v;\mathbb R)}\bigr)
 \leq C_h\epsilon r^{-2-2h}.
\]
The constants \(C_h\) depend only on \(S,A_0,s,\eta_0,\mu,h\);
they and \(C\) are polynomially bounded in \(1+\mu\) with the other
data fixed. The Hamiltonian identity
\(d\chi(u)[h]=\Omega(X_\chi(u),h)\) holds for
\(u\in\mathcal U_s(3)\) and \(h\in H^s\).
\end{proposition}

\begin{proof}
For an ordered block \((k,m,q)\), let \(p\) be its coefficient,
\(\lambda(v)=\lim_{N\to\infty}p(N,-N,v)\), and let \(T\lambda\) be
its extension to tuples with a strict highest opposite pair. Write
\(p^0=p-T\lambda\) and
\[
 \gamma_{k,m,q}=\max\left\{
 \sup_n(1+|n|_\infty)|p^0(n)|,\ \sup_v|\lambda(v)|\right\}.
\]
By the definition of the generating norm,
\[
 \sum_{k,m,q}We^{\eta(|k|+2|m|+q)}r^{2|m|+q}\gamma_{k,m,q}
 \le\epsilon.
\]
The fourth-index cutoff gives $L_2p=0$, as in
Section~\ref{sec:coefficient-definitions}.

Put \(\pi(k)=\sum_{j\in S}jk_j\) and
\[
 K(k)=2+\max(S\cup\{0\})+2\mu+|\pi(k)|.
\]
If $\lambda(v)\ne0$, the core has at most three entries and at most
one has modulus greater than $\mu$. Momentum conservation therefore gives
\[
 |v|_\infty\le|\pi(k)|+2\mu<K(k).
\]

Define the scalar block
\[
 A_{k,m,q}=q(q-1)e^{ik\cdot\theta}y^m
       \sum_{v}\lambda(v)u_{v_1}\cdots u_{v_{q-2}}
 \quad(q\geq2),
\]
and put it equal to zero for \(q<2\). The factor $q(q-1)$ counts the ordered positions of the marked pair.
Above $K(k)$ this pair is unique, also when the core has repeated
entries. Thus
\[
 (T\lambda)_{k,m,q}=M_{>K(k)}A_{k,m,q}+C_{k,m,q},
\]
where \(C_{k,m,q}\) retains only tuples with all indices at most \(K(k)\).
Writing \(M_N=\sum_{n>0,n\notin S}|u_n|^2\) and
\(M_{N,\le K}=\sum_{0<n\le K,\,n\notin S}|u_n|^2\), let
\[
 A=\sum_{k,m,q}A_{k,m,q},\qquad
 B=\sum_{k,m,q}\bigl(C_{k,m,q}-M_{N,\le K(k)}A_{k,m,q}\bigr).
\]
Consequently,
\[
 \chi=\chi_0+M_NA+B,\qquad
 a=A,\qquad V=X_{\chi_0}+M_NX_A+X_B.
\]
Here \(\chi_0\) has coefficients \(p^0\). The scalar blocks of \(A\)
have weighted degree two less than the original block and coefficients at most
\(20\gamma_{k,m,q}\). The blocks of \(B\) have the original degree,
coefficients at most \(21\gamma_{k,m,q}\), and normal indices at most
\(K(k)\). These bounds follow by symmetrizing the product with the finite mass.
Momentum conservation and reality are preserved.

The signed normal sequences satisfy \(\|u\|_{\ell^1}\leq c_s\|u\|_{H^s}\), by
Cauchy--Schwarz. For at most five factors, weighted convolution obeys
\[
 \|b_1*\cdots*b_l(\cdot-\pi)\|_{\ell^2_v}
 \leq C_\sigma\langle\pi\rangle^v
 \sum_a\|b_a\|_{\ell^2_v}\prod_{b\ne a}\|b_b\|_{\ell^1}.
\]
This follows from Young's inequality after distributing the weight over
the factors. With no factors, use the unit mass at $\pi$. For a scalar
block $H$ of weighted degree $d$ and coefficient bound $\gamma$, this yields
\[
 |H|\le C\gamma r^d,\qquad
 \|Q\nabla H\|_v\le C(1+|k|)^\sigma\gamma r^{d-1},\qquad
 |H_y|\le C\gamma r^{d-2}.
\]
For $X_{\chi_0}$, coefficient decay compensates for the Gardner
multiplier; for the finite-core terms, the multiplier is bounded by
$K(k)$. The angular weight absorbs all polynomial factors in $k$,
leaving polynomial dependence on $1+\mu$.

To prove the two-point estimates, set \(z_j=u_j/|u_j|\).
On \(\mathcal U_v(3)\), \(9r<|u_j|\leq A_0+3\), and
\begin{align*}
 |y_j(u)-y_j(w)|&\leq2(A_0+3)|u_j-w_j|,\\
 |z_j(u)-z_j(w)|&\leq C r^{-1}|u_j-w_j|,\\
 |u_{-j}^{-1}-w_{-j}^{-1}|&\leq C r^{-2}|u_j-w_j|.
\end{align*}
Telescoping integer powers on the unit circle gives
\(\left|e^{ik\cdot\theta(u)}-e^{ik\cdot\theta(w)}\right|
 \leq C|k|r^{-1}\|u-w\|_v\).
The tangential component of the vector field is
\[
 (X_\chi)_j=iju_j\chi_{y_j}-\frac{j\chi_{\theta_j}}{2u_{-j}}.
\]
For $H$ affine in $y$, $H_{y_jy_b}=0$. Writing $\Delta F=F(u)-F(w)$,
we use
\[
 \begin{aligned}
 \Delta(u_jH_{y_j})&=(u_j-w_j)H_{y_j}(u)+w_j\Delta H_{y_j},\\
 \Delta(H_{\theta_j}/u_{-j})
   &=\frac{\Delta H_{\theta_j}}{u_{-j}}
     +H_{\theta_j}(w)\frac{w_{-j}-u_{-j}}{u_{-j}w_{-j}}.
 \end{aligned}
\]
Together with the convolution estimate, these identities bound normal
and tangential vector-field differences by
$C\gamma r^{d-3}\|u-w\|_v$, up to polynomial factors in $k$.
No convexity of the annulus is required.

Summing the block bounds gives
\[
 \|QX_{\chi_0}\|_v,\ \|QX_B\|_v\leq C\epsilon/r,
 \qquad
 \|X_{\chi_0}\|_v,\ \|X_B\|_v\leq C\epsilon/r^2,
\]
with Lipschitz constants at most \(C\epsilon/r^3\). Removing the pair in
\(A\) removes two powers of \(r\), giving
\[
 |A|\leq C\epsilon/r^2,\quad
 \operatorname{Lip}_{H^s}A\leq C\epsilon/r^4,\quad
 \|QX_A\|_v\leq C\epsilon/r^3,\quad
 \|X_A\|_v\leq C\epsilon/r^4,\quad
 \operatorname{Lip}_{H^v}X_A\leq C\epsilon/r^5.
\]
Since \(|M_N|\leq9r^2\) and
\(|M_N(u)-M_N(w)|\leq6r\|u-w\|_v\), the product \(M_NX_A\) has precisely
the required size and Lipschitz bounds. The action identity
\(2\operatorname{Re}(\overline{u_j}V_j)=-j\chi_{\theta_j}\)
gives the stated bound on the action variation by the scalar majorant.

Complexify the signed Fourier variables as in
Proposition~\ref{prop:flows-evaluation}, now choosing $c$ so that the
ball of radius $cr^2$ also satisfies
$|\operatorname{Im}\theta|_\infty<\eta_0/4$. The same majorants converge
on this ball, where $|a|+\|W\|_v\le C\epsilon r^{-2}$.
The Cauchy estimates give analyticity and the higher derivative bounds.
The gradient and its Gardner field converge in $H^s$ and $H^{s-1}$,
respectively. By continuity of $\Omega$, the Hamiltonian identity passes
from finite sums to the limit. The symmetries follow from those of the
coefficients.
\end{proof}

\subsection{Flow estimates}
\label{subsec:flows-finite}

\begin{proposition}[Hamiltonian flow]\label{prop:flows-finite}
Suppose that compatible real $C^2$ maps
$a:\mathcal U_v(3)\to\mathbb R$ and
$V:\mathcal U_v(3)\to H^v$ are given for $s\leq v\leq\sigma$.
Assume that $a$ is invariant and $V$ equivariant under full spatial
translation.  Let $M\geq1$, $\epsilon\geq0$, and assume, on these domains,
\begin{align*}
 \|V(u)\|_s&\leq M\epsilon r^{-2},&
 \|QV(u)\|_v&\leq M\epsilon r^{-1},&
 |a(u)|&\leq M\epsilon r^{-2},\\
 \sum_{j\in S}|2\operatorname{Re}(\overline{u_j}V_j(u))|
 &\leq M\epsilon.
\end{align*}
For every pair $u,w\in\mathcal U_v(3)$ assume also
\[
 \|V(u)-V(w)\|_v\leq M\epsilon r^{-3}\|u-w\|_v,
 \qquad
 |a(u)-a(w)|\leq M\epsilon r^{-4}\|u-w\|_s.
\]
Suppose also that a real-valued $C^2$ Hamiltonian $\chi$ on $\mathcal U_s(3)$ satisfies
\[
 d\chi(u)[h]=\Omega(X(u),h),\qquad
 X(u)=V(u)+a(u)Du,\qquad h\in H^s.
\]
If $M\epsilon r^{-3}\leq1/16$, every $u_0\in\mathcal U_\sigma(1)$
has a unique solution of $\dot u=X(u)$ on $[-1,1]$ in
$C(H^\sigma)\cap C^1(H^{\sigma-1})$ remaining in
$\mathcal U_\sigma(3)$; this solution lies in $\mathcal U_\sigma(2)$.
For \(|t|\leq1\) and \(u,w\in\mathcal U_\sigma(1)\), its flow satisfies
\begin{align*}
 \|T^tu-u\|_s&\leq3M\epsilon |t|r^{-2},\\
 e^{-3M\epsilon |t|r^{-3}}\|u-w\|_v
 &\leq\|T^tu-T^tw\|_v
 \leq e^{3M\epsilon |t|r^{-3}}\|u-w\|_v,
 \qquad s\leq v\leq\sigma-1.
\end{align*}
It is continuous into $H^\sigma$ and $C^1$ from $\mathcal U_\sigma(1)$ into
$H^{\sigma-1}$.  Along each such trajectory its derivative extends
to an invertible evolution operator $E(t,t_0)$ on each such $H^v$, with
\[
 \|E(t,t_0)\|,\ \|E(t,t_0)^{-1}\|
 \leq e^{3M\epsilon |t-t_0|r^{-3}},
 \qquad
 \Omega(E(t,0)h,E(t,0)k)=\Omega(h,k).
\]
The operator estimate holds for \(-1\leq t,t_0\leq1\), and the
symplectic identity holds for \(h,k\in H^s\).
\end{proposition}

\begin{proof}
Set $W=V-aD_Su$. Since $D_S$ has finite rank, $W$ is a bounded
$C^2$ vector field with bounded first derivative on each
$\mathcal U_v(3)$. Let $\psi^t$ be its local flow. The subtracted
field only rotates tangential phases, so
\begin{align*}
 \|Qz(t)\|_v&\leq\|Qz(0)\|_v+M\epsilon |t|/r,\\
 \sum_{j\in S}\bigl||z_j(t)|^2-\xi_j\bigr|
 &\leq\sum_{j\in S}\bigl||z_j(0)|^2-\xi_j\bigr|+M\epsilon |t|.
\end{align*}
For initial data in $\mathcal U_v(2)$, these bounds keep the solution
strictly inside $\mathcal U_v(3)$ for $|t|\le1$. Boundedness of $W$
and the usual continuation argument extend $\psi^t$ to that interval.
Put
\[
 \beta^t(u)=\int_0^t a(\psi^q(u))\,dq.
\]
The variational equations give $C^2$ dependence at each Sobolev level;
uniqueness gives equivariance of $\psi^t$ and invariance of $\beta^t$.

The solution of the original equation is
\[
 T^t(u)=\tau_{\beta^t(u)}\psi^t(u).
\]
Strong continuity of translations gives $T^tu\in C(H^\sigma)$, and
differentiation in $H^{\sigma-1}$ yields
\[
 \frac{d}{dt}T^tu=W(T^tu)+a(T^tu)\partial_xT^tu=X(T^tu).
\]
For two solutions write $z=u-w$.  Then
\[
 \dot z=a(u)Dz+V(u)-V(w)+[a(u)-a(w)]Dw.
\]
On $H^v$, $s\leq v\leq\sigma-1$, the last two terms have norm at most
$3M\epsilon r^{-3}\|z\|_v$, because $\|Dw\|_v\leq2r$.
Conjugating by $\exp(D\int a(u))$ and applying Gronwall's inequality
forward and backward gives the metric bounds and uniqueness. The estimate
$\|X(u)\|_s\leq M\epsilon r^{-2}(1+2r)$ gives the displacement estimate.

Differentiability from the $H^\sigma$ domain into $H^{\sigma-1}$
follows from the factorization:
\[
 E(t,0)h=\tau_{\beta^t(u)}
 \bigl(D\psi^t(u)h+D\beta^t(u)[h]\,\partial_x\psi^t(u)\bigr).
\]
Here $(b,z)\mapsto\tau_bz$ is $C^1$ from
$\mathbb R\times H^\sigma$ to $H^{\sigma-1}$. The variational equation is
\[
 \dot h=a(u)Dh+B(t)h,
 \qquad B(t)h=DV(u)h+Da(u)[h]Du,
 \qquad \|B(t)\|_{H^v\to H^v}\leq3M\epsilon r^{-3}.
\]
Conjugation by the translation flow and Gronwall's inequality give the
forward and inverse bounds. Differentiating the factorization identifies
the evolution on directions in $H^\sigma$; density gives the bounded extension.
For $h,k\in H^\sigma$, both variations are continuous in
$H^{\sigma-1}$ and differentiable in $H^{\sigma-2}$.  These spaces
embed in $H^s$.  Differentiating the Hamiltonian identity yields
\[
 \frac{d}{dt}\Omega(h,k)
 =d^2\chi(u)[h,k]-d^2\chi(u)[k,h]=0.
\]
Density and the bounded variational evolution prove the assertion for all
$h,k\in H^s$.
\end{proof}

\subsection{Convergence of the transformations}
\label{subsec:flows-inner}

\begin{proposition}[Convergence of transformations]\label{prop:flows-inner}
Let \((a_\nu,V_\nu,\chi_\nu)\) satisfy all hypotheses of
Proposition~\ref{prop:flows-finite} with the same $S,r,A_0,\xi,M$ and sizes $\epsilon_\nu$.
Set $W_\nu=V_\nu-a_\nu D_Su$.  Suppose
\[
 \mathcal E:=\sum_{\nu\geq0}\epsilon_\nu<\infty,
 \qquad M\mathcal E r^{-3}\leq1/16,
\]
and suppose there are nonnegative \(d_\nu\), independent of \(v\), such
that, on every \(\mathcal U_v(3)\), \(s\leq v\leq\sigma\),
\[
 \|W_\nu\|_{C^2(H^v,H^v)}+
 \|a_\nu\|_{C^2(H^v,\mathbb R)}\leq d_\nu,
 \qquad \sum_{\nu\geq0}d_\nu<\infty.
\]
Here $\|f\|_{C^2}=\sum_{h=0}^2\sup\|D^hf\|$, with $H^v$ input
directions.
For \(T_\nu=T_\nu^1\), the compositions
\(\Phi_N=T_0\circ\cdots\circ T_N\) are defined on
\(\mathcal U_\sigma(1/2)\). Every intermediate point and every trajectory of a constituent flow
remain in \(\mathcal U_\sigma(1)\).
There are compatible \(C^1\) maps
\[
 P:\mathcal U_v(1/2)\to H^v,\qquad
 B:\mathcal U_v(1/2)\to\mathbb R,\qquad s\leq v\leq\sigma,
\]
with \(P\) equivariant and \(B\) invariant, such that
\[
 \Phi_N(u)\longrightarrow\Phi_\infty(u):=\tau_{B(u)}P(u)
\]
pointwise in \(H^\sigma\) and uniformly in \(H^{\sigma-1}\).
Writing the finite compositions as $\Phi_N(u)=\tau_{B_N(u)}P_N(u)$,
we have $P_N\to P$ and $B_N\to B$ uniformly in $C^1$.
The limit is continuous into \(H^\sigma\),
maps the smaller domain into $\mathcal U_\sigma(1)$, and obeys
\begin{align*}
 \|\Phi_\infty(u)-u\|_s&\leq3M\mathcal E r^{-2},\\
 e^{-3M\mathcal E r^{-3}}\|u-w\|_v
 &\leq\|\Phi_\infty(u)-\Phi_\infty(w)\|_v
 \leq e^{3M\mathcal E r^{-3}}\|u-w\|_v,
 \quad s\leq v\leq\sigma-1.
\end{align*}
For $u\in\mathcal U_\sigma(1/2)$ its variational extension is
\[
 E_\infty h=\tau_B\bigl(DPh+DB[h]\,\partial_xP(u)\bigr).
\]
It is invertible on each such $H^v$, and it and its inverse have norm at
most $e^{3M\mathcal E r^{-3}}$.  The finite variational products converge
to it on every fixed vector.  It preserves $\Omega$ on $H^s$ and is
the bounded extension of the derivative from $H^\sigma$ into
$H^{\sigma-1}$.

\end{proposition}

\begin{proof}
The increases in the action deviation and the normal norm under any
finite composition of the flows of the bounded vector fields $W_\nu$ are at most $M\mathcal E$ and $M\mathcal E/r$.
Consequently an input with radius $r/2$ stays within radius $r$, and
an input with radius $r$ stays within radius $2r$.  The same estimates
hold for the actual flows, since translations preserve Fourier moduli.
In particular each constituent Hamiltonian flow starts in the domain of
Proposition~\ref{prop:flows-finite}.

Introduce the augmented bounded vector field and its flow
\[
 Z_\nu(u,b)=(W_\nu(u),a_\nu(u)),\qquad
 F_\nu^t(u,b)=(\psi_\nu^t(u),b+\beta_\nu^t(u)).
\]
For $d=\sum_\nu d_\nu$, the variational equations give
\[
 \|DF_\nu^t\|,\|(DF_\nu^t)^{-1}\|\leq e^{d_\nu},
 \qquad \|D^2F_\nu^t\|\leq d_\nu e^{2d_\nu},
 \qquad \|F_\nu^t-I\|\leq d_\nu,
 \quad |t|\leq1.
\]
Thus $G_N=F_0^1\circ\cdots\circ F_N^1$ has uniformly bounded first
and second derivatives on the larger input annulus: bounds $e^d$ and
$de^{5d}$ suffice.  Write
$G_N(u,b)=(P_N(u),b+B_N(u))$.
We compare consecutive compositions along a trajectory, so that no
convexity of the annulus is needed:
\[
 G_{N+1}(x)-G_N(x)
 =\int_0^1DG_N(F_{N+1}^t x)
       Z_{N+1}(F_{N+1}^t x)\,dt.
\]
The path lies in the larger input annulus. The preceding first and
second derivative bounds imply
\[
 \|G_{N+1}-G_N\|_{C^0}\le e^d d_{N+1},\qquad
 \|D(G_{N+1}-G_N)\|_{C^0}\le C(d)d_{N+1}.
\]
Thus $P_N,B_N$ converge in $C^1$. Compatibility and translation symmetry
pass to the limit.

The derivatives $DP_N$ converge in operator norm and have uniformly
bounded inverses.  For large $N$, the factorization
\[
 DP=DP_N\bigl(I+(DP_N)^{-1}(DP-DP_N)\bigr)
\]
and a geometric series prove that $DP$ is invertible.
For global injectivity of $P$, set
\[
 J=\max(S\cup\{0\}),\qquad
 K_S=(A_0+3)\left(\sum_{j\in S}j^2\langle j\rangle^{2\sigma}\right)^{1/2},
 \qquad L=(1+J+K_S)M r^{-4}\sum_\nu\epsilon_\nu.
\]
On \(\mathcal U_v(3)\), the two-point estimates on the annulus and
\(\|D_Su\|_v\le K_S\), \(\|D_S\|_{H^v\to H^v}\le J\) give
\[
 \|W_\nu(u)-W_\nu(w)\|_v
 \le(1+J+K_S)M\epsilon_\nu r^{-4}\|u-w\|_v.
\]
Gronwall's inequality, applied forward and backward and then composed
over the factors, gives in the limit
\[
 e^{-L}\|u-w\|_v\le\|P(u)-P(w)\|_v
                   \le e^L\|u-w\|_v
 \quad (u,w\in\mathcal U_v(1/2)).
\]
Hence $P$ is a $C^1$ diffeomorphism onto its open image at every
level $s\le v\le\sigma$. No smallness of $L$ is required.

Translation symmetry gives exactly
\[
 \Phi_N=\tau_{B_N}P_N,
 \qquad B_N=B_{N-1}\circ\psi_N+\beta_N.
\]
Strong continuity of translation proves pointwise convergence in
$H^\sigma$.  The estimate
$\|(\tau_b-\tau_c)f\|_{\sigma-1}\leq|b-c|\|f\|_\sigma$
proves uniform convergence one derivative lower.  Summing finite
displacements and multiplying finite two-point estimates proves the stated
bounds, which pass to the limit.

For two pairs $(P_1,b_1)$ and $(P_2,b_2)$, the composed pair is
$(P_1\circ P_2,b_2+b_1\circ P_2)$.  Differentiating equivariance and
invariance gives
\[
 DP_1(w)\partial_xw=\partial_xP_1(w),\qquad
 Db_1(w)[\partial_xw]=0.
\]
These identities give the product rule for the variational operators.
Induction yields the finite forward and inverse bounds and preservation
of $\Omega$.

For $u\in\mathcal U_\sigma(1/2)$, convergence of $DP_N,DB_N$ at level $v$ and of $P_N(u)$
at level $v+1$ gives convergence of the vector inside
\[
 E_Nh=\tau_{B_N}\bigl(DP_Nh+DB_N[h]\,\partial_xP_N(u)\bigr).
\]
Strong continuity of translation on the fixed limiting vector gives the
convergence of $E_Nh$. For $z=\partial_xu$ and $\lambda=DB(u)$, symmetry gives
$DP(u)z=\partial_xP(u)$ and $\lambda[z]=0$.  Hence
\[
 E_\infty=\tau_BDP(I+z\otimes\lambda),\qquad
 E_\infty^{-1}=(I-z\otimes\lambda)(DP)^{-1}\tau_{-B}.
\]
The finite metric bounds give the inverse norm, and convergence on pairs
of vectors preserves \(\Omega\). The factorization identifies this extension with the derivative from
$H^\sigma$ into $H^{\sigma-1}$.
\end{proof}

\subsection{Conjugacy by Lie series}
\label{subsec:flows-pullbacks}

\begin{proposition}[Lie-series representation]\label{prop:flows-lie}
Let \(T^t\) be a flow satisfying Proposition~\ref{prop:flows-finite},
with field \(X\). Let \(F_j:\mathcal U_s(3)\to\mathbb R\), \(j\geq0\),
be \(C^1\) functions with bounded \(C^1\) norms. Assume
\(F_{j+1}=dF_j[X]\) on \(\mathcal U_\sigma(2)\) and
\[
 \sum_{j\geq0}\frac{\|F_j\|_{C^1(\mathcal U_s(3))}}{j!}<\infty,
 \qquad
 \sum_{j\geq0}\frac{\|F_{j+1}\|_{C^0(\mathcal U_s(3))}}{j!}<\infty.
\]
The derivative norm is taken in the dual of \(H^s\). Then, for
\(u\in\mathcal U_\sigma(1)\) and \(|t|\leq1\),
\[
 \sum_{j\geq0}\frac{t^j}{j!}F_j(u)=F_0(T^tu).
\]
If \(H_2=\sum_{n>0}n^2|u_n|^2\), and \(C^1\) functions
\(Q_j:\mathcal U_s(3)\to\mathbb R\) satisfy the same two summability
conditions, \(Q_{j+1}=dQ_j[X]\), and the identity for the first Poisson bracket
\(Q_0=dH_2[X]\) on \(\mathcal U_\sigma(2)\), then
\[
 H_2(T^tu)=H_2(u)+\sum_{j\geq0}\frac{t^{j+1}}{(j+1)!}Q_j(u).
\]
\end{proposition}

\begin{proof}
The assumptions allow termwise differentiation of
$L(t,u)=\sum_{j\ge0}t^jF_j(u)/j!$. Since $\partial_tL=dL[X]$,
\[
 \frac{d}{dq}L(t-q,T^qu)=0,\qquad
 L(t,u)=L(0,T^tu)=F_0(T^tu).
\]
The domain inclusions justify this computation. For $H_2\in C^1(H^s)$,
\[
 \frac{d}{dt}H_2(T^tu)=Q_0(T^tu)
   =\sum_{j\ge0}\frac{t^j}{j!}Q_j(u).
\]
Termwise integration proves the second identity.
\end{proof}

\medskip\noindent\textit{Hamiltonian pullbacks.}
Fix \(\eta\geq\eta_*\), \(0<\delta\leq1/2\), and \(\mu\geq0\).
Let \(P\in\mathcal H_{\eta+\delta,r}(S)\) and
\(\chi\in\mathcal C_{G,\eta+2\delta,r}(S;\mu)\) be real.
Assume that the Hamiltonian vector field \(X_\chi\) has the flow of
Proposition~\ref{prop:flows-finite} on the same annuli. Set
\[
 \operatorname{ad}_\chi P=\{P,\chi\}_G,\qquad
 K=\frac{C_S(1+\mu)^2}{\delta^3r^2}
       \|\chi\|_{\mathcal C_{G,\eta+2\delta,r}},
 \qquad eK\leq1/2,
\]
where \(C_S\) is the coefficient Lie constant of
Theorem~\ref{thm:max-coefficient-algebra}.
The coefficient estimate is
\[
 \|\operatorname{ad}_\chi^nP\|_{\mathcal H_{\eta,r}}
 \leq(Kn)^n\|P\|_{\mathcal H_{\eta+\delta,r}},
 \qquad n\geq1.
\]
By Proposition~\ref{prop:flows-evaluation} and $n!\ge(n/e)^n$, the
scalar Lie series converges in $C^1$, and its time derivative is
controlled by $\sum_{n\ge1}(n+1)(eK)^n<\infty$. The same proposition
gives $F_{n+1}=dF_n[X_\chi]$ for $F_n=\operatorname{ad}_\chi^nP$.
Proposition~\ref{prop:flows-lie} therefore yields
\[
 P(T^tu)=\sum_{n\geq0}\frac{t^n}{n!}\operatorname{ad}_\chi^nP(u)
 \quad (u\in\mathcal U_\sigma(1),\ |t|\leq1).
\]
For the expansion of \(H_2\), assume in addition that
\[
 Q_0=\{H_2,\chi\}_G\in\mathcal H_{\eta+\delta,r}(S).
\]
Propositions~\ref{prop:flows-evaluation} and~\ref{prop:flows-lie} then
apply to the iterated brackets of $Q_0$. In the normal-form steps below,
the homological equation gives $Q_0\in\mathcal H_{\eta+\delta,r}(S)$.

\begin{proposition}[Conjugacy]\label{prop:flows-conjugacy}
Assume the hypotheses of Proposition~\ref{prop:flows-inner}. Let
\(H:\mathcal U_s(3)\to\mathbb R\) and
\(K:\mathcal U_s(1/2)\to\mathbb R\) be \(C^1\).
\begin{enumerate}[label=(\roman*)]
\item If functions \(K_N\) on \(\mathcal U_\sigma(1/2)\)
satisfy \(K_N=H\circ\Phi_N\) and \(K_N(u)\to K(u)\) for each $u\in\mathcal U_\sigma(1/2)$, then \(K=H\circ\Phi_\infty\) there.
\item Let \(T\) be one of the finite compositions or \(\Phi_\infty\),
and suppose \(K=H\circ T\) on \(\mathcal U_\sigma(1/2)\).
Assume there are continuous fields
\[
 X_H:\mathcal U_\sigma(3)\to H^{s+1},\qquad
 X_K:\mathcal U_\sigma(1/2)\to H^{s+1}
\]
satisfying \(dH(u)[h]=\Omega(X_H(u),h)\) and
\(dK(w)[h]=\Omega(X_K(w),h)\) on the respective domains for
every \(h\in H^s\). Then
\[
 X_H(Tu)=E_T(u)X_K(u)\quad\hbox{in }H^{s+1}.
\]
For any interval \(I\) and curve
\(u\in C(I;\mathcal U_\sigma(1/2))\cap C^1(I;H^{s+1})\)
solving \(\dot u=X_K(u)\), its image belongs to
\(C(I;H^\sigma)\cap C^1(I;H^{s+1})\) and solves
\((Tu)'=X_H(Tu)\).
\end{enumerate}
\end{proposition}

\begin{proof}
Assertion (i) follows by continuity of $H$. For (ii), differentiation
in a direction $h\in H^\sigma$ gives
\[
 dK(u)[h]=\Omega(X_H(Tu),E_T(u)h)
 =\Omega(E_T(u)^{-1}X_H(Tu),h).
\]
Density and nondegeneracy of \(\Omega\) give the identity;
boundedness of \(E_T\) and its inverse on \(H^{s+1}\) shows that the identity holds in
\(H^{s+1}\).

For the trajectory, write \(T(u)=\tau_{b(u)}P(u)\).
The maps \(P,b\) are \(C^1\) on \(H^{s+1}\), so the translation product
rule first gives
\[
 \frac{d}{dt}T(u(t))=E_T(u(t))\dot u(t)=X_H(T(u(t)))
 \quad\hbox{in }H^s.
\]
Since the right-hand side is continuous in $H^{s+1}$, integration gives
$Tu\in C^1(I;H^{s+1})$.
\end{proof}

For a Hamiltonian \(H_2+P\) whose real \(L^2\) gradient
\(\nabla P\) is continuous from the annulus in $H^\sigma$ into \(H^\sigma\),
the vector field is
\[
 -\partial_x^3u+\partial_x\nabla P(u)\in H^{\sigma-3}=H^{s+1}.
\]
For the arrays of Proposition~\ref{prop:flows-evaluation}, this verifies
the vector-field hypothesis of Proposition~\ref{prop:flows-conjugacy}.

\subsection{Uniform estimates}
\label{subsec:flows-scope}

Fix $S,A_0,\mu$ and a positive lower bound $\eta_0$ for the analyticity
widths. Consider real generating Hamiltonians that are affine in the
actions and have weighted degree at most five. The constants in the
vector-field estimates
are uniform in the radius, action centers and normal Fourier indices.
Propositions~\ref{prop:flows-affine-realization}--\ref{prop:flows-inner}
therefore apply with
\[
 \|D^h(V_\nu-a_\nu D_Su)\|+
 \|D^h a_\nu\|\leq C_h\epsilon_\nu r^{-2-2h},
 \qquad h=0,1,2,
\]
and
\[
 M r^{-3}\sum_\nu\epsilon_\nu\leq1/16.
\]
At fixed $r$, the sum of the $C^2$ norms is bounded by
$Cr^{-6}\sum_\nu\epsilon_\nu$; the displacement and metric bounds
retain the powers $r^{-2}$ and $r^{-3}$.

For widths at least $\eta_*$, coefficient convergence implies $C^1$
scalar convergence by Proposition~\ref{prop:flows-evaluation}. Together
with the Lie smallness condition and
$\{H_2,\chi\}_G\in\mathcal H$, this justifies the pullback of the Hamiltonians and conjugacy of
the PDEs by Propositions~\ref{prop:flows-lie}--\ref{prop:flows-conjugacy}.
Parameter dependence is treated below; the required
uniform data are supplied by Proposition~\ref{prop:coherent-production}.

\subsection{Parameter dependence}
\label{subsec:flows-parameters}
Fix a finite tangential set $A$ and $0<R\le1$. In the preceding notation, set $Q_A=\operatorname{Id}-\Pi_A$,
$D_A=\partial_xQ_A$ and $D_A^{\rm tan}=\partial_x\Pi_A$.
For $\eta>0$ and $\mu\ge0$, consider real symmetric arrays preserving
momentum, with $|m|_1\le1$, $2|m|_1+q\le5$, $\mu_4\le\mu$, and first
action limit $\lambda$. With $p-T\lambda$ as in
Proposition~\ref{prop:flows-affine-realization} and
$N(\boldsymbol n)=\max|n_i|$, define
\begin{align*}
 \gamma_{k,m,q}&=\max\left\{
  \sup_{\boldsymbol n}(1+N(\boldsymbol n))|p_{k,m,q}-T\lambda_{k,m,q}|,
  \sup_{\boldsymbol v}|\lambda_{k,m,q}(\boldsymbol v)|\right\},\\
 \physnorm{\chi}_{A,\eta,R}
 &=\sum_{k,m,q}\max(1,|k|_1,|m|_\infty,q)
 e^{\eta(|k|_1+2|m|_1+q)}R^{2|m|_1+q}\gamma_{k,m,q}.
\end{align*}
For a parameter set $E$ with its $\ell^1$ metric, define
\[
 \physnorm{\chi}_{E,A,\eta,R}^{\mathrm{tot}}
 =\sup_E\physnorm{\chi_\xi}_{A,\eta,R}
    +R^2\operatorname{Lip}_E(\chi;\physnorm{\cdot}_{A,\eta,R}).
\]
This norm is bounded by the affine generating norm.

Suppose $10R\le\sqrt{\xi_j}\le A_0$, with $A_0\ge1$. For
\(0<b\le3\) and \(s\le v\le\sigma\), set
\[
 U_{\xi,v}(b)=\left\{u\in H^v:
 \sum_{j\in A}\bigl||u_j|^2-\xi_j\bigr|<b^2R^2,
 \quad\physnorm{Q_Au}_v<bR\right\}.
\]
\begin{lemma}[Parameter estimates]\label{lem:physical-parameter}
Let $E\subset\mathbb R^A$ be a nonempty Borel set, and let
$\chi_\xi$ be a family as above with
\[
 \epsilon=\sup_E\physnorm{\chi_\xi}_{A,\eta,R}<\infty,
 \qquad \Lambda=\operatorname{Lip}_E(\chi;\physnorm{\cdot}_{A,\eta,R})<\infty,
 \qquad \eta\ge\eta_0>0.
\]
Use the splitting $X_{\chi_\xi}=V_\xi+a_\xi D_Au$ from
Proposition~\ref{prop:flows-affine-realization}. There is $C$,
depending only on $A,A_0,s,\eta_0,\mu$, such that, for
$u\in U_{\xi,v}(3)$ and $w\in U_{\xi',v}(3)$,
\begin{align}
 \physnorm{V_\xi(u)-V_{\xi'}(w)}_v
 &\le C\epsilon R^{-3}\physnorm{u-w}_v
    +C(\Lambda R^{-2}+\epsilon R^{-3})\physnorm{\xi-\xi'}_1,\label{eq:physical-V}\\
 |a_\xi(u)-a_{\xi'}(w)|
 &\le C\epsilon R^{-4}\physnorm{u-w}_s
    +C(\Lambda R^{-2}+\epsilon R^{-4})\physnorm{\xi-\xi'}_1.\label{eq:physical-a}
\end{align}
If $C\epsilon/R^3$ is sufficiently small, the Hamiltonian flows exist
on $U_{\xi,\sigma}(1)$ for $|t|\le1$. For
$u\in U_{\xi,\sigma}(1)$, $w\in U_{\xi',\sigma}(1)$ and
$s\le v\le\sigma-1$,
\begin{equation}\label{eq:physical-flowparameter}
 \physnorm{T_\xi^tu-T_{\xi'}^tw}_v
 \le e^{C\epsilon|t|/R^3}
 \left(\physnorm{u-w}_v+C|t|(\Lambda/R^2+\epsilon/R^3)\physnorm{\xi-\xi'}_1\right).
\end{equation}
The map $(\xi,t,u)\mapsto T_\xi^tu$ is jointly continuous into
$H^\sigma$ on its parameter-dependent domain, endowed with the
relative product topology inherited from $E\times[-1,1]\times H^\sigma$.
In particular, it is Borel measurable.
\end{lemma}

\begin{proof}
For $j\in A$, we have $|u_j|,|w_j|>9R$ and
\begin{align*}
 y_\xi(u)-y_{\xi'}(w)
 &=|u_A|^2-|w_A|^2-(\xi-\xi'),\\
 |u_j/|u_j|-w_j/|w_j||&\le
       2|u_j-w_j|/\min(|u_j|,|w_j|),\\
 |u_{-j}^{-1}-w_{-j}^{-1}|&\le
       |u_j-w_j|/(|u_j||w_j|).
\end{align*}
Proposition~\ref{prop:flows-affine-realization} gives the linear splitting
$\chi=\chi_0+M_N A_{\rm ext}+B$ and
$V=X_{\chi_0}+M_NX_{A_{\rm ext}}+X_B$, $a=A_{\rm ext}$.
For an input block of weighted degree $d$, the terms of $\chi_0,B$ have
degree $d$ and those of $A_{\rm ext}$ have degree $d-2$; all are affine in $y$.
An action difference costs $R^{-2}$, giving a normal vector-field bound
of order $R^{d-3}$. For the tangential component,
\[
 (X_H)_j=ij u_jH_{y_j}-\frac{jH_{\theta_j}}{2u_{-j}},
\]
the first term has no action variation because $H_{yy}=0$, and the second
has the same $R^{d-3}$ bound. For $A_{\rm ext}$ and $X_{A_{\rm ext}}$ the
powers are $R^{d-4}$ and $R^{d-5}$, respectively; the factor $M_N=O(R^2)$
restores $R^{d-3}$. Coefficient differences contribute
$\Lambda R^{-2}\physnorm{\xi-\xi'}_1$.
Summation with the analytic weights, as in
Proposition~\ref{prop:flows-affine-realization}, gives
\eqref{eq:physical-V}--\eqref{eq:physical-a} on the respective annuli.

Proposition~\ref{prop:flows-finite} applies at each parameter. Along the
two flow trajectories the difference $h=u-w$ solves
\[
 h'=a_\xi(u)D_Ah+[V_\xi(u)-V_{\xi'}(w)]
       +[a_\xi(u)-a_{\xi'}(w)]D_Aw.
\]
The first term generates a unitary partial translation and
$\physnorm{D_Aw}_v\le2R$ for $v\le\sigma-1$. Thus variation of constants,
\eqref{eq:physical-V}--\eqref{eq:physical-a}, and Gronwall's inequality
yield \eqref{eq:physical-flowparameter}.

To prove strong continuity, write the flow as a translation composed with the flow of $W_\xi$, where
\[
 W_\xi=V_\xi-a_\xi D_A^{\rm tan}u,
 \quad (\psi_\xi^t)'=W_\xi(\psi_\xi^t),
 \quad \beta_\xi^t=\int_0^t a_\xi(\psi_\xi^q)\,dq.
\]
Since $D_A^{\rm tan}$ is bounded, \eqref{eq:physical-V}--\eqref{eq:physical-a}
give, for $s\le v\le\sigma$,
\[
 \physnorm{W_\xi(u)-W_{\xi'}(w)}_v+|a_\xi(u)-a_{\xi'}(w)|
 \le \frac{C\epsilon}{R^4}\physnorm{u-w}_v
 +C\left(\frac{\Lambda}{R^2}+\frac{\epsilon}{R^4}\right)
       \physnorm{\xi-\xi'}_1.
\]
The integral equations give joint continuity of $(\psi,\beta)$ in
$H^\sigma\times\mathbb R$. The assertion follows from
$T_\xi^tu=\tau_{\beta_\xi^t(u)}\psi_\xi^t(u)$ and strong continuity
of translations.
\end{proof}

\begin{theorem}[Parameter estimates for the limit]\label{thm:physical-inner}
Let $\chi_{\nu,\xi}$, $\nu\ge0$, be families as above at radius $R$
on a common Borel set $E$. Their positive widths and finite cutoffs may
vary with $\nu$. Denote their supremum norms and global Lipschitz
seminorms by $\epsilon_\nu$ and $\Lambda_\nu$. Choose
$c_\nu\ge1$ larger than the constants in Lemma~\ref{lem:physical-parameter},
the vector-field estimates, and all the realization constants $C_h$
with $h=0,1,2$. Suppose
\[
 E_0=\sum_\nu c_\nu\epsilon_\nu<\infty,
 \qquad E_1=\sum_\nu c_\nu\Lambda_\nu<\infty,
 \qquad E_0/R^3\le c_* 
\]
for a sufficiently small numerical $c_*>0$. Then the inner compositions
$\Phi_{\xi,N}=T_{0,\xi}\circ\cdots\circ T_{N,\xi}$ on
$U_{\xi,\sigma}(1/2)$ have all the conclusions of Proposition~\ref{prop:flows-inner}. Their limit is jointly continuous into $H^\sigma$ and, for inputs $u,w$ in their respective $H^\sigma$ domains,
\[
 \physnorm{\Phi_{\xi,\infty}u-\Phi_{\xi',\infty}w}_v
 \le e^{CE_0/R^3}\left(\physnorm{u-w}_v+
 C(E_1/R^2+E_0/R^3)\physnorm{\xi-\xi'}_1\right),
 \quad s\le v\le\sigma-1.
\]
For qualitative joint continuity it suffices that each $\Lambda_\nu$ is finite and the uniform $E_0$ bound holds.
\end{theorem}

\begin{proof}
The pairs $(W_{\nu,\xi},a_{\nu,\xi})$ have summable $C^2$ norms, bounded by
$3R^{-6}c_\nu\epsilon_\nu$. Proposition~\ref{prop:flows-inner}
therefore yields $C^1$ convergence of the bounded flow components
$P_N$ and the associated translation parameters $B_N$ with respect to
the state variables, uniformly in $\xi$.
Lemma~\ref{lem:physical-parameter} gives joint continuity of the finite
factors, hence of their limits. The translation factorization gives
continuity into $H^\sigma$.

Iterating \eqref{eq:physical-flowparameter} from the innermost factor
and letting $N\to\infty$ proves the stated estimate. The domain
inclusions are those of Proposition~\ref{prop:flows-inner}; we omit
the repeated composition argument. Qualitative continuity requires only
the uniform convergence and continuity of each finite factor.
\end{proof}

\begin{remark}
If the analyticity widths have a fixed positive lower bound and the
cutoffs have a common finite upper bound, one may take $c_\nu=C$.
Geometric growth $c_\nu\le C_1 2^{c\nu}$ is allowed by the weighted
generator estimate of Proposition~\ref{prop:max3-newton}.
\end{remark}

\section{KAM iteration}
\label{sec:verified-newton}
For a fixed tangential set, we first eliminate the nonintegrable jet of
weighted degree at most two by a triangular KAM iteration and then
remove the cubic terms by one normal-form step, following
\cite{BGR25,BP26}. We keep track of the integrable correction produced
by both operations and realize the transformations by the Hamiltonian
flows of Section~\ref{sec:verified-flows}. The frequency hypotheses for
the successive applications are verified in
Section~\ref{sec:max3-retained-bootstrap}.

\subsection{Quadratic jet}
\label{sec:max3-jet-convention}

For a jet $J=J_0+J_1+J_2$ of weighted degree at most two, set
$\epsilon=\sum_{d=0}^2\|J_d\|^{\mathrm{tot}}$.
By Theorem~\ref{thm:coefficient-initial},
\[
 \sum_{d=0}^2\|(1-\Pi_{\mathrm{int}})\Pi_d E\|^{\mathrm{tot}}
 \leq3\|E\|^{\mathrm{tot}}.
\]
Accordingly, with $B_\nu,T_\nu$ and $t_0,\Lambda$ as in
Proposition~\ref{prop:max3-newton}, the error majorant is
\[
 C_0=9B_0T_0(1+t_0)^5,
 \qquad \epsilon_{\nu+1}=C_0\Lambda^\nu\epsilon_\nu^2.
\]
The factor $9$ accounts for the three projected jet norms.

\subsection{Homological system}
\label{sec:max3-newton}

Fix a finite tangential set $S\subset\mathbb N$. For $n>0$, $n\notin S$,
write $z_n=u_n$, $\bar z_n=u_{-n}$ and $\rho_n=z_n\bar z_n$.
The variables $z_n,\bar z_n$ are independent in the formal calculus.
With the degree convention of Section~\ref{sec:coefficient-definitions},
the Poisson bracket reads
\[
 \{U,V\}=\sum_{j\in S}j(U_{\theta_j}V_{y_j}-U_{y_j}V_{\theta_j})
 +i\sum_{n>0,\,n\notin S}n
 (U_{z_n}V_{\bar z_n}-U_{\bar z_n}V_{z_n}).
\]
Write $I=\Pi_{\mathrm{int}}$ for the integrable projection
\eqref{eq:coefficient-integrable-projection}. Integrable terms of degree at most two are diagonal,
and there are no integrable terms of odd degree.
Let $B\subset\mathbb R^M$, $M\geq1$, be a bounded rectangular
parameter box of positive volume and let $\mu$ be its normalized Lebesgue
measure. Let $E_0\subset B$ be an initial Borel parameter set of positive measure.
Fix $r>0$. For a family on a Borel set $E$, the total parameter norm is
\[
 \|U\|_{\eta;E}
 =\sup_{x\in E}\|U(x)\|_{\mathcal H_\eta}
 +r^2\sup_{\substack{x,x'\in E\\x\ne x'}}
       \frac{\|U(x)-U(x')\|_{\mathcal H_\eta}}{\|x-x'\|_1}.
\]
The same convention applies to $\mathcal G_\eta$, including the empty
set and singleton conventions of Section~\ref{sec:coefficient-definitions}.

\begin{proposition}[Abstract KAM iteration]
\label{prop:max3-newton}
Fix $\epsilon_0>0$ and $K_0\geq0$.
Assume that the coefficient and generating spaces are Banach in both
the coefficient and total parameter norms. Restrictions to smaller
widths or parameter sets, degree and integrable projections, their
complements and commuting compositions are contractions preserving the
prescribed closed subspaces. Assume also that coefficient evaluation is
continuous, $\|F\|_{\mathcal H}\le\|F\|_{\mathcal G}$, and the brackets
and Lie series below converge absolutely coefficientwise.

Choose
$\eta_{\nu,0}>\cdots>\eta_{\nu,8}=\eta_{\nu+1,0}>\eta_\infty>0$.
Suppose that $B_\nu,T_\nu>0$ satisfy the following estimates uniformly
on Borel subsets of $E_0$:
\begin{enumerate}
\item For $i=1,3,5,6$ and generators of degree at most two,
\[
 \|\{U,F\}\|_{\eta_{\nu,i+1}}
 \leq B_\nu\|U\|_{\eta_{\nu,i}}\|F\|_{\mathcal G,\eta_{\nu,i}}.
\]
If $B_\nu\|F\|_{\mathcal G,\eta_{\nu,5}}\leq1/16$, then
$e^{t\operatorname{ad}_F}$ has norm at most $2$ from $\eta_{\nu,7}$ to
$\eta_{\nu,8}$ for $|t|\leq1$ and satisfies Taylor's formula with a
second-order integral remainder. The formula remains valid for an
unbounded diagonal plus a bounded perturbation when its first bracket
with $F$ is bounded.
\item Write
\[
 N_0=e_0+\sum_{j\in S}\frac{\omega_{j,0}}{j}y_j
              +\sum_{n>0,\,n\notin S}\frac{\Omega_{n,0}}{n}\rho_n,
\]
with Borel coefficient functions. There are numbers $\ell_\nu\geq0$ with
$\sum_\nu\ell_\nu<\mu(E_0)$ such that every diagonal
$N=N_0+\Delta$ on a Borel $E'\subset E_0$, with integrable
$\deg\Delta\leq2$ and $\|\Delta\|_{\eta_{\nu,0};E'}\leq4\epsilon_0$,
admits a Borel $E''\subset E'$ satisfying
$\mu(E'\setminus E'')\leq\ell_\nu$.
On $E''$, for $i=0,2,4$, the homological equation for every
nonintegrable homogeneous source $U$ of degree zero, one, or two admits
a solution depending linearly on the source, with
\[
 \{N,F\}=-U,\qquad
 \|F\|_{\mathcal G,\eta_{\nu,i+1};E''}
       \leq T_\nu\|U\|_{\eta_{\nu,i};E''}.
\]
The set $E''$ is independent of $U$ and $i$. The inverse acts on sources
defined on $E''$ and preserves support, degree, reality, momentum and evenness.
For sources on $E'$, restriction gives the same estimate with the
$E'$-norm on the right.
\end{enumerate}

Let $H_0=N_0+J_0+P_0$, with $I J_0=0$, $\deg J_0\leq2$,
$\Pi_{\geq3}P_0=P_0$, and
\[
 \sum_{k=0}^2\|\Pi_kJ_0\|\leq\epsilon_0,\qquad
 \|P_0\|\leq K_0,\qquad \|\Pi_3P_0\|\leq\epsilon_0.
\]
Suppose $B_\nu\leq B_0 2^{a\nu}$, $T_\nu\leq T_0 2^{b\nu}$,
where $a,b\geq0$. For $c\geq0$ choose $d\geq a+b+c$ and put
\[
 D_0=B_0T_0,\quad K_*=2(K_0+\epsilon_0),\quad
 t_0=D_0(K_*+\epsilon_0),\quad
 C_0=9D_0(1+t_0)^5,\quad\Lambda=2^{6d}.
\]
If $t_0\leq1/64$ and $x_0=C_0\Lambda\epsilon_0\leq1/16$,
there are coefficient Hamiltonians
$H_\nu=N_\nu+J_\nu+P_\nu$ on decreasing Borel sets $E_\nu$, with
\[
 \mu(E_\infty)\geq\mu(E_0)-\sum_\nu\ell_\nu>0,
 \qquad E_\infty=\bigcap_\nu E_\nu.
\]
On $E_\infty$ the sequence converges at width $\eta_\infty$ to
$N_\infty+P_\infty$, where
\[
 \Pi_{\leq2}P_\infty=0,\quad
 \|N_\infty-N_0\|\leq3\epsilon_0,\quad
 \|P_\infty\|\leq K_*,\quad
 \|\Pi_3P_\infty\|\leq\tfrac32\epsilon_0.
\]
The generators satisfy
\[
 \sum_\nu2^{c\nu}\|F_\nu\|_{\mathcal G}
 \leq2T_0(1+t_0)^2\epsilon_0.
\]
Here the norm of $F_\nu$ is taken on $E_{\nu+1}$ at width
$\eta_{\nu,5}$, and limiting norms on $E_\infty$ at width
$\eta_\infty$. Moreover,
\[
 I\Pi_4(P_\infty-P_0)=\sum_\nu I\Pi_4(Q_\nu+\mathcal E_\nu),
 \qquad
 \left\|I\Pi_4(P_\infty-P_0)\right\|
 \leq\left(2t_0(1+t_0)^2+\tfrac18\right)\epsilon_0,
\]
where $Q_\nu,\mathcal E_\nu$ are defined in the proof.
\end{proposition}

\begin{proof}
At step $\nu$, choose $E_{\nu+1}$ by the inverse hypothesis and restrict
all coefficients to this set. Suppressing $\nu$, put
$P_k=\Pi_kP$, $J_k=\Pi_kJ$ and $LF=\{N,F\}$.
The bracket lowers weighted degree by two, and $F_0$ depends only on
$\theta$. Hence
\[
 \Pi_{\leq2}\{P,F_0+F_1+F_2\}
 =\Pi_1\{P_3,F_0\}
  +\Pi_2\bigl(\{P_3,F_1\}+\{P_4,F_0\}\bigr).
\]
Define successively
\begin{align*}
 S_0&=J_0,& Z_0&=IS_0,& LF_0&=-(S_0-Z_0),\\
 S_1&=J_1+\Pi_1\{P_3,F_0\},& Z_1&=IS_1,& LF_1&=-(S_1-Z_1),\\
 S_2&=J_2+\Pi_2\bigl(\{P_3,F_1\}+\{P_4,F_0\}\bigr),
 & Z_2&=IS_2,& LF_2&=-(S_2-Z_2).
\end{align*}
The source and solution widths are $\eta_{\nu,2k}$ and
$\eta_{\nu,2k+1}$. Each $F_k$ satisfies $|m|_1\le1$, $q\le2$,
and $Z_1=0$ since there are no integrable degree-one monomials.
Put $F=\sum F_k$, $Z=\sum Z_k$, and
$C=\Pi_{\leq2}\{P,F\}$. Then
\[
 \{N,F\}=-J-C+Z.
\]
With $\epsilon=\sum_k\|J_k\|$, $K=\|P\|$ and $t_k=BT\|P_k\|$,
successive substitution gives
\[
 \|F\|_{\mathcal G}\leq TA\epsilon,\quad
 \|Z\|\leq A\epsilon,\quad
 A=1+t_3+t_3^2+t_4\leq(1+t)^2,\quad t=BT(K+\epsilon).
\]
Let $Q=\Pi_{\geq3}\{P,F\}$. The first bracket has the bounded expression
$W=\{H,F\}=-J+Z+Q+\{J,F\}$, hence the exact pullback is
\[
 e^{\operatorname{ad}_F}H=N+Z+P+Q+\mathcal E,
 \qquad
 \mathcal E=\{J,F\}
 +\int_0^1(1-t)e^{t\operatorname{ad}_F}\{W,F\}\,dt.
\]
Estimating $\{W,F\}$ from $\eta_{\nu,6}$ to $\eta_{\nu,7}$ and its
pullback at $\eta_{\nu,8}$, we obtain
\[
 \|\mathcal E\|
 \leq BTA[2+A(1+t)]\epsilon^2
 \leq3BT(1+t)^5\epsilon^2.
\]
Set
\[
 N^+=N+Z+I\Pi_{\leq2}\mathcal E,\quad
 J^+=(1-I)\Pi_{\leq2}\mathcal E,\quad
 P^+=P+Q+\Pi_{\geq3}\mathcal E.
\]
Consequently,
\[
 \|P^+-P\|\leq BK\|F\|_{\mathcal G}+\|\mathcal E\|,
 \qquad
 \sum_{k=0}^2\|\Pi_kJ^+\|\leq3\|\mathcal E\|.
\]

Set $\epsilon_{\nu+1}=C_0\Lambda^\nu\epsilon_\nu^2$ and
$x_\nu=C_0\Lambda^{\nu+1}\epsilon_\nu$.
Then $x_{\nu+1}=x_\nu^2$ and
$\epsilon_{\nu+1}/\epsilon_\nu=x_\nu/\Lambda$.
Assuming inductively that the higher-order perturbation is bounded by
$K_*$, we obtain
$t_\nu\leq t_0 2^{d\nu}$ and therefore
\[
 B_\nu\|F_\nu\|_{\mathcal G}
 \leq D_0(1+t_0)^2 2^{3d\nu}\epsilon_\nu,
 \qquad
 \|\mathcal E_\nu\|\leq\epsilon_{\nu+1}/3.
\]
The first majorant has initial value
$x_0/[9\Lambda(1+t_0)^3]\leq1/144$ and ratio at most
$2^{-3d}x_\nu\leq1/16$. Thus the Lie-series smallness condition holds at
every step. Since $d\geq b+c$, the majorant for
$2^{c\nu}\|F_\nu\|_{\mathcal G}$ has the same ratio bound.
Summing the resulting geometric majorants gives
\[
 \sum_\nu\|Z_\nu\|\leq2(1+t_0)^2\epsilon_0,\qquad
 \sum_\nu B_\nu\|F_\nu\|_{\mathcal G}
 \leq2D_0(1+t_0)^2\epsilon_0,\qquad
 \sum_\nu\|\mathcal E_\nu\|\leq\epsilon_0/8.
\]
Thus $\|N_\nu-N_0\|<3\epsilon_0<4\epsilon_0$, and the inverse
remains available. Moreover,
\[
 \|P_\nu\|\leq K_0+2K_*D_0(1+t_0)^2\epsilon_0+\epsilon_0/8<K_*,
\]
and its cubic projection is bounded by
$\epsilon_0+2t_0(1+t_0)^2\epsilon_0+\epsilon_0/8<3\epsilon_0/2$.
This closes the induction. Summing the exclusions gives the measure
bound; completeness gives convergence on $E_\infty$. Applying $I\Pi_4$
to the update yields the quartic correction. In particular,
\[
 \Pi_4Q_\nu=
 \{P_{6,\nu},F_{0,\nu}\}
 +\{P_{5,\nu},F_{1,\nu}\}
 +\{P_{4,\nu},F_{2,\nu}\}.
\]
Their integrable part is retained until the final normal form is extracted.
\end{proof}

\begin{corollary}[Cubic normalization]
\label{cor:max3-cubic}
Assume that there is a Borel set $E_*\subset E_\infty$ satisfying
$\mu(E_\infty\setminus E_*)\leq\ell_*$, where
$\ell_*<\mu(E_0)-\sum\ell_\nu$.
Suppose that, at suitable smaller widths below $\eta_\infty$, the
homological operator associated with the limiting diagonal has a
support-preserving inverse on degree-three sources, with norm $T_*$,
and that the Poisson bracket and bounds for the Lie-series pullback hold
with constant $B_*$. Identities involving $N_\infty$ use its bounded
first bracket as in Proposition~\ref{prop:max3-newton}. If
$(3/2)B_*T_*\epsilon_0\leq1/16$, then, on $E_*$,
\[
 \{N_\infty,\Gamma_3\}=-P_3,\qquad P_3=\Pi_3P_\infty,
 \qquad\|\Gamma_3\|_{\mathcal G}\leq\tfrac32T_*\epsilon_0,
\]
and
\[
 e^{\operatorname{ad}_{\Gamma_3}}(N_\infty+P_\infty)
 =N_\infty+\Pi_4P_\infty+\tfrac12\{P_3,\Gamma_3\}+R_{\geq5}.
\]
Consequently the new quartic action polynomial is
$Z_{\mathrm{cub}}=\tfrac12 I\{P_3,\Gamma_3\}$, with norm at most
$\tfrac98B_*T_*\epsilon_0^2$.
The generator is affine in the actions and has at most three normal factors.
After the full integrable quartic part is included in the normal form,
the bounded remainder has degree at least four, zero integrable quartic
projection, and norm at most
$2K_*+3\epsilon_0$.
\end{corollary}

\begin{proof}
Since $\{N_\infty,\Gamma_3\}=-P_3$ and
$\operatorname{ad}_{\Gamma_3}$ raises weighted degree by one, the
terms of degree at most four are
\[
 N_\infty+P_3+\Pi_4P_\infty-P_3
 +\{P_3,\Gamma_3\}-\tfrac12\{P_3,\Gamma_3\}.
\]
The bracket estimate gives the bound for $Z_{\mathrm{cub}}$, while
\[
 e^{\operatorname{ad}_{\Gamma_3}}N_\infty
 =N_\infty-\int_0^1e^{t\operatorname{ad}_{\Gamma_3}}P_3\,dt,
\]
which bounds the final bounded perturbation by $2K_*+3\epsilon_0$.

For the quartic correction, write
$P_3=\sum_\lambda p_\lambda e^{ik\theta}y^m z^\alpha\bar z^\beta$,
where $\lambda=(k,m,\alpha,\beta)$ and
$2|m|+|\alpha|+|\beta|=3$.
Put $D_\lambda=k\cdot\omega+\sum_n\Omega_n(\alpha_n-\beta_n)$.
Then $D_\lambda\ne0$ on the support and
$(\Gamma_3)_\lambda=p_\lambda/(iD_\lambda)$.
For $\kappa=(k',m',\alpha',\beta')$, impose
$k'=-k$ and $\gamma=\alpha+\alpha'=\beta+\beta'$.
Differentiating the monomials gives
\begin{align*}
 Z_{\mathrm{cub}}=\frac12
 \sum_{\substack{\lambda,\kappa:\ k'=-k\\
                 \alpha+\alpha'=\beta+\beta'}}
 \frac{p_\lambda p_\kappa}{D_\kappa}
 \bigg[&\sum_{j\in S}j(k_jm'_j-m_jk'_j)
             y^{m+m'-e_j}\rho^\gamma\\
 &+\sum_{n>0,\,n\notin S}n(\alpha_n\beta'_n-\beta_n\alpha'_n)
             y^{m+m'}\rho^{\gamma-e_n}\bigg].
\end{align*}
Zero differentiation terms are omitted. Absolute convergence justifies
the sum over ordered pairs and its integrable projection.
\end{proof}

\subsection{KAM normal form}
\label{sec:max3-concrete-inner}

Let $\mathcal H$ be the real, even, momentum-preserving subspace of
Section~\ref{sec:verified-coefficients}, and let $G$ denote the norm on
$\mathcal C_G^{\mathrm{aff}}(S;0)$. The generators below have weighted
degree at most three. All parameter norms below are total norms at
radius $r$.

\begin{proposition}[KAM normal form]
\label{prop:max3-concrete-inner}
Fix a nonempty finite set \(S\subset\mathbb N\), \(d=|S|\), \(s\geq2\),
\(\sigma=s+4\), \(A\geq1\), and
\[
 C_s=4\left(1+\sum_{n\geq1}(1+n^2)^{-s}\right)^{1/2},
 \qquad \eta_*= \sigma+3+\log(8C_s).
\]
Let \(\eta_f\geq\eta_*\), \(0<h\leq1/32\), \(0<r\leq1\),
\(\tau>d+8\), and \(0<\kappa\leq1\).
Let \(B\subset[100r^2,A^2]^S\) be a rectangular box with positive
volume and shortest side \(\ell>0\), and let \(\mu\) be its normalized
Lebesgue measure. Let \(E_0\subset B\) be Borel with positive measure.
Suppose the initial diagonal is
\[
 N_0=e_0+\sum_{j\in S}\frac{\omega_{j,0}}{j}y_j
                  +\sum_{n>0,n\notin S}\frac{\Omega_{n,0}}{n}\rho_n,
\]
\[
 \omega_{j,0}=j^3+ja_0-6j\xi_j+jf_{j,0},
 \qquad \Omega_{n,0}=n^3+na_0+g_{n,0}.
\]
The coefficients are real and Borel. Assume, for some \(X\geq1\),
\[
 \operatorname{Lip}(f_0;\ell^1)\leq\frac1{24},
 \quad \sup_n\operatorname{Lip}g_{n,0}\leq\frac1{24},
 \quad \sup_{\xi,n}|g_{n,0}|\leq\frac1{10},
 \quad \sup_\xi(\|\xi\|_1+\|f_0\|_1)\leq X-1.
\]
The dispersive diagonal is treated separately; $a_0,e_0$ need only be
finite and Borel.

There are constants
\[
 C_A=C_A(S),\qquad C_D=C_D(S,X,\tau),\qquad
 C_P=C_P(S,A,s,\eta_f),
\]
all at least one and independent of $r,h,\kappa,B,E_0$, the initial
coefficients, and their amplitudes, such that the following holds. Set
\[
 B_0=C_Ar^{-2}h^{-3},\qquad
 T_0=C_D\left(\kappa^{-1}h^{-\tau-4}
                 +r^2\kappa^{-2}h^{-2\tau-5}\right),\qquad D_0=B_0T_0,
\]
\[
 d_*=2\tau+10,\quad \Lambda=2^{6d_*},\quad
 K_*=2(K+\epsilon),\quad t_0=D_0(K_*+\epsilon),\quad
 C_0=9D_0(1+t_0)^5,\quad x_0=C_0\Lambda\epsilon,
\]
where \(\epsilon>0\) and \(K\geq0\). Assume that
\begin{equation}
 4(d+1)\epsilon r^{-4}\leq\frac1{24},\quad
 t_0\leq\frac1{64},\quad x_0\leq\frac1{16},\quad
 \frac{3C_D\kappa}{\ell}<\mu(E_0),\quad
 \frac{4C_PT_0\epsilon}{r^3}\leq\frac1{16}.
 \label{eq:max3-ci-smallness}
\end{equation}
At width \(\eta_{\rm in}=\eta_f+24h\), suppose
\[
 H_0=N_0+J_0+P_0,\qquad
 \Pi_{\leq2}J_0=J_0,\quad IJ_0=0,\quad \Pi_{\geq3}P_0=P_0,
\]
\begin{equation}
 \sum_{j=0}^2\|\Pi_jJ_0\|^{\rm tot}\leq\epsilon,\qquad
 \|P_0\|^{\rm tot}\leq K,\qquad
 \|\Pi_3P_0\|^{\rm tot}\leq\epsilon.
 \label{eq:max3-ci-seed}
\end{equation}
Here \(I\) is the integrable projection. The perturbations obey reality,
momentum conservation, and even symmetry.

Then there exists a Borel set \(E_*\subset E_0\) such that
\begin{equation}
 \mu(E_*)\geq\mu(E_0)-3C_D\kappa/\ell>0.
 \label{eq:max3-ci-loss}
\end{equation}
On this set, the triangular iteration followed by the cubic
homological equation yields generating Hamiltonians \(F_\nu\) of
weighted degree at most two and a generator
\(\Gamma\) of degree three. All are affine in the actions, with fourth
normal index zero. Put
\[
 h_\nu=h2^{-\nu},\qquad
 \eta_{\nu,i}=\eta_f+8h+(16-i)h_\nu\quad(0\leq i\leq8).
\]
The iteration uses \(\kappa_\nu=\kappa2^{-\nu}\) for the current
diagonal, and the threshold \(\kappa\) for the limiting diagonal. Its error
majorant is
\[
 \epsilon_\nu=\frac{x_0^{\,2^\nu}}{C_0\Lambda^{\nu+1}}.
\]
The diagonal increments and higher-order perturbations converge in total coefficient
norm at width \(\eta_f+8h\) to \(\Delta_\infty\) and \(P_\infty\), with
\[
 \|\Delta_\infty\|^{\rm tot}<3\epsilon,\quad
 \|P_\infty\|^{\rm tot}\leq K_*,\quad
 \|\Pi_3P_\infty\|^{\rm tot}\leq3\epsilon/2.
\]
At width \(\eta_f\), the transformed Hamiltonian is \(N_0+Z+R\), where
\[
 Z=\Delta_\infty+I\Pi_4P_\infty
                      +\tfrac12I\{\Pi_3P_\infty,\Gamma\},
 \qquad \Pi_{\geq4}R=R,\qquad I\Pi_4R=0,
\]
\[
 \|R\|^{\rm tot}\leq2K_*+3\epsilon,
\]
\begin{equation}
 \|Z-\Delta_\infty-I\Pi_4P_0\|^{\rm tot}
 \leq[2t_0(1+t_0)^2+1/8]\epsilon+\tfrac98D_0\epsilon^2.
 \label{eq:max3-ci-correction}
\end{equation}
The norms of \(F_\nu\) below are at width \(\eta_{\nu,5}\), and that of
\(\Gamma\) is at width \(\eta_f+7h\):
\begin{equation}
 \sum_{\nu\geq0}\|F_\nu\|_G^{\rm tot}
       \leq2T_0(1+t_0)^2\epsilon,\qquad
 \|\Gamma\|_G^{\rm tot}\leq\tfrac32T_0\epsilon.
 \label{eq:max3-ci-generators}
\end{equation}

For each \(\xi\in E_*\), define
\[
 U_v(R_0;\xi)=
 \left\{u\in H^v_0:\sum_{j\in S}\big||u_j|^2-\xi_j\big|<R_0^2r^2,\
             \|(I-\Pi_S)u\|_v<R_0r\right\}.
\]
Here $s\leq v\leq\sigma$, $0<R_0\leq3$, $\Pi_S$ retains both signs of
$S$, and $\|u\|_v^2=\sum_{n>0}(1+n^2)^v|u_n|^2$, with
$u_{-n}=\overline{u_n}$. In the original Fourier variables we have
$y_j=|u_j|^2-\xi_j$ and $e^{i\theta_j}=u_j/|u_j|$.
The Hamiltonian time-one flows \(\mathcal T_\nu,\mathcal T_\Gamma\) exist.
On \(U_\sigma(1/4;\xi)\) their compositions
\[
 \Phi_N=\mathcal T_0\circ\cdots\circ\mathcal T_N\circ\mathcal T_\Gamma
\]
converge pointwise in \(H^\sigma\) and uniformly in \(H^{\sigma-1}\) to
\(\Phi\). All intermediate flow trajectories remain in \(U_\sigma(1;\xi)\).
The map is continuous into \(H^\sigma\) and \(C^1\) into
\(H^{\sigma-1}\). With
\[
 E_G=\left(2(1+t_0)^2+\tfrac32\right)T_0\epsilon
       \leq4T_0\epsilon,
\]
one has, uniformly in the parameter,
\[
 \|\Phi(u)-u\|_s\leq3C_PE_G/r^2,
\]
\begin{equation}
 e^{-3C_PE_G/r^3}\|u-w\|_v
 \leq\|\Phi(u)-\Phi(w)\|_v
 \leq e^{3C_PE_G/r^3}\|u-w\|_v,\qquad s\leq v\leq\sigma-1.
 \label{eq:max3-ci-metric}
\end{equation}
The variational extension is invertible on these spaces, obeys the
same operator and inverse bounds, and preserves the Gardner form on $H^s$.
It extends the derivative at each point of $U_\sigma(1/4;\xi)$ from the dense subspace $H^\sigma$.
The limiting map satisfies the scalar and vector-field identities
\begin{equation}
 H_0\circ\Phi=N_0+Z+R,
 \qquad
 X_{H_0}(\Phi(u))=E_\Phi(u)X_{N_0+Z+R}(u)
 \quad\hbox{in }H^{s+1}.
 \label{eq:max3-ci-conjugacy}
\end{equation}
Thus the image of the reference torus \(y=0,\ u_{S^c}=0\) is a compact
invariant torus for the initial Hamiltonian \(H_0\). Its induced rotation
has the tangential frequencies of \(N_0+\Delta_\infty\), and every
image orbit belongs to
\(C(\mathbb R;H^\sigma)\cap C^1(\mathbb R;H^{s+1})\).
\end{proposition}

\begin{proof}
We first verify the inverse estimates throughout the required
neighborhood of the initial diagonal. Write
\[
 \Delta=e+\sum_{j\in S}d_jy_j+\sum_{n>0,n\notin S}b_n\rho_n,
 \qquad \|\Delta\|^{\rm tot}\leq4\epsilon.
\]
The ordered normal coefficient is \(b_n/2\). If \(c=\lim_n b_n\),
the degree-two norm implies
\[
 |c|+\textstyle\sum_j|d_j|
       \leq8\epsilon/r^2,\qquad
 \sup_n n|b_n-c|\leq4\epsilon/r^2.
\]
Individually, $|c|,\sum_j|d_j|\le4\epsilon/r^2$; parameter
differences give the same bounds with an extra factor $r^{-2}$.
The frequency correction is
\[
 \delta a=c,\qquad \delta f_j=d_j-c,\qquad
 \delta g_n=n(b_n-c).
\]
Condition \eqref{eq:max3-ci-smallness} keeps both reduced Lipschitz constants at most
\(1/12\), the normal supremum below \(1/5\), and
\(\sup(\|\xi\|_1+\|f\|_1)<X\).
Propositions~\ref{prop:max4-exclusion} and~\ref{prop:max4-inverse}
therefore give an inverse for the homological operator at each updated diagonal,
after a single Borel parameter exclusion common to the three equations.
The loss is \(C_D\kappa2^{-\nu}/\ell\), and the inverse norm is bounded by
\[
 T_\nu\leq T_0\,2^{(2\tau+7)\nu}.
\]
Theorems~\ref{thm:coefficient-initial} and~\ref{thm:max-coefficient-algebra}
give the projection, bracket and Lie estimates, with
$B_\nu=B_0 2^{3\nu}$. Proposition~\ref{prop:max3-newton} applies with
$a=3$, $b=2\tau+7$, $c=0$, $d=d_*$, proving convergence and the
first generator bound.

The limiting diagonal lies in the same \(4\epsilon\) neighborhood.
Apply Proposition~\ref{prop:max4-inverse} once more at threshold \(\kappa\), using
the reserved interval of width \(8h\). This gives the cubic inverse
with constant \(T_0\); the smallness condition for its Lie series follows from
\[
 \tfrac32D_0\epsilon\leq t_0/2\leq1/128.
\]
Corollary~\ref{cor:max3-cubic} proves
\eqref{eq:max3-ci-correction}--\eqref{eq:max3-ci-generators}.
Summing the exclusions gives \eqref{eq:max3-ci-loss}.

All generators are real, affine in $y$, of weighted degree at most three,
and satisfy $\mu_4=0$. Proposition~\ref{prop:flows-affine-realization}
and the bound $\sum_\nu\|F_\nu\|_G+\|\Gamma\|_G\leq E_G$
verify the hypotheses of Propositions~\ref{prop:flows-finite}
and~\ref{prop:flows-inner}. These give the finite flows and their limiting
composition on $U_\sigma(1/2;\xi)$.
The final cubic flow maps \(U_\sigma(1/4;\xi)\) into that domain:
its action change is at most \(r^3/16\), and its normal norm change
is at most \(r^2/16\). Thus
\[
 r^2/16+r^3/16<r^2/4,\qquad r/4+r^2/16<r/2.
\]
Composition of the metric and variational estimates gives
\eqref{eq:max3-ci-metric} and the displacement bound.

To prove the scalar conjugacy identity, fix \(\xi\). The bounded diagonal
\(\mathcal B_\nu=N_\nu-H_2\) belongs pointwise to the coefficient
space: its normal coefficient is a common constant plus \(g_n/n\).
Its first action remainder is therefore \(O(1/n)\).
The separate dispersive bracket is the bounded expression
\[
 \{H_2,F_\nu\}
     =\{N_\nu,F_\nu\}-\{\mathcal B_\nu,F_\nu\}.
\]
It belongs to the coefficient space at width \(\eta_{\nu,6}\).
Propositions~\ref{prop:flows-evaluation} and~\ref{prop:flows-lie}
identify the brackets and pullbacks, giving
\[
 H_{\nu+1}=H_\nu\circ\mathcal T_\nu.
\]
Composition on the common domain in $H^\sigma$, followed by coefficient
evaluation and strong convergence, yields
\[
 H_0\circ\left(\lim_N\mathcal T_0\circ\cdots\circ\mathcal T_N\right)
      =N_\infty+P_\infty.
\]
Applying the same argument to $\Gamma$ proves the scalar identity in
\eqref{eq:max3-ci-conjugacy}.

By the evaluation theorem the bounded gradients lie in \(H^\sigma\),
and the full Hamiltonian fields lie in \(H^{s+1}\).
Proposition~\ref{prop:flows-conjugacy} then gives the vector identity
in \eqref{eq:max3-ci-conjugacy}. On the reference torus all first
derivatives of the quartic action polynomial and of \(R\) vanish.
The remaining field is the rotation generated by $N_\infty$.
The metric bounds and strong continuity give the asserted invariant
torus and global solution curves.
\end{proof}

Proposition~\ref{prop:repeatable-opening} constructs Hamiltonians
satisfying these hypotheses. At each finite stage, all constants are
fixed before the next radius is chosen.

\section{Opening estimates}
\label{sec:verified-opening}
Following \cite{BGR25,BP26}, we collect the estimates required when
normal modes are added to the tangential set. We first estimate the
opening substitution in the coefficient and total parameter norms.
We then control the frequency corrections and the small-divisor
hypotheses under successive openings. The last subsection provides
the domain estimates for the associated Hamiltonian transformations.
Throughout, we use the spaces and projections of
Section~\ref{sec:coefficient-definitions}, and write
$\|P\|_{S,\eta,R}=\|P\|_{\mathcal H_{\eta,R}(S)}$.

\subsection{Opening map}
\label{sec:max3-opening}

Let $S,T\subset\mathbb N$ be finite and disjoint, with $T\ne\varnothing$,
and put $S'=S\cup T$ and $L=\max(1,\max S')$.
Let $E\subset\mathbb R^S$ be a nonempty Borel parameter set, and fix
$a,\rho,\varrho>0$. For $\xi\in B_T=[a^2/2,a^2]^T$, define
\[
 u_t=(\xi_t+y_t)^{1/2}e^{i\theta_t},\qquad
 u_{-t}=(\xi_t+y_t)^{1/2}e^{-i\theta_t}\quad(t\in T),
\]
using the Taylor expansion of the positive square root at $y_t=0$.
We denote this substitution by $\mathcal O_T$. Total norms use the global
$\ell^1$ metric and the square of the indicated radius. In particular,
the metric on $E\times B_T$ is
$\|(x,\xi)-(x',\xi')\|_1=\|x-x'\|_1+\|\xi-\xi'\|_1$.

\begin{proposition}[Estimates for a finite opening]
\label{prop:max3-opening}
Let $P$ be a family in $\mathcal H_{\eta,\rho}(S)$ with finite total norm
on $E$. Let $0\leq\eta'\leq\eta$ and suppose
\[
 t_*=\frac{2e^{2\eta'}\varrho^2}{a^2}\leq\frac14,
 \qquad
 \vartheta=\frac{\varrho}{\rho}
       +\frac{2|T|a}{\rho\sqrt{1-t_*}}\leq\frac12.
\]
Then
\[
 \|\mathcal O_TP\|^{\mathrm{tot}}_{S',\eta',\varrho;E\times B_T}
 \leq C(1+L)\|P\|^{\mathrm{tot}}_{S,\eta,\rho;E},
\]
where $C$ is absolute. Both action limits commute with the coefficient
sums defining $\mathcal O_TP$. The substitution preserves reality,
momentum conservation and evenness.
\end{proposition}

\begin{proof}
Fix a coefficient block $p=p_{k,m,q}$ before opening. For multiplicities
$A_t,B_t$ of the opened modes, put
\[
 h=\sum_t(A_t+B_t),\quad q'=q-h,\quad
 \nu_t=A_t-B_t,\quad \lambda_t=(A_t+B_t)/2.
\]
For an ordered tuple $K$ of the remaining $q'$ normal indices, the contribution
to the coefficient indexed by $(k,\nu),(m,\beta),K$ is
\begin{equation}
\label{eq:max3-kernel}
 \frac{q!}{q'!\prod_t A_t!B_t!}
 p\bigl(K,\{t^{A_t},(-t)^{B_t}\}_{t\in T}\bigr)
 \prod_t \binom{\lambda_t}{\beta_t}\xi_t^{\lambda_t-\beta_t}.
\end{equation}
This follows by counting the ordered positions and expanding
$(\xi_t+y_t)^{\lambda_t}$.
If $O$ is any fixed tuple in $\pm T$ and $f(K)=p(K,O)$, then
\begin{equation}
\label{eq:max3-slice}
 a_{S'}(f)\leq5(1+L)a_S(p).
\end{equation}
Indeed, with $\mathfrak a=a_S(p)$,
\[
 |p|\le3\mathfrak a,\qquad |L_1p|\le2\mathfrak a,\qquad
 |L_2p|\le\mathfrak a,\qquad
 L_bf(V)=L_bp(V,O)\quad(b=1,2).
\]
Here $O$ is fixed in both limits.
For strict highest pairs in the remaining normal tuple,
\[
 \begin{aligned}
 (1+n)|f(n,-n,V)-L_1f(V)|
 &\le\begin{cases}
       \mathfrak a,&n>L,\\
       5(1+L)\mathfrak a,&n\le L,
      \end{cases}\\
 (1+h)|L_1f(h,-h,V)-L_2f(V)|
 &\le\begin{cases}
       \mathfrak a,&h>L,\\
       3(1+L)\mathfrak a,&h\le L.
      \end{cases}
 \end{aligned}
\]
For $n,h>L$, use the remainder bounds for $p$; for $n,h\le L$, use
the uniform bounds. Off the pairing sets, the three components are
bounded by $3\mathfrak a$, $2\mathfrak a$ and $\mathfrak a$.
This proves \eqref{eq:max3-slice}.

For $\lambda\ge0$,
\[
 \left|\binom{\lambda}{b}\right|
 \le\frac{\lambda(\lambda+1)\cdots(\lambda+b-1)}{b!},
\]
so $(1-t)^{-\lambda}$ majorizes the absolute Taylor coefficients of
$(1+t)^\lambda$. Set $b=|\beta|_1$,
$c_\beta=\prod_t|\binom{\lambda_t}{\beta_t}|$ and
$M_0=(1-t)^{-h/2}$. Applying $t\partial_t$ and $(t\partial_t)^2$ to
this majorant gives
\[
 \sum_\beta b c_\beta t^b\leq \frac{ht}{2(1-t)}M_0,
 \qquad
 \sum_\beta b^2c_\beta t^b
 \leq\left(\frac{ht}{2(1-t)}
 +\frac{(h/2)(h/2+1)t^2}{(1-t)^2}\right)M_0.
\]
Let $W,W'$ be the input and output block weights. Since $h\le q\le W$,
$|\nu|_1\le h$,
$\xi_t^{\lambda_t-\beta_t}\le a^{2\lambda_t}(2/a^2)^{\beta_t}$ and
$W'\le2W+b$, the zeroth and first moments bound the Taylor sum by a
constant times
\[
 \sum_{q'+\sum_t(A_t+B_t)=q}
 \frac{q!}{q'!\prod_tA_t!B_t!}
 \left(\frac{\varrho}{\rho}\right)^{q'}
 \left(\frac{a}{\rho\sqrt{1-t_*}}\right)^h
 =\vartheta^q.
\]
Together with $(\varrho/\rho)^{2|m|_1}\le1$ and
\eqref{eq:max3-slice}, this proves the supremum estimate.

For the new parameter derivative, the scalar Taylor coefficient $g_\beta$
satisfies
\[
 \sum_t|\partial_{\xi_t}g_\beta|
 \leq \frac{2}{a^2}(h/2+b)|g_\beta|.
\]
The two moment bounds give a factor
$C(1+L)a^{-2}(1+q)\vartheta^q$ relative to the original block norm.
Now $\sup_q(1+q)\vartheta^q<\infty$, and line segments in the new
parameters remain in $B_T$. Split a two-point difference at
$(x',\xi)\in E\times B_T$. Applying the same substitution estimate to
the old-parameter difference and using
\[
 \varrho^2\le\rho^2,\qquad \varrho^2/a^2\le1/8
\]
proves the total-norm estimate.

The preceding majorants are summable uniformly in each marked index,
so both action limits pass through the coefficient sums. Momentum and
parity are preserved because
\[
 \sum_{t\in T}t(A_t-B_t)=\sum_{t\in T}t\nu_t,\qquad
 q'+\sum_{S'}k'_j-q-\sum_Sk_j=-2\sum_{t\in T}B_t.
\]
Reality follows by exchanging $A$ and $B$ and reversing the angles.
\end{proof}

\begin{corollary}[Jets and integrable quartics]
\label{cor:max3-opening-support}
A monomial of weighted degree $D$ containing $h$ factors at the opened
normal sites gives terms of degrees
$D-h+2|\beta|_1$, and
\[
 \Pi_{\mathrm{int}}\mathcal O_TP
 =\mathcal O_T\Pi_{\mathrm{int}}P.
\]
Suppose that $R$ has weighted degree at least four and
$\Pi_4\Pi_{\mathrm{int}}R=0$. Decompose its terms of degrees four and five
as $B+G$, where $B$ contains the monomials with $h\geq D-3$ and $G$ those
with $h<D-3$. The sum of their total norms is at most
$C(1+L)\|\Pi_{4,5}R\|^{\mathrm{tot}}$.
Moreover, $B$ is affine in the old actions, its second action limit
vanishes, and
\[
 \mu_4(B)\le L,\qquad \Pi_{\mathrm{int}}B=0,\qquad
 \Pi_{\le3}\mathcal O_TG=\Pi_{\mathrm{int}}\mathcal O_TG=0.
\]
The opened jet of degree at most three is affine in all actions and has
at most three normal factors. Its fourth normal index and second action
limit therefore vanish.

For a retained integrable polynomial
$Z=c+\sum_i d_iY_i+\sum_{i,j}B_{ij}Y_iY_j$, with $B_{ij}=B_{ji}$,
assume finite supremum and global Lipschitz coefficient norms with respect to the old parameters
on $E$ at width $\eta$ and radius one. Here $Y_i=y_i$ for $i\in S$
and $Y_i=u_iu_{-i}$ for positive $i\notin S$; the sums are interpreted
in the ordered coefficient class. Put $Y'_t=y_t$ for $t\in T$ and
$Y'_i=Y_i$ otherwise. Opening gives the exact translation $Y=Y'+\zeta$, with
$\zeta_t=\xi_t$ for $t\in T$ and zero otherwise. In particular,
\[
 \Pi_4\mathcal O_TZ=\sum_{i,j}B_{ij}Y'_iY'_j,
 \qquad \|\Pi_4\mathcal O_TZ\|^{\mathrm{tot}}_{\eta',r}\leq C_Zr^4
 \quad(0<r\leq1).
\]
Here $C_Z$ depends only on the original coefficient bounds and their
Lipschitz seminorms, the finite sets and the widths. The quartic
coefficients are independent of the new parameters, and $C_Z$ is
independent of $a$. The dispersive quadratic Hamiltonian is transformed
by exact substitution.
\end{corollary}

\begin{proof}
The degree formula follows from \eqref{eq:max3-kernel}. An opened
monomial is integrable precisely when $k=0$, $A_t=B_t$ for every $t\in T$,
and the remaining normal indices are balanced. These are also the
integrability conditions before opening.
The indicator defining $B$ is unchanged when the index of a marked pair exceeds $L$,
so the proof of \eqref{eq:max3-slice} applies to it and its complement.
On $B$, $q-h\le D-h\le3$, whence $\mu_4(B)\le L$ and the second
action limit vanishes.
Also $|m|_1\ge2$ would imply $h\le q\le D-4$, a contradiction.
Integrable monomials have even degree, so $\Pi_4\Pi_{\mathrm{int}}R=0$
gives the asserted integrable projections. On $G$, $D-h\ge4$.

Finally, $\Pi_4Z(Y'+\zeta)=\sum_{i,j}B_{ij}Y'_iY'_j$.
The coefficient estimate follows from \eqref{eq:max3-slice}; the
supremum and weighted Lipschitz terms carry factors $r^4$ and $r^6$.
We omit the remaining details.
\end{proof}

Opening may lower degree six to degree zero. We therefore estimate the
opened remainder at an intermediate radius by
Proposition~\ref{prop:max3-opening}, before restricting the radius.

\subsection{Frequency corrections}
\label{sec:frequency-replacement}

We use the frequency decomposition of Bernier and Procesi~\cite{BP26},
with the factor $n$ arising from the Gardner bracket. Thus an action
derivative $F_n$ contributes $nF_n$ to the frequency.

Fix an infinite set \(S_\infty\subset\mathbb N\) with
\[
 L=\sum_{n\in S_\infty}\frac1n<\infty.
\]
Let \(A\subset S_\infty\) be finite, \(d=|A|\),
\(J=\max(A\cup\{0\})\), and \(N=\mathbb N\setminus A\).
For the parameters \(\alpha\) associated with \(A\), ranging over a nonempty Borel set
\(E\subset\mathbb R^A\), write
\[
 Y_j=\rho_j-\alpha_j\quad(j\in A),\qquad
 Y_n=\rho_n\quad(n\in N).
\]
Parameter Lipschitz seminorms are taken with respect to the $\ell^1$
metric, as in Section~\ref{sec:coefficient-definitions}.

\subsubsection{Action--Toeplitz decomposition}

Let \(Z_\alpha\) be a real integrable formal polynomial of weighted
degree at most four, with ordered expansion
\[
 Z_\alpha=
 \sum_{\substack{m,q\\2|m|_1+q\le4}}
 \sum_{\mathbf j\in(\pm N)^q}
 p_{m,q}(\alpha;\mathbf j)Y_A^m u_{j_1}\cdots u_{j_q}.
\]
The coefficients follow the convention of \eqref{eq:coefficient-expansion}.
Integrability identifies $Z$ with a polynomial of degree at most two in
$Y$, through $u_nu_{-n}=Y_n$.

For \(q=0\), put \(a_{m,0}=|p_{m,0}|\).  For \(q=2\), define
\[
 \lambda_m=\lim_{n\to\infty,\ n\in N}p_{m,2}(n,-n),\qquad
 a_{m,2}=\max\left\{|\lambda_m|,
 \sup_{n\in N}(1+n)|p_{m,2}(n,-n)-\lambda_m|\right\}.
\]
For \(q=4\), necessarily \(m=0\); write
\[
 \lambda_h^{[4]}=\lim_{n\to\infty,\ n\in N}
 p_{0,4}(n,-n,h,-h),\qquad
 \mu=\lim_{h\to\infty,\ h\in N}\lambda_h^{[4]},
\]
and set
\[
\begin{split}
 a_{0,4}=\max\bigg\{&
 \sup_{\substack{n>h\\n,h\in N}}(1+n)
 |p_{0,4}(n,-n,h,-h)-\lambda_h^{[4]}|,\\
 &\sup_{h\in N}(1+h)|\lambda_h^{[4]}-\mu|,
 |\mu|,\ \sup_{n\in N}|p_{0,4}(n,n,-n,-n)|\bigg\}.
\end{split}
\]
The limits are assumed to exist.  Define
\[
 \|Z\|_{\eta,r}=
 \sum_{2|m|_1+q\le4}
 \max(1,|m|_\infty,q)e^{\eta(2|m|_1+q)}
 r^{2|m|_1+q}a_{m,q},\qquad \eta\ge0,\quad 0<r\le1,
\]
and
\[
 Q=\sup_{\alpha\in E}\|Z_\alpha\|_{\eta,r},\qquad
 D=\operatorname{Lip}_E(Z;\|\cdot\|_{\eta,r}).
\]
These are precisely the norms of the integrable blocks of weighted
degree at most four in the coefficient space of Section~\ref{sec:verified-coefficients}.

\begin{proposition}[Frequency estimates under opening]
\label{prop:frequency-interface}
Let $Q,D<\infty$. Then $Z$ admits the linear decomposition
\[
\begin{split}
 Z={}&c_0+c_2\sum_{i\in N}Y_i+\sum_i a_iY_i
 +c_3\sum_{\substack{j<i\\i,j\in N}}Y_iY_j\\
 &+\sum_{i,j}b_{ij}Y_iY_j
 +\sum_{\substack{j<i\\i\in N}}e_jY_iY_j,\qquad b_{ij}=b_{ji}.
\end{split}
\]
Put \(k_{ii}=1\), \(k_{ij}=(1+i+j)^{-1}\) for \(i\ne j\), and
\[
 B_2=|c_2|+\sup_i(1+i)|a_i|,\qquad
 B_4=|c_3|+\sup_{i,j}\frac{|b_{ij}|}{k_{ij}}+\sup_i(1+i)|e_i|.
\]
With \(C_A=100(1+J)\), one has
\[
 B_2\le C_AQr^{-2},\qquad B_4\le C_AQr^{-4}.
\]
The same estimates apply to coefficient differences with \(Q\) replaced
by \(D\|\alpha-\alpha'\|_1\).

Let \(A\subset S\subset S_\infty\), initially with \(S\) finite, and let
\(\zeta\ge0\) be supported in \(T=S\setminus A\), with
\(t=\|\zeta\|_1\le r^2\).  Set
\[
 a_Z=c_2+c_3t+\sum_i e_i\zeta_i,\qquad
 F_n=\partial_{Y_n}Z(0,\zeta)-a_Z,\qquad
 G_n=nF_n\quad(n\notin S).
\]
Then \(a_Z=\lim_{n\to\infty}\partial_{Y_n}Z(0,\zeta)\), and
\[
 F_n=a_n+2\sum_i b_{ni}\zeta_i-c_3\zeta_n
       +e_n\sum_{i>n}\zeta_i-\sum_{i\ge n}e_i\zeta_i
       \quad(n\notin A),
\]
\[
 F_n=a_n-c_2+2\sum_i b_{ni}\zeta_i
       +e_n\sum_{i>n}\zeta_i-c_3t-\sum_i e_i\zeta_i
       \quad(n\in A).
\]
On every Borel subset of the admissible pairs \((\alpha,\zeta)\),
with metric \(\|\Delta\alpha\|_1+\|\Delta\zeta\|_1\), one has
\[
 \sup\left(\|F_S\|_1+\sup_{n\notin S}|G_n|+|a_Z|\right)
 \le K_A Qr^{-2},
\]
\[
 \operatorname{Lip}(F_S;\ell^1)
 +\sup_{n\notin S}\operatorname{Lip}G_n
 +\operatorname{Lip}a_Z
 \le K_A(Q+r^2D)r^{-4},
 \qquad
 K_A=200(1+J)(10+3L+2d).
\]
The constants are uniform in $S\setminus A$, the lower amplitude bounds,
the Borel set $E$ and $\eta\ge0$.
The assertions also hold for \(S=S_\infty\) and
\(\zeta\in\ell^1(S_\infty\setminus A)\), with the same mass bound;
then \(F_{S_\infty}\in\ell^1\), and \(G\) is indexed by
\(\mathbb N\setminus S_\infty\).
\end{proposition}

\begin{proof}
Write \(\lambda=\lambda_0\) for the quadratic normal limit and
\(\lambda^{(j)}=\lambda_{\epsilon_j}\) for the block with one tangential action and two normal
variables, where \(\epsilon_j\) is the unit multiindex at \(j\in A\).
Choose
\[
 c_0=p_{0,0},\qquad c_2=2\lambda,\qquad c_3=24\mu.
\]
For \(n,h\in N\), define
\[
\begin{split}
 a_n&=2(p_{0,2}(n,-n)-\lambda),\\
 e_n&=24(\lambda_n^{[4]}-\mu),\\
 b_{nn}&=6p_{0,4}(n,n,-n,-n),\\
 b_{nh}&=12\bigl(p_{0,4}(n,-n,h,-h)
                  -\lambda_{\min(n,h)}^{[4]}\bigr)\quad(n\ne h).
\end{split}
\]
For \(j\in A,n\in N\), define
\[
 a_j=p_{\epsilon_j,0},\qquad e_j=2\lambda^{(j)},\qquad
 b_{jn}=b_{nj}=p_{\epsilon_j,2}(n,-n)
                      -\lambda^{(j)}\mathbf1_{\{n>j\}}.
\]
For \(j,k\in A\), set
\[
 b_{jj}=p_{2\epsilon_j,0},\qquad
 b_{jk}=\tfrac12p_{\epsilon_j+\epsilon_k,0}\quad(j\ne k).
\]
The multiplicities of $Y_n$, $Y_nY_h$ ($n\ne h$), and $Y_n^2$ in the
ordered expansion are $2$, $24$, and $6$, respectively. Hence, for
$n>h$ in $N$,
\[
 c_3+e_h+2b_{nh}=24p_{0,4}(n,-n,h,-h);
\]
and, for $j\in A$, $n\in N$,
\[
 2b_{jn}+e_j\mathbf1_{\{n>j\}}=2p_{\epsilon_j,2}(n,-n).
\]
The other coefficients agree by definition.

If \(\|Z\|_{\eta,r}\le H\), every degree-\(\ell\) block seminorm is
at most \(Hr^{-\ell}\).  Hence
\[
 B_2\le[2+\max(2,1+J)]Hr^{-2}.
\]
For normal off-diagonal pairs,
\[
 |b_{nh}|/k_{nh}
 \le12Hr^{-4}\frac{1+n+h}{1+\max(n,h)}
 \le24Hr^{-4}.
\]
The diagonal, mixed and tangential coefficients satisfy
\[
 \begin{aligned}
 \sup_{n\in N}|b_{nn}|&\le6Hr^{-4},\\
 \sup_{j\in A,\,n\in N}\frac{|b_{jn}|}{k_{jn}}
     &\le2(1+2J)Hr^{-4},\\
 \sup_{j,k\in A}\frac{|b_{jk}|}{k_{jk}}
     &\le(1+J)Hr^{-4}.
 \end{aligned}
\]
Together with $|c_3|\le24Hr^{-4}$ and
$\sup_i(1+i)|e_i|\le24(1+J)Hr^{-4}$, these estimates give
$B_4\le72(1+J)Hr^{-4}$. Applying the linear decomposition to
$Z_\alpha-Z_{\alpha'}$ proves the coefficient-difference estimates.

Differentiation of the displayed representation gives the two
formulas for \(F_n\).  The coefficient of the new action \(\zeta_m\), \(m\in T\), in this
affine expression is
\[
 R_{nm}=
 \begin{cases}
 2b_{nm}-c_3\delta_{nm}
       +e_n\mathbf1_{\{m>n\}}-e_m\mathbf1_{\{m\ge n\}},
       &n\notin A,\\
 2b_{nm}+e_n\mathbf1_{\{m>n\}}-c_3-e_m,&n\in A.
 \end{cases}
\]
For each $m\in T$, the $b$-terms satisfy
\[
 \sum_{n\in S}2|b_{nm}|\le2(1+L)B_4,\qquad
 \sum_{\substack{n\in S\\n\ne m}}\frac1{1+n+m}\le L.
\]
The $e_n$-terms contribute at most $LB_4$. The $e_m$-terms occur in at most
$m+d$ rows, and the $c_3$-terms in at most $1+d$ rows; each contributes
at most $(1+d)B_4$. Hence
\[
 \sum_{n\in S}|R_{nm}|\le(4+3L+2d)B_4.
\]
If \(n\notin S\), then \(n\ne m\) and the diagonal is absent.
The inequalities
\[
 \frac n{1+n+m}\le1,\qquad
 \frac n{1+n}\le1,\qquad
 \frac n{1+m}\le1\quad(m>n)
\]
give \(n|R_{nm}|\le4B_4\).  The common slope is
\(c_3+e_m\), of size at most \(2B_4\).
The constant terms satisfy
\[
 \sum_{n\in S}|a_n-c_2\mathbf1_{\{n\in A\}}|
 \le(L+d)B_2,\qquad
 \sup_{n\notin S}n|a_n|\le B_2.
\]
Putting \(K_0=10+3L+2d\), these separate bounds give
\[
 \|F_S\|_1+\sup_{n\notin S}|G_n|+|a_Z|
 \le K_0(B_2+tB_4),
\]
and the sum of the three Lipschitz constants with respect to the new actions is
at most \(K_0B_4\).

At fixed \(\zeta\), coefficient differences give a total
Lipschitz constant with respect to \(\alpha\) at most
\(K_0(C_ADr^{-2}+r^2C_ADr^{-4})\).
Splitting the difference into a coefficient difference and an action
difference gives the joint Lipschitz bound
\[
 C_AK_0(2Dr^{-2}+Qr^{-4})
 \le2C_AK_0(Q+r^2D)r^{-4}.
\]
The size bound is \(2C_AK_0Qr^{-2}\), proving the claim.

For $\zeta\in\ell^1$, the same column bounds give absolute convergence
and the asserted estimates by passage through finite subsets. Moreover,
\[
 \left|\sum_i b_{ni}\zeta_i\right|
 \le B_4\left(|\zeta_n|+\frac{\|\zeta\|_1}{1+n}\right)\longrightarrow0.
\]
Together with $|a_n|\le B_2/(1+n)$, $|e_n|\le B_4/(1+n)$, and the
summability of $\zeta$, the formula for $F_n$ implies
$\partial_{Y_n}Z(0,\zeta)\to a_Z$.
\end{proof}

\subsubsection{Action centers}

Set
\[
 h_n(\alpha,Y)=\partial_{Y_n}Z_\alpha(Y)
 -\left(c_2(\alpha)+c_3(\alpha)\sum_{i\in N}Y_i
                         +\sum_i e_i(\alpha)Y_i\right).
\]
At $\rho_A=\alpha$ and $\rho_A'=\alpha'$, we have $Y_A=Y_A'=0$.
Thus the difference at fixed normal actions $\zeta$ is
\[
 h_n(\alpha,0,\zeta)-h_n(\alpha',0,\zeta).
\]
The coefficient-difference estimate of
Proposition~\ref{prop:frequency-interface} applies directly, without a
restriction on $\|\alpha-\alpha'\|_1/r^2$.

\subsubsection{Summation}

Let \(A_p\uparrow S_\infty\) be an increasing sequence of finite
sets. Suppose that a correction \(Z^{(p)}\) is obtained at a stage
with tangential set \(A_p\), parameter set \(E_p\), radius
\(0<r_p\le1\), and norms \(Q_p,D_p\) as above.  Define
\[
 K_p=200(1+\max(A_p\cup\{0\}))(10+3L+2|A_p|),\qquad
 \delta_p=K_p(Q_p+r_p^2D_p)r_p^{-4}.
\]
For each finite \(P\), let \(X_P\subset[0,\infty)^{A_P}\) be Borel,
and extend its elements by zero outside $A_P$. Assume, on $X_P$, that
\[
 \sum_p\delta_p<\infty,\qquad
 \xi|_{A_p}\in E_p,\qquad
 \sum_{j\in A_P\setminus A_p}\xi_j\le r_p^2\quad(p\le P).
\]
At subsequent stages, $Z^{(p)}$ is evaluated at the new center, with
its coefficients unchanged as functions of $\xi|_{A_p}$.
These conditions are verified in Proposition~\ref{prop:coherent-production}.

For \(p\le P\), evaluate \(Z^{(p)}\) at its later center and denote
the resulting quantities by \(a_{p,P},F^{p,P},G^{p,P}\).  Set
\[
 a_P=12\sum_{j\in A_P}\xi_j+\sum_{p\le P}a_{p,P},\qquad
 f_{j,P}=\sum_{p\le P}F_j^{p,P},\qquad
 g_{n,P}=\sum_{p\le P}G_n^{p,P}.
\]
Proposition~\ref{prop:frequency-interface}, applied separately
to each summand on any admissible Borel parameter set, gives
\[
\begin{split}
 &\operatorname{Lip}(f_{A_P,P};\ell^1)
 +\sup_{n\notin A_P}\operatorname{Lip}g_{n,P}\\
 &\hspace{15mm}
 +\operatorname{Lip}\left(a_P-12\sum_{j\in A_P}\xi_j\right)
 \le\sum_{p\le P}\delta_p,
\end{split}
\]
and the corresponding combined supremum bound is
\(\sum_{p\le P}\delta_pr_p^2\).
In particular, the bound
\(Q_p+r_p^2D_p\le r_p^{4+1/3}\) is sufficient whenever
\(\sum_pK_pr_p^{1/3}\) is suitably small.
For intermediate integrable updates, we include their total norms in the
same sum.

By \eqref{ham:frequencies}, the tangential and normal frequencies at stage $P$ are
\[
 \omega_{j,P}=j^3+ja_P-6j\xi_j+jf_{j,P},\qquad j\in A_P,
\]
\[
 \Omega_{n,P}=n^3+na_P+g_{n,P},\qquad n\notin A_P.
\]
For a divisor satisfying momentum conservation,
\(\sum_{j\in A_P}jk_j+\sum_{a=1}^qn_a=0\),
the common term cancels.  Extend the normal remainders by oddness to negative indices. The
frequency correction in the divisor then has global Lipschitz constant
at most
\[
 \left(\max_{j\in A_P}|jk_j|+q\right)\sum_{p\le P}\delta_p.
\]

Finally, let \(\xi\in\ell^1(S_\infty)\) be nonnegative,
with \(\xi|_{A_p}\in E_p\) and
\(\sum_{j\in S_\infty\setminus A_p}\xi_j\le r_p^2\) for every
\(p\).  Each correction can then be
evaluated at that entire tail.  The bounds
\(\delta_pr_p^2\) imply absolute uniform convergence of the
common correction in \(\mathbb R\), the normalized tangential
correction in \(\ell^1(S_\infty)\), and the weighted correction to the normal frequencies
outside \(S_\infty\) in \(\ell^\infty(\mathbb N\setminus S_\infty)\).
The sums have the same two-point bounds with
\(\sum_p\delta_p\).

To identify these sums with the finite-stage limits, put
$\xi^{[P]}=\mathbf1_{A_P}\xi$. For fixed $p$, the tails used to evaluate
$Z^{(p)}$ satisfy
\[
 \|\xi^{[P]}|_{S_\infty\setminus A_p}
       -\xi|_{S_\infty\setminus A_p}\|_1\longrightarrow0.
\]
Proposition~\ref{prop:frequency-interface} therefore gives convergence
of each retained correction. For its tangential part, writing $F^p(\xi)$
for the $\ell^1(S_\infty)$ correction evaluated at the full tail, we have
\[
 \|\mathbf1_{A_P}F^p(\xi^{[P]})-F^p(\xi)\|_1
 \le\|F^p(\xi^{[P]})-F^p(\xi)\|_1
    +\|\mathbf1_{S_\infty\setminus A_P}F^p(\xi)\|_1
 \longrightarrow0.
\]
The remaining sum is bounded uniformly in $P$ by
$\sum_{p>M}\delta_pr_p^2\to0$. Letting first $P\to\infty$ and then
$M\to\infty$ proves convergence, in particular of each fixed tangential
frequency.

\subsection{Uniform frequency bounds}
\label{sec:max3-retained-bootstrap}

We use the coefficient and total parameter norms of
Section~\ref{sec:coefficient-definitions}, with the $\ell^1$ metric on
parameter sets. For an integrable diagonal
\[
 V=e+\sum_{j\in S}v_jy_j+
        \sum_{\substack{n>0\\n\notin S}}b_n\rho_n,
 \qquad \rho_n=u_nu_{-n},\qquad c=\lim_{n\to\infty}b_n,
\]
the norm is exactly
\begin{equation}\label{eq:max3-rb-diagonal-norm}
 \|V\|_{\mathcal H,\eta,r}
 =|e|+e^{2\eta}r^2\left[
   \sum_{j\in S}|v_j|+
   \max\left(|c|,\sup_{n\notin S}(1+n)|b_n-c|\right)\right].
\end{equation}
We use this identity to estimate the frequency corrections.

\begin{proposition}[Uniform frequency estimates]
\label{prop:max4-retained-bootstrap}
Let \(S_*\subset\mathbb N\) be infinite with
\(L=\sum_{j\in S_*}j^{-1}<\infty\). Fix a nonempty finite
\(S\subset S_*\), put \(d=|S|\), and fix finitely many nested
sets \(A_1\subset\cdots\subset A_P\subset S\); allow \(P=0\).
For every \(p\), let \(Z^{(p)}_\alpha\) be a real integrable
correction of weighted degree at most four on a nonempty Borel
domain \(E_p\subset\mathbb R^{A_p}\), measured in the ordered
coefficient norm of Proposition~\ref{prop:frequency-interface}
at its original radius \(0<r_p\leq1\) and width \(\eta_p\geq0\).
Set
\begin{equation}\label{eq:max3-rb-retained-budgets}
 \begin{gathered}
 Q_p=\sup_{E_p}\|Z^{(p)}\|_{\eta_p,r_p},\qquad
 D_p=\operatorname{Lip}_{E_p}(Z^{(p)};\|\cdot\|_{\eta_p,r_p}),\\
 K_p=200(1+\max(A_p\cup\{0\}))(10+3L+2|A_p|),\\
 \delta_p=K_p\frac{Q_p+r_p^2D_p}{r_p^4},\qquad
 \beta=\sum_{p=1}^P\delta_p,\qquad
 m_{\mathrm{old}}=\sum_{p=1}^P\delta_pr_p^2.
 \end{gathered}
\end{equation}
Empty sums are zero. Let
\(B=\prod_{j\in S}[b_j/2,b_j]\), with \(b_j>0\), let
\(\ell=\min_{j\in S}b_j/2\), and let \(E_0\subset B\) be
Borel. Assume, for every $\xi\in E_0$ and $1\le p\le P$, that
\begin{equation}\label{eq:max3-rb-prefix-tail}
 \xi|_{A_p}\in E_p,\qquad
 \sum_{j\in S\setminus A_p}\xi_j\leq r_p^2.
\end{equation}
Evaluate the coefficients of $Z^{(p)}$ at $\alpha=\xi|_{A_p}$ and set
\[
 Y_n^{(p)}=\begin{cases}
 \rho_n-\xi_n,&n\in A_p,\\
 \rho_n,&n\notin A_p.
 \end{cases}
\]
Define
\begin{equation}\label{eq:max3-rb-reference}
 \begin{split}
 H_{\mathrm{ref}}(\rho;\xi)
 ={}&\sum_{n\geq1}n^2\rho_n+
       6\left(\sum_{n\geq1}\rho_n\right)^2
       -3\sum_{n\geq1}\rho_n^2\\
   &+\sum_{p=1}^P Z^{(p)}_{\xi|_{A_p}}(Y^{(p)}).
 \end{split}
\end{equation}
Writing $y_j=\rho_j-\xi_j$ for $j\in S$, decompose
\[
 N_{\mathrm{ret}}=\Pi_{\le2}H_{\mathrm{ref}},\qquad
 W_{\mathrm{ret}}=\Pi_4H_{\mathrm{ref}},\qquad
 H_{\mathrm{ref}}=N_{\mathrm{ret}}+W_{\mathrm{ret}}.
\]
Polynomial identities are first read on finitely supported actions.
Frequencies are obtained by differentiating in the actions at fixed
$\xi$ and evaluating at $y=0$, $\rho_{S^c}=0$.

Fix \(0<r\leq1\), \(\epsilon>0\), and a bound \(u\geq0\). Let \(U\) be a real integrable diagonal family on
\(E_0\) with
\(\|U\|^{\mathrm{tot}}_{\mathcal H,\eta_{\mathrm{in}},r}\leq u\),
where \(\eta_{\mathrm{in}}\geq0\), and put
\(N_0=N_{\mathrm{ret}}+U\). Fix \(M\geq0\) with
\(\sum_{j\in S}b_j\leq M\), put \(X=M+2\), and assume
\begin{equation}\label{eq:max3-rb-small-budgets}
 \beta\leq\frac1{48},\qquad
 \frac{(d+1)u}{r^4}\leq\frac1{48},\qquad
 \frac{4(d+1)\epsilon}{r^4}\leq\frac1{24}.
\end{equation}

\smallskip
\noindent\textup{A. Frequency perturbations.}
There are Borel functions $a_0,f_0,g_0$ such that
\[
 \begin{gathered}
 \omega_{j,0}=j^3+ja_0-6j\xi_j+jf_{j,0},\qquad
 \Omega_{n,0}=n^3+na_0+g_{n,0},\\
 N_0=e_0+\sum_{j\in S}\frac{\omega_{j,0}}j y_j
         +\sum_{\substack{n>0\\n\notin S}}\frac{\Omega_{n,0}}n\rho_n.
 \end{gathered}
\]
They satisfy
\begin{equation}\label{eq:max3-rb-seed-frequencies}
 \begin{gathered}
 \operatorname{Lip}(f_0;\ell^1)\leq\frac1{24},\qquad
 \sup_n\operatorname{Lip}g_{n,0}\leq\frac1{24},\\
 \sup_{E_0,n}|g_{n,0}|\leq\frac1{10},\qquad
 \sup_{E_0}(\|\xi\|_1+\|f_0\|_1)\leq X-1.
 \end{gathered}
\end{equation}
For every Borel \(F\subset E_0\), every
\(0\leq\eta\leq\eta_{\mathrm{in}}\), and every real integrable
diagonal \(\Delta\) on \(F\) with
\(\|\Delta\|^{\mathrm{tot}}_{\mathcal H,\eta,r;F}\leq4\epsilon\),
the diagonal \(N_0|_F+\Delta\) has such a representation and obeys
\begin{equation}\label{eq:max3-rb-full-neighborhood}
 \begin{aligned}
 \operatorname{Lip}_F(f;\ell^1)
   &\leq\beta+(d+1)(u+4\epsilon)/r^4\leq1/12,\\
 \sup_n\operatorname{Lip}_F g_n
   &\leq\beta+(u+4\epsilon)/r^4\leq1/12,\\
 \sup_{F,n}|g_n|
   &\leq m_{\mathrm{old}}+(u+4\epsilon)/r^2<1/5,\\
 \sup_F(\|\xi\|_1+\|f\|_1)
   &\leq M+m_{\mathrm{old}}+(d+1)(u+4\epsilon)/r^2<X.
 \end{aligned}
\end{equation}
These bounds are independent of $\min_j b_j$.
The smallness assumptions \eqref{eq:max3-rb-small-budgets} may instead
be replaced here by the corresponding numerical bounds in
\eqref{eq:max3-rb-seed-frequencies}--\eqref{eq:max3-rb-full-neighborhood},
with the four initial quantities bounded respectively by
\[
 \beta+(d+1)u/r^4,\quad \beta+u/r^4,\quad
 m_{\mathrm{old}}+u/r^2,\quad
 M+m_{\mathrm{old}}+(d+1)u/r^2.
\]

\smallskip
\noindent\textup{B. Nonresonance estimates.}
Fix \(\tau>d+8\), \(0<h\leq1/32\), \(\eta_f\geq0\),
\(\eta_{\mathrm{in}}=\eta_f+24h\), and \(0<\kappa\leq1\).
For \(\nu\geq0\), set
\[
 h_\nu=h2^{-\nu},\qquad \kappa_\nu=\kappa2^{-\nu},\qquad
 \eta_{\nu,i}=\eta_f+8h+(16-i)h_\nu\quad(0\leq i\leq8).
\]
There is \(C_D=C_D(S,X,\tau)\geq1\), fixed independently of
\(r,h,\kappa\), the side lengths, parameter subsets, and coefficient
values satisfying the preceding bounds, such that, with
\begin{equation}\label{eq:max3-rb-inverse-constants}
 \begin{gathered}
 T_0=C_D\bigl[\kappa^{-1}h^{-\tau-4}
            +r^2\kappa^{-2}h^{-2\tau-5}\bigr],\\
 T_\nu=T_0\,2^{(2\tau+7)\nu},\qquad
 \lambda_\nu=C_D\kappa2^{-\nu}/\ell,
 \end{gathered}
\end{equation}
the following holds. For every diagonal in the neighborhood of
part A at width \(\eta_{\nu,0}\), on every Borel
\(F\subset E_0\), there is a Borel \(F'\subset F\) with
\begin{equation}\label{eq:max3-rb-original-box-loss}
 \frac{|F\setminus F'|}{|B|}\leq\lambda_\nu,
\end{equation}
on which the three inverses in the triangular iteration of
Proposition~\ref{prop:max3-newton}, from width \(\eta_{\nu,i}\)
to \(\eta_{\nu,i+1}\) for \(i=0,2,4\), have norm at most
\(T_\nu\). The sources have weighted degree at most two and zero integrable
projection; they may be defined only on $F'$. The set $F'$ may depend on
the diagonal and $\nu$, but is common to all sources and widths.
Here $\mu_4=0$, and no amplitude lower bound is needed.

Suppose the other hypotheses of
Proposition~\ref{prop:max3-newton} hold in the coefficient spaces and widths above, with initial jet and
cubic terms bounded by \(\epsilon\), higher-order perturbation bounded
by \(K\geq0\), and Poisson bracket estimates
\(B_\nu=B_0 2^{3\nu}\), \(B_0=C_A r^{-2}h^{-3}\), where
\(C_A\) is a fixed constant in the coefficient estimates. Put
\begin{equation}\label{eq:max3-rb-newton-constants}
 \begin{gathered}
 D_0=B_0T_0,\qquad K_*=2(K+\epsilon),\qquad
 t_0=D_0(K_*+\epsilon),\\
 \Lambda=2^{6(2\tau+10)},\qquad
 C_0=9D_0(1+t_0)^5,\qquad x_0=C_0\Lambda\epsilon.
 \end{gathered}
\end{equation}
Assume \(t_0\leq1/64\), \(x_0\leq1/16\), and
\(\sum_{\nu\geq0}\lambda_\nu<|E_0|/|B|\).
Then the successive diagonals \(N_\nu\) and their limit
\(N_\infty\) satisfy \eqref{eq:max3-rb-full-neighborhood}.
Moreover, every partial diagonal sum formed from all updates
through step \(\nu\), retaining arbitrary degree subsets of
each \(Z_j\) and \(I\Pi_{\leq2}\mathcal E_j\),
\(0\leq j\leq\nu\), lies in the \(4\epsilon\) neighborhood
at width \(\eta_{\nu,0}\), after restriction to \(E_{\nu+1}\).
Here \(Z_j\) and \(\mathcal E_j\) denote the integrable update and
Lie-series remainder in Proposition~\ref{prop:max3-newton}.

For the limiting diagonal, a further exclusion of normalized measure
at most $C_D\kappa/\ell$ gives the inverse of the cubic homological operator from
$\eta_f+8h$ to $\eta_f+7h$, with norm at most $T_0$.
All measure losses are normalized by $|B|$. If
\begin{equation}\label{eq:max3-rb-final-positive-measure}
 3C_D\kappa/\ell<|E_0|/|B|,
\end{equation}
the final normalized measure is at least
\(|E_0|/|B|-3C_D\kappa/\ell>0\).
Cubic normalization uses the bracket and Lie estimates of
Proposition~\ref{prop:max3-concrete-inner} under its hypotheses.

\smallskip
\noindent\textup{C. The integrable correction.}
Suppose that cubic normalization applies with Poisson bracket constant \(B_0\)
and inverse constant \(T_0\), as in
Proposition~\ref{prop:max3-concrete-inner}. Let the initial Hamiltonian be
\(H_0=N_0+J_0+P_0\), and assume the exact identity
\begin{equation}\label{eq:max3-rb-seed-quartic}
 I\Pi_4P_0=W_{\mathrm{ret}}+V_4,
\end{equation}
where \(V_4\) is an integrable homogeneous quartic polynomial, with total norm
at most \(v\) at \(\eta_{\mathrm{in}}\). On the final set define
\begin{equation}\label{eq:max3-rb-new-correction}
 \begin{split}
 Z_{\mathrm{new}}={}&U+\Delta_\infty+V_4
                +I\Pi_4(P_\infty-P_0)\\
                &+\tfrac12 I\{\Pi_3P_\infty,\Gamma\},
 \end{split}
\end{equation}
where \(\Delta_\infty=N_\infty-N_0\) and \(\Gamma\) is the
cubic generator. The final integrable action Hamiltonian is
exactly \(H_{\mathrm{ref}}+Z_{\mathrm{new}}\); the remaining
perturbation has weighted degree at least four and zero integrable
quartic projection. At width \(\eta_f\),
\begin{equation}\label{eq:max3-rb-new-correction-bound}
 \begin{aligned}
 \|Z_{\mathrm{new}}\|^{\mathrm{tot}}
 &\leq u+v+\bigl[3+2t_0(1+t_0)^2+1/8\bigr]\epsilon
                  +\tfrac98 D_0\epsilon^2\\
 &\leq u+v+4\epsilon.
 \end{aligned}
\end{equation}
The homogeneous integrable quartics have zero first action derivatives
at \(y=0\), \(\rho_{\mathrm{normal}}=0\). Their contributions to later
frequencies are given by Proposition~\ref{prop:frequency-interface}.

Set
\begin{equation}\label{eq:max3-rb-future-budget}
 K_S=200(1+\max S)(10+3L+2d),\qquad
 \delta_{\mathrm{new}}=K_S(u+v+4\epsilon)/r^4.
\end{equation}
On every Borel restriction and after every finite or infinite opening
with new action mass at most $r^2$, Proposition~\ref{prop:frequency-interface}
applies to $Z_{\mathrm{new}}$. The resulting corrections satisfy
\[
 \begin{gathered}
 \operatorname{Lip}(f_{\mathrm{new}};\ell^1)
 +\sup_n\operatorname{Lip}g_{\mathrm{new},n}
 +\operatorname{Lip}a_{\mathrm{new}}\le\delta_{\mathrm{new}},\\
 \sup\bigl(\|f_{\mathrm{new}}\|_1+
       \sup_n|g_{\mathrm{new},n}|+|a_{\mathrm{new}}|\bigr)
       \le\delta_{\mathrm{new}}r^2.
 \end{gathered}
\]
The coefficients retain their dependence on the present parameters.
\end{proposition}

\begin{proof}
\emph{Reference frequencies.}
\label{par:max3-rb-retained-frequencies}
Apply Proposition~\ref{prop:frequency-interface} to each
\(Z^{(p)}\), with old set \(A_p\), present set \(S\), and
future actions \(\zeta=\xi|_{S\setminus A_p}\).
Its hypotheses follow from \eqref{eq:max3-rb-prefix-tail}. Since
restriction is contractive and the splitting
$\xi\mapsto(\xi|_{A_p},\xi|_{S\setminus A_p})$ preserves the $\ell^1$
metric, the frequency corrections satisfy
\[
 \begin{gathered}
 \operatorname{Lip}(f^p;\ell^1)+\sup_n\operatorname{Lip}g_n^p
                  +\operatorname{Lip}a_p\leq\delta_p,\\
 \sup\bigl(\|f^p\|_1+\sup_n|g_n^p|+|a_p|\bigr)
                    \leq\delta_pr_p^2.
 \end{gathered}
\]
Summing over the retained corrections yields
\begin{equation}\label{eq:max3-rb-retained-frequency-bounds}
 \begin{gathered}
 \operatorname{Lip}(f_{\mathrm{ret}};\ell^1)\leq\beta,\qquad
 \sup_n\operatorname{Lip}g_{\mathrm{ret},n}\leq\beta,\\
 \sup\bigl(\|f_{\mathrm{ret}}\|_1+
                 \sup_n|g_{\mathrm{ret},n}|\bigr)\leq m_{\mathrm{old}}.
 \end{gathered}
\end{equation}
Since $\partial_{\rho_n}(H_2+Z_4)=n^2+12\sum\rho-6\rho_n$,
the frequencies of $N_{\mathrm{ret}}$ are
\[
 \begin{aligned}
 \omega_{\mathrm{ret},j}
      &=j^3+ja_{\mathrm{ret}}-6j\xi_j+jf_{\mathrm{ret},j},\\
 \Omega_{\mathrm{ret},n}
      &=n^3+na_{\mathrm{ret}}+g_{\mathrm{ret},n},\\
 a_{\mathrm{ret}}&=12\sum_{j\in S}\xi_j+\sum_pa_p.
 \end{aligned}
\]
Recentering is treated as in Section~\ref{sec:frequency-replacement}.
The coefficients are Borel by passage through pointwise action limits.
For $P=0$, all retained corrections vanish.

\medskip\noindent\textit{Frequency shifts.}
\label{par:max3-rb-diagonal-extraction}
The ordered coefficient of $b_n\rho_n$ is $b_n/2$, with first action
limit $c/2$ and block weight $W=2$. Together with the scalar and
$y$-linear blocks, this gives \eqref{eq:max3-rb-diagonal-norm}.

Define the linear maps
\[
 \mathcal A(V)=c,\qquad
 \mathcal F(V)=(v_j-c)_{j\in S},\qquad
 \mathcal G(V)=\bigl(n(b_n-c)\bigr)_{n>0,\ n\notin S}.
\]
Set \(H=\|V\|_{\mathcal H,\eta,r}\). The formula for this norm and
\(\eta\geq0\) give
\begin{equation}\label{eq:max3-rb-extraction-bounds}
 \begin{gathered}
 |\mathcal A(V)|\leq H/r^2,\qquad
 \|\mathcal F(V)\|_1
       \leq\sum_j|v_j|+d|c|\leq(d+1)H/r^2,\\
 \|\mathcal G(V)\|_\infty\leq H/r^2.
 \end{gathered}
\end{equation}
By linearity, the same bounds apply to two-point differences. Thus
$\|V\|^{\mathrm{tot}}\le q$ implies
\[
 \begin{aligned}
 \sup|\mathcal A(V)|,\ \sup\|\mathcal G(V)\|_\infty
      &\le q/r^2,&
 \operatorname{Lip}\mathcal A(V),\ \operatorname{Lip}\mathcal G(V)
      &\le q/r^4,\\
 \sup\|\mathcal F(V)\|_1&\le(d+1)q/r^2,&
 \operatorname{Lip}(\mathcal F(V);\ell^1)&\le(d+1)q/r^4.
 \end{aligned}
\]
The action limits are pointwise limits of Borel coefficients and hence
are Borel.

Adding \(V\) to a diagonal changes its representation by
\begin{equation}\label{eq:max3-rb-frequency-update}
 a\longmapsto a+\mathcal A(V),\qquad
 f\longmapsto f+\mathcal F(V),\qquad
 g\longmapsto g+\mathcal G(V).
\end{equation}
Indeed, $j(\mathcal A(V)+\mathcal F(V)_j)=jv_j$ and
$n\mathcal A(V)+\mathcal G(V)_n=nb_n$; the scalar term contributes zero.

Apply these estimates to $U$ and use
\eqref{eq:max3-rb-retained-frequency-bounds}. Since
$m_{\mathrm{old}}\le\beta$, \eqref{eq:max3-rb-small-budgets} gives
\[
 \begin{aligned}
 \operatorname{Lip}(f_0;\ell^1),\ \sup_n\operatorname{Lip}g_{n,0}
    &\le\beta+(d+1)u/r^4\le1/24,\\
 \sup_{E_0,n}|g_{n,0}|&\le m_{\mathrm{old}}+u/r^2\le1/24<1/10,\\
 \sup(\|\xi\|_1+\|f_0\|_1)
    &\le M+m_{\mathrm{old}}+(d+1)u/r^2\le M+1/24<X-1.
 \end{aligned}
\]
This proves \eqref{eq:max3-rb-seed-frequencies}.

For part A, apply \eqref{eq:max3-rb-extraction-bounds} to
$V=U|_F+\Delta$, whose total norm at width $\eta$ is at most
$u+4\epsilon$. The numerical bounds in
\eqref{eq:max3-rb-full-neighborhood} follow from
\[
 (d+1)(u+4\epsilon)/r^4\le1/16,\qquad
 m_{\mathrm{old}}\le\beta\le1/48,\qquad r\le1,
\]
since $1/16+1/48=1/12<1/5$ and $M+1/12<X$.
The alternative smallness conditions follow from the same estimates.

\medskip\noindent\textit{Small divisors and measure estimates.}
\label{par:max3-rb-inverse}
By \eqref{eq:max3-rb-full-neighborhood}, Proposition~\ref{prop:max4-inverse}
applies with frequency bounds $1/12,1/5,X$. Momentum conservation
cancels the common shift. Moreover,
\[
 2|m|_1+q\le2\quad\Longrightarrow\quad
 |m|_1\le1,\quad q\le2,\quad\mu_4=0,
\]
so every nonintegrable low-order source belongs to the required class.

Impose all admissible nonresonance conditions with threshold
$\kappa_\nu$ and exponent $\tau$. Proposition~\ref{prop:max4-exclusion}
gives a Borel set $F'\subset F$ satisfying
\eqref{eq:max3-rb-original-box-loss}. The set $F'$ is common to the three equations, since coefficient
division is independent of the source and width. Antisymmetry gives
$\{N,F\}=-U$.

With consecutive width loss \(h_\nu\), the inverse bound is
\[
 C_D\bigl[\kappa_\nu^{-1}h_\nu^{-\tau-4}
             +r^2\kappa_\nu^{-2}h_\nu^{-2\tau-5}\bigr].
\]
Since $2^{(\tau+5)\nu}\le2^{(2\tau+7)\nu}$, this proves
\eqref{eq:max3-rb-inverse-constants}, after enlarging $C_D$ to include
the measure constant. Uniformity follows from
Propositions~\ref{prop:max4-exclusion} and~\ref{prop:max4-inverse}.

\medskip\noindent\textit{Convergence of the diagonal terms.}
\label{par:max3-rb-partial-diagonals}
Under the remaining assumptions of B,
Proposition~\ref{prop:max3-newton} applies, since its
inverse assumption has just been verified throughout the required
neighborhood.
Its parameters are \(a=3\), \(b=2\tau+7\), \(c=0\), and
\(d_{\mathrm{Newton}}=2\tau+10\). It gives
\(\|N_\nu-N_0\|^{\mathrm{tot}}_{\eta_{\nu,0}}<3\epsilon\)
and convergence to an increment \(\Delta_\infty\) at
\(\eta_\infty=\eta_f+8h\), with norm less than
\(3\epsilon\). Part A applies at every such width.

For the partial-sum assertion, write the integrable update at
step \(j\) as \(Z_j+I\Pi_{\leq2}\mathcal E_j\), in the
notation of Proposition~\ref{prop:max3-newton}. Its estimates are
\begin{equation}\label{eq:max3-rb-update-sums}
 \sum_{j\geq0}\|Z_j\|^{\mathrm{tot}}_{\eta_{j,4}}
          \leq2(1+t_0)^2\epsilon,\qquad
 \sum_{j\geq0}\|\mathcal E_j\|^{\mathrm{tot}}_{\eta_{j,8}}
          \leq\epsilon/8.
\end{equation}
Every integrable update term has weighted degree zero or two.
Equation~\eqref{eq:max3-rb-diagonal-norm} shows that increasing
the width by \(t\geq0\) multiplies the norm by at most \(e^{2t}\).
This also holds in total norm. For $j<\nu$, both available widths
are at least $\eta_{\nu,0}$; for $j=\nu$, the losses are $4h_\nu$ and
$8h_\nu$. Thus passage to $\eta_{\nu,0}$ costs at most
$e^{16h}\le e^{1/2}$ in \eqref{eq:max3-rb-update-sums}.

By Theorem~\ref{thm:coefficient-initial}, the degree projections are
contractive. Thus every partial sum in the statement has norm at most
\begin{equation}\label{eq:max3-rb-partial-bound}
 e^{1/2}\bigl[2(1+t_0)^2+1/8\bigr]\epsilon<4\epsilon.
\end{equation}
Indeed, \(2(65/64)^2+1/8<9/4\) and \(e^{1/2}<5/3\).
All sums are increments relative to the fixed \(N_0\).

The extraction maps \eqref{eq:max3-rb-extraction-bounds} are bounded
linear maps into $\mathbb R$, $\ell^1(S)$ and $\ell^\infty(\mathcal N_+)$,
with norms at most $(d+1)/r^2$. They therefore commute with the limit
in both the supremum and global Lipschitz norms.

Apply Propositions~\ref{prop:max4-exclusion} and~\ref{prop:max4-inverse}
to the limiting diagonal with threshold $\kappa$, cutoff zero and width
loss $h$. Cubic sources satisfy $|m|_1\le1$, $q\le3$ and have zero
integrable projection. The total measure loss is
\[
 \sum_{\nu\ge0}\frac{C_D\kappa2^{-\nu}}\ell
       +\frac{C_D\kappa}\ell=\frac{3C_D\kappa}\ell,
\]
which proves B.

\medskip\noindent\textit{The quartic correction.}
\label{par:max3-rb-new-quartic}
By Corollary~\ref{cor:max3-cubic} and
\eqref{eq:max3-rb-seed-quartic}, the final integrable part is
\[
 \begin{aligned}
 N_0+\Delta_\infty+I\Pi_4P_\infty
       +\tfrac12 I\{\Pi_3P_\infty,\Gamma\}
 &=N_{\mathrm{ret}}+W_{\mathrm{ret}}+Z_{\mathrm{new}}\\
 &=H_{\mathrm{ref}}+Z_{\mathrm{new}}.
 \end{aligned}
\]
The remaining perturbation $R$ satisfies $\Pi_{\le3}R=I\Pi_4R=0$
by the same corollary.

The three newly generated terms satisfy
\begin{equation}\label{eq:max3-rb-quartic-term-bounds}
 \begin{aligned}
 \|\Delta_\infty\|^{\mathrm{tot}}&<3\epsilon,\\
 \|I\Pi_4(P_\infty-P_0)\|^{\mathrm{tot}}
   &\leq\bigl[2t_0(1+t_0)^2+1/8\bigr]\epsilon,\\
 \bigl\|\tfrac12I\{\Pi_3P_\infty,\Gamma\}\bigr\|^{\mathrm{tot}}
   &\leq\tfrac98 D_0\epsilon^2.
 \end{aligned}
\end{equation}
Adding the bounds for \(U\) and \(V_4\), restricted to the final
parameter set, proves the
first inequality in \eqref{eq:max3-rb-new-correction-bound}.
Since \(K_*+\epsilon=2K+3\epsilon\),
\(D_0\epsilon\leq t_0/3\leq1/192\). Therefore
\[
 2t_0(1+t_0)^2+1/8+\tfrac98 D_0\epsilon
 \leq\tfrac1{32}(65/64)^2+1/8+9/1536<1,
\]
which proves the second inequality.

Since an integrable quartic is quadratic in $(y,\rho_{S^c})$, its first
action derivatives vanish at the current torus. Its later frequency
contributions are retained in $Z_{\mathrm{new}}$.

Let \(Q_{\mathrm{new}}=\sup\|Z_{\mathrm{new}}\|_{\eta_f,r}\)
and \(D_{\mathrm{new}}\) be its global coefficient Lipschitz
seminorm. Equation~\eqref{eq:max3-rb-new-correction-bound} gives
\(Q_{\mathrm{new}}+r^2D_{\mathrm{new}}\leq u+v+4\epsilon\).
Proposition~\ref{prop:frequency-interface}, applied at radius $r$ with
\eqref{eq:max3-rb-future-budget}, proves C for every admissible opening,
including an infinite nonnegative $\ell^1$ tail.
\end{proof}

The finite-stage construction below uses
\(K=r^{4-\delta/4}\), with \(\delta=10^{-3}\).

\subsection{Transformation estimates}
\label{subsec:opening-domains}
We use the norms for the generating Hamiltonians and the parameter estimates of
Section~\ref{subsec:flows-parameters}. Fix $s\ge2$ and $\sigma=s+4$.
Fix disjoint finite sets $S,T$, $T\ne\varnothing$, and put $A=S\cup T$. Fix an upper amplitude bound \(A_0\ge1\), positive width bounds, a finite preliminary cutoff, an old lower amplitude \(0<a_{\rm old}\le A_0\), an old radius \(0<R_{\rm old}\le1\), and a coefficient bound \(K\ge1\). The old parameter set is any nonempty Borel set
$E_{\rm old}\subset[a_{\rm old}^2,A_0^2]^S$.
For $\delta=10^{-3}$ and a new radius $r$, put
\[
 \rho=r^{1-1.01\delta},\qquad a=r^{1-\delta},\qquad
 E\subset E_{\rm old}\times[a^2/2,a^2]^T,
\]
where $E$ is any nonempty Borel set. Write $\beta=(\alpha,\zeta)$.

Consider a preliminary generating Hamiltonian on $E_{\rm old}$,
a sequence of inner generating Hamiltonians on $E$, and a final
generating Hamiltonian on $E$, all real and affine in the actions.
With the weights of Theorem~\ref{thm:physical-inner}, assume
\begin{align}
 P&=\physnorm{\chi_{\rm pre}}^{\mathrm{tot}}_{E_{\rm old},S,\eta_{\rm pre},\rho}
       \le K\rho^4r^{-0.05\delta},\label{eq:physical-prebudget}\\
 Q&=\sum_\nu c_\nu\physnorm{F_\nu}^{\mathrm{tot}}_{E,A,\eta_\nu,r}
       +c_\Gamma\physnorm{\Gamma}^{\mathrm{tot}}_{E,A,\eta_\Gamma,r}
       \le Kr^{4-5.10\delta}.\label{eq:physical-postbudget}
\end{align}
All generators have weighted degree at most five. The realization
constants for the varying post-opening widths and cutoffs are included
in $c_\nu,c_\Gamma$.
Let $w_j=\avg{j}^{2\sigma}$ and $m(\beta)=\sum_{j\in T}w_j\zeta_j$. Define
\[
 D_\beta(r)=\left\{u\in H^\sigma:
 \sum_{j\in A}w_j\bigl||u_j|^2-\beta_j\bigr|\le r^2,
 \quad\physnorm{Q_Au}_\sigma\le r\right\}.
\]

\begin{theorem}[Transformation for one opening]\label{thm:physical-opening}
There exist constants $C_0\ge1$ and $r_*>0$, depending only on the
fixed data and uniform over the admissible old parameters, such that,
for every $0<r<r_*$ satisfying
\eqref{eq:physical-prebudget}--\eqref{eq:physical-postbudget}, the following
map is well defined:
\[
 \tau_\beta=T_{{\rm pre},\alpha}\circ\Phi_{\beta,\infty}\circ T_{\Gamma,\beta}
 \quad\hbox{on }D_\beta(r).
\]
It is jointly continuous into $H^\sigma$ on its
parameter-dependent $H^\sigma$ domain and is Borel measurable. Set
\begin{align*}
 d&=C_0(P/\rho^2+Q/r^2),&
 \ell&=C_0(P/\rho^3+Q/r^3),\\
 A_{\rm path}&=C_0(P+Q),&
 B_{\rm path}&=C_0(P/\rho+Q/r^2).
\end{align*}
Then for $u,w\in D_\beta(r)$,
\begin{align*}
 \physnorm{\tau_\beta u-u}_s&\le d,\\
 e^{-\ell}\physnorm{u-w}_s&\le\physnorm{\tau_\beta u-\tau_\beta w}_s
       \le e^\ell\physnorm{u-w}_s.
\end{align*}
For different parameters and inputs in their respective domains,
\[
 \physnorm{\tau_\beta u-\tau_{\beta'}w}_s
 \le e^\ell\left(\physnorm{u-w}_s+
 C_0(P/\rho^4+Q/r^4)\physnorm{\beta-\beta'}_1\right).
\]
For every finite approximation and every constituent flow trajectory,
the following integral estimates hold:
\begin{align*}
 \sum_{\rm paths}\int\sum_{j\in S}w_j
       \physabs{\frac{d}{dt}|v_j|^2}\,dt&\le A_{\rm path},\\
 \sum_{\rm paths}\int\physnorm{Q_SV(t)}_\sigma\,dt&\le B_{\rm path}.
\end{align*}
Here $V$ denotes the bounded part of the vector field for the corresponding
factor, and $Q_S$ includes all newly opened sites $T$. Consequently,
at every intermediate point and at the limiting endpoints,
\begin{align}
 \physnorm{Q_Av}_\sigma&\le r+B_{\rm path}<3r,\label{eq:physical-currentnormal}\\
 \sum_{j\in A}w_j\bigl||v_j|^2-\beta_j\bigr|
 &\le r^2+A_{\rm path}+2\sqrt{m(\beta)+r^2}\,B_{\rm path}
                   +B_{\rm path}^2<4r^2.\label{eq:physical-currentaction}
\end{align}
The threshold $r_*$ may also be chosen so that
\[
 r^2+A_{\rm path}\le R_{\rm old}^2,
 \qquad\sqrt{m(\beta)+2r^2}+B_{\rm path}\le R_{\rm old}.
\]
Thus all finite flow trajectories and limiting endpoints belong to the
old $H^\sigma$ domain, with weighted deviation of the old actions at
most $R_{\rm old}^2$ and old normal norm at most $R_{\rm old}$. Moreover,
\[
 d\le Cr^{2-6\delta},\quad\ell\le Cr^{1-6\delta},\quad
 A_{\rm path}\le Cr^{4-6\delta},\quad B_{\rm path}\le Cr^{2-6\delta}.
\]
\end{theorem}

\begin{proof}
The post-opening working radius is $4r$. Since every weighted degree is at most five,
\[
 \physnorm{\chi}_{\eta,4r}\le4^5\physnorm{\chi}_{\eta,r},\qquad
 \physnorm{\chi}_{\eta,4r}^{\mathrm{tot}}\le4^7\physnorm{\chi}_{\eta,r}^{\mathrm{tot}}.
\]
The factor $4^2$ comes from the parameter weight. For small $r$,
$4r\le1$, the old amplitudes exceed $\max(10\rho,40r)$, and
$a/\sqrt2\ge40r$. Thus the flow constants are uniform in $r$.

The final flow changes the current actions by at most $CQ$ and the
normal norm by at most $CQ/r$. Choose $r_*$ so that
$CQ<r^2/100$ and $CQ/r<r/100$. Its trajectory then lies in the
half annulus at working radius $4r$, where
Theorem~\ref{thm:physical-inner} applies. Along every finite composition
after opening, the weighted action deviation and normal norm are
bounded by $2r^2$ and $2r$, respectively.

The initial norm of the old normal component is at most $\sqrt{m+2r^2}$. For the bounded part of a vector field after opening,
\[
 \physnorm{Q_SV}_\sigma
 \le\physnorm{\Pi_TV}_\sigma+\physnorm{Q_AV}_\sigma
 \le C_T\physnorm{V}_s+C\epsilon/r\le C\epsilon/r^2.
\]
Here $Q_S=\Pi_T+Q_A$, and the scalar transport is unitary on the old
normal subspace. Before the preliminary flow,
\[
 \sum_{j\in S}w_j\bigl||v_j|^2-\alpha_j\bigr|\le r^2+CQ,
 \qquad
 \physnorm{Q_Sv}_\sigma\le\sqrt{m+2r^2}+CQ/r^2.
\]
Since $a/\rho=r^{0.01\delta}\to0$, these bounds place the endpoint in the
radius-$\rho$ half annulus for small $r$. Apply
Proposition~\ref{prop:flows-finite} to the preliminary flow. Its action
variation and normal displacement are bounded by $CP$ and $CP/\rho$,
respectively, proving the integral estimates.

The map estimates follow by summing displacements and multiplying the
metric bounds, then passing to the strong limit. For the parameter
estimate, apply Lemma~\ref{lem:physical-parameter} and
Theorem~\ref{thm:physical-inner} with coefficient Lipschitz bounds
$P/\rho^2$ and $Q/r^2$. The established domain inclusions also justify
composition of the jointly strongly continuous factors.

To prove the domain estimates in $H^\sigma$, let $b(t)$ be the
integral of the old normal norm of the bounded vector field along a
finite concatenation of flow trajectories, and set
$R_T(t)=\physnorm{\Pi_Tv(t)}_\sigma$. Then
$R_T(t)\le\sqrt{m+r^2}+b(t)$ and
\[
 \sum_{j\in T}w_j\physabs{\frac d{dt}|v_j|^2}
 \le2R_T(t)b'(t).
\]
Integrating gives
\[
 \sum_{\rm paths}\int\sum_{j\in T}w_j
       \physabs{\frac d{dt}|v_j|^2}\,dt
 \le 2\sqrt{m+r^2}B_{\rm path}+B_{\rm path}^2.
\]
Together with the old-action bound, this proves
\eqref{eq:physical-currentaction}. Removing the scalar rotations and
integrating $Q_AV$ gives \eqref{eq:physical-currentnormal}; the estimate
for $Q_Sv$ follows in the same way. Strong convergence extends these
bounds to the limiting endpoints and the preliminary path.

The threshold estimates follow from
\begin{align*}
 P&\le Kr^{4-4.09\delta},& Q&\le Kr^{4-5.10\delta},\\
 P/\rho^2&\le Kr^{2-2.07\delta},&P/\rho^3&\le Kr^{1-1.06\delta},\\
 P/\rho&\le Kr^{3-3.08\delta},& Q/r^3&\le Kr^{1-5.10\delta}.
\end{align*}
Since $T$ is fixed and $1-6.10\delta>0$,
\[
 \frac{B_{\rm path}}r+\frac{A_{\rm path}}{r^2}
 +\frac{\sqrt{m+r^2}B_{\rm path}}{r^2}\longrightarrow0.
\]
Decreasing $r_*$ gives the remaining domain inclusions.
\end{proof}

\medskip\noindent\textit{Application to the inner iteration.}
The generators in Section~\ref{sec:verified-newton} have degree at
most three and fourth-index cutoff zero. Since their analyticity widths have a uniform positive lower bound,
the quantity in \eqref{eq:physical-postbudget} is bounded by a fixed
multiple of the sum of their total norms.

\section{Preparatory normal forms}
\label{sec:preparatory-normal-forms}
We now construct the normal forms needed before opening new modes.
The first result treats the initial stage, when the tangential set is
empty; the second treats an already constructed torus. We normalize the
relevant quartic and quintic terms before opening, and use
Section~\ref{sec:verified-opening} to estimate the resulting low-order
jet. This provides the input for the KAM theorem of
Section~\ref{sec:verified-newton}.

\subsection{At the origin}
We first normalize the quartic terms involving the modes to be opened.
\begin{theorem}[Initial normal form]
\label{thm:max-initial-cutoff}
Let $s\ge2$ and $\sigma=s+4$. With the conventions of
Section~\ref{sec:hamiltonian}, set $\rho_n=u_nu_{-n}$,
$M=\sum_{n>0}\rho_n$, and
\[
H=H_2+H_4+H_6,\quad H_2=\sum_{n>0}n^2\rho_n,
\quad H_4=\tfrac12\langle u^4\rangle,\quad H_6=\tfrac16\langle u^6\rangle,
\quad Z_4=6M^2-3\sum_{n>0}\rho_n^2.
\]
The tangential set is empty. Fix a finite nonempty set
$T\subset\mathbb N$ and put $L=\max T$. For an ordered zero-momentum
quartet $J=(j_1,j_2,j_3,j_4)$ of nonzero integers, set
$\Delta(J)=\sum_{a=1}^4j_a^3$ and define
\[
B_T=\frac12\sum_{\substack{j_1+\cdots+j_4=0,\ \Delta(J)\ne0\\
\{|j_1|,\ldots,|j_4|\}\cap T\ne\varnothing}}u_{j_1}\cdots u_{j_4},
\qquad G_T=H_4-Z_4-B_T,
\]
\[
\chi_T=\sum_{\substack{j_1+\cdots+j_4=0,\ \Delta(J)\ne0\\
\{|j_1|,\ldots,|j_4|\}\cap T\ne\varnothing}}
\frac{u_{j_1}\cdots u_{j_4}}{2i\Delta(J)}.
\]
The sums follow the ordered convention. Every monomial of $G_T$ is
supported outside $T\cup(-T)$.
There exists $\eta_s>0$, depending only on $s$, such that for every
$\eta\ge\eta_s$ there are $0<R_*=R_*(T,s,\eta)\le1$ and
$C=C(T,s,\eta)$ for which the following assertions hold whenever
$0<R\le R_*$.

\begin{enumerate}[label=(\roman*)]
\item The Hamiltonian $\chi_T\in\mathcal C_{G,\eta+1,R}(\varnothing;L)$
is real, even and translation invariant, has vanishing action limits,
and satisfies
\[
\{H_2,\chi_T\}=-B_T,\qquad
\|\chi_T\|_{G,\varnothing,\eta+1,R}
\leq32(1+L)e^{4(\eta+1)}R^4.
\]
No parameter exclusion is required.

\item For $s\le v\le\sigma$, the vector field
$X_{\chi_T}=\partial_x\nabla\chi_T:H^v_0\to H^v_0$ is a bounded
cubic analytic map. Its real-analytic flow $\Phi_T^t$, $|t|\le1$, is
compatible across Sobolev levels and maps $\{\|u\|_\sigma<R\}$ into
$\{\|u\|_\sigma<2R\}$. It preserves $M$ and the weak symplectic form
$\Omega$, and commutes with spatial translations and sign reversal.
For $s\le v\le\sigma$, $|t|\le1$ and
$\|u\|_\sigma,\|w\|_\sigma<R$,
\[
\|\Phi_T^t(u)-u\|_v\leq C_sR^2|t|\,\|u\|_v,
\]
\[
e^{-C_sR^2|t|}\|u-w\|_v
\leq\|\Phi_T^t(u)-\Phi_T^t(w)\|_v
\leq e^{C_sR^2|t|}\|u-w\|_v.
\]
Here $C_s$ and the smallness threshold for these flow properties depend
only on $s$.

\item There exists a real even momentum-preserving Hamiltonian
\(\mathcal R_T\) of normal degree at least six, with formal coefficients
independent of \(R\), such that
\begin{equation}
H\circ\Phi_T^1=H_2+Z_4+G_T+\mathcal R_T,
\qquad \|\mathcal R_T\|_{\varnothing,\eta,R}\leq C R^6.
\label{max-init-A}
\end{equation}
The identity holds coefficientwise and on $\{\|u\|_\sigma<R\}$.
With $\operatorname{ad}_{\chi_T}P=\{P,\chi_T\}$, the remainder is
\begin{equation}
\mathcal R_T=H_6+
\sum_{n\geq1}\left(
\frac{\operatorname{ad}_{\chi_T}^{n}(H_4+H_6)}{n!}
-\frac{\operatorname{ad}_{\chi_T}^{n}B_T}{(n+1)!}
\right),
\label{max-init-B}
\end{equation}
with absolute convergence in $\mathcal H_{\eta,R}(\varnothing)$.

\item Fix \(\delta=10^{-3}\). For a new radius \(0<r<1\), set
\[
\rho=r^{1-1.01\delta},\qquad a=r^{1-\delta},\qquad
\varrho=r^{1-0.99\delta},\qquad
\mathcal B_T=[a^2/2,a^2]^T.
\]
For all sufficiently small $r>0$, depending only on $T,s,\eta$, take
$R=\rho$ and apply
the opening $\mathcal O_T$ with $\xi\in\mathcal B_T$. Set
$A_\xi=\sum_{t\in T}\xi_t$ and
$M_T=\sum_{n>0,n\notin T}u_nu_{-n}$. Then
\begin{equation}
\mathcal O_T(H\circ\Phi_T^1)
=N_\xi+Q_T+G_T+E_\xi,
\label{max-init-C}
\end{equation}
where
\[
\begin{split}
N_\xi={}&c_\xi+
\sum_{t\in T}(t^2+12A_\xi-6\xi_t)y_t
+\sum_{n>0,n\notin T}(n^2+12A_\xi)u_nu_{-n},\\
c_\xi={}&\sum_{t\in T}t^2\xi_t+6A_\xi^2-3\sum_{t\in T}\xi_t^2,\\
Q_T={}&6\left(M_T+\sum_{t\in T}y_t\right)^2
-3\sum_{n>0,n\notin T}(u_nu_{-n})^2-3\sum_{t\in T}y_t^2,\\
E_\xi={}&\mathcal O_T\mathcal R_T.
\end{split}
\]
In the global total norm on \(\mathcal B_T\),
\begin{equation}
\|E\|^{\mathrm{tot}}_{T,\eta,r}
\leq\|E\|^{\mathrm{tot}}_{T,\eta,\varrho}
\leq C\rho^6=C r^{6-6.06\delta},
\quad
\|Q_T+G_T\|^{\mathrm{tot}}_{T,\eta,r}\leq C r^4.
\label{max-init-D}
\end{equation}
Moreover, $\mathcal O_TG_T=G_T$ and $\Pi_{\mathrm{int}}G_T=0$.
Apart from the unbounded quadratic part of $N_\xi$, all families in
\eqref{max-init-C} belong to the indicated coefficient class and are
Borel, real, momentum preserving and even.

Write \(I=\Pi_{\mathrm{int}}\) and set
\[
\widetilde N_\xi=N_\xi+I\Pi_{\leq2}E_\xi,\quad
J_\xi=(1-I)\Pi_{\leq2}E_\xi,\quad
P_\xi=Q_T+G_T+\Pi_{\geq3}E_\xi.
\]
Then $\widetilde N_\xi$ is diagonal,
$\mathcal O_T(H\circ\Phi_T^1)=\widetilde N_\xi+J_\xi+P_\xi$, and
\begin{equation}
\|\widetilde N-N\|^{\mathrm{tot}}_{T,\eta,r}
+\sum_{d=0}^2\|\Pi_dJ\|^{\mathrm{tot}}_{T,\eta,r}
+\|\Pi_3P\|^{\mathrm{tot}}_{T,\eta,r}\leq C\rho^6,
\qquad \|P\|^{\mathrm{tot}}_{T,\eta,r}\leq Cr^4.
\label{max-init-E}
\end{equation}
The jet and cubic part satisfy $|m|_1\le1$, $q\le3$, and are bounded
by $r^{6-6.5\delta}$ for small $r$.

Finally, \eqref{max-init-C} holds as a scalar identity on the physical annulus
\[
U_{\sigma,\xi}(3r)=
\left\{u:\sum_{t\in T}\big||u_t|^2-\xi_t\big|<9r^2,
\quad \|\Pi_{\mathbb N\setminus T}u\|_\sigma<3r\right\},
\]
where the normal projection retains both signs. For small $r$, these annuli lie uniformly in $\{\|u\|_\sigma<\rho\}$,
and all opened actions are positive.
\end{enumerate}
\end{theorem}
\begin{proof}
Throughout the proof, $C=C(T,s,\eta)$ and $C_s=C(s)$, uniformly for
$v\in[s,s+4]$.

\medskip\noindent\textit{The quartic generating function.}

By the quartic resonance calculation in Section~\ref{sec:hamiltonian},
$H_4-Z_4$ has coefficient $1/2$ on nonintegrable zero-momentum quartets
and vanishes otherwise. Such quartets contain no opposite pair. Hence
$B_T,G_T$ have vanishing action limits, block norms at most $1/2$, and
\begin{equation}
\|B_T\|_{\varnothing,\zeta,R},\ \|G_T\|_{\varnothing,\zeta,R}
\leq 2e^{4\zeta}R^4. \label{max-init-proof-1}
\end{equation}
Since $\mu_4(B_T)\le L$, Proposition~\ref{prop:max4-empty-inverse}
with $P=B_T$, $a=g_n=0$ and $Q=H_2$ gives (i).
The symmetries follow from the support and $\Delta(-J)=-\Delta(J)$.

Theorem~\ref{thm:coefficient-initial} also gives
\(\|H_4\|_{\varnothing,\zeta,R}\leq2e^{4\zeta}R^4\) and
\(\|H_6\|_{\varnothing,\zeta,R}=e^{6\zeta}R^6\).

\medskip\noindent\textit{Sobolev estimates for the flow.}

For a nonresonant quartet $(b,c,d,-b-c-d)$, set
$x=b+c$, $y=b+d$, $z=c+d$. These are nonzero integers. If $N$ denotes
the largest absolute entry of the quartet, then
$2N\le|x|+|y|+|z|\le3\max(|x|,|y|,|z|)$, and therefore
\begin{equation}
|\Delta|=3|xyz|\geq3\max(|x|,|y|,|z|)\geq2N. \label{max-init-proof-2}
\end{equation}
For \(n\ne0\), differentiation in the ordered convention gives
\begin{equation}
(X_{\chi_T}(u))_n=
\sum_{b+c+d=n}\kappa_T(n;b,c,d)u_bu_cu_d,
\quad
\kappa_T=\frac{2n}{b^3+c^3+d^3-n^3}
\label{max-init-proof-3}
\end{equation}
when all four entries are nonzero, the denominator is nonzero and the quartet \((-n,b,c,d)\) contains an entry in \(T\cup(-T)\); put \(\kappa_T=0\) otherwise. Equation \eqref{max-init-proof-2} implies \(|\kappa_T|\leq1\), with no dependence on \(T\).

On the complexified Sobolev space, use \(\|u\|_v^2=\frac12\sum_{n\ne0}\langle n\rangle^{2v}|u_n|^2\), and set
\[
A_s=\left(2\sum_{n\ne0}\langle n\rangle^{-2s}\right)^{1/2}<\infty.
\]
Using $\|u\|_{\ell^1}\le A_s\|u\|_s$,
$\langle b+c+d\rangle^v\le3^{v-1}(\langle b\rangle^v+
\langle c\rangle^v+\langle d\rangle^v)$ and Young's inequality, the
trilinear map in \eqref{max-init-proof-3} satisfies
\begin{equation}
\|X_T(u_1,u_2,u_3)\|_v
\leq C_s\sum_{i=1}^3\|u_i\|_v\prod_{j\ne i}\|u_j\|_s.
\label{max-init-proof-4}
\end{equation}
The constant is uniform for $s\le v\le s+4$. In particular,
$X_{\chi_T}(u)=X_T(u,u,u)$ is analytic and
\begin{equation}
\|X_{\chi_T}(u)\|_v\leq C_s\|u\|_s^2\|u\|_v,
\qquad \|DX_{\chi_T}(u)\|_{H^v\to H^v}\leq C_s\|u\|_v^2.
\label{max-init-proof-5}
\end{equation}

The Picard theorem gives a holomorphic flow for $|t|<2$ on a small
$H^v$ ball, uniformly in $T,v$. Uniqueness gives compatibility, and
Gronwall's inequality yields
\[
 \|\Phi_T^tu\|_v\le e^{C_sR^2|t|}\|u\|_v.
\]
For small $R$, this is less than $2R$. Integration and the forward and
backward variational estimates give the displacement and metric bounds.

Bounded coefficients and \(H^s\hookrightarrow\ell^1\) make \(\chi_T\) analytic. By ordered differentiation and continuity,
\(d\chi_T(u)[h]=\Omega(X_{\chi_T}(u),h)\). Differentiation and symmetry of \(d^2\chi_T\) yield
\[
\Omega(DX_{\chi_T}(u)h,k)+\Omega(h,DX_{\chi_T}(u)k)=0.
\]
Thus the flow preserves $\Omega$. Translation and sign equivariance
follow by uniqueness, and $\{M,\chi_T\}=0$ gives conservation of $M$.

\medskip\noindent\textit{The transformed remainder.}

Apply Theorem~\ref{thm:max-coefficient-algebra} with widths $\eta$,
$\eta+1/2$, $\eta+1$, radius $R$ and cutoff $L$. Its Lie constant
satisfies $K\le CR^2$. Choose $R_*$ so that $eK\le1/2$. For
$F\in\mathcal H_{\eta+1/2,R}(\varnothing)$,
\begin{equation}
\sum_{n\geq1}\frac{\|\operatorname{ad}_{\chi_T}^nF\|_{\varnothing,\eta,R}}{n!}
\leq\sum_{n\geq1}(eK)^n\|F\|_{\varnothing,\eta+1/2,R}
\leq2eK\|F\|_{\varnothing,\eta+1/2,R}.
\label{max-init-proof-6}
\end{equation}
Applying \eqref{max-init-proof-6} to $H_4+H_6$ and $B_T$, whose norms
are $O(R^4)$, proves absolute convergence of \eqref{max-init-B} and
$\|\mathcal R_T\|_{\varnothing,\eta,R}\le CR^6$. Each bracket raises
degree by two; hence only even degrees at least six occur, with finitely
many terms at each degree. The coefficients are independent of $R$,
and the stated symmetries are preserved.

For scalar convergence, let $F_q$ be a zero-momentum homogeneous
polynomial of degree $q$ with block seminorm $a_q$. Since
$|p_q(J)|\le3a_q$,
\begin{equation}
|F_q(u)|\leq3a_q\left(\sum_{n\ne0}|u_n|\right)^q.
\label{max-init-proof-7}
\end{equation}
Choose $e^{\eta_s}>8(1+A_s)$. Since $2A_se^{-\eta}<1/4$, the
majorant \eqref{max-init-proof-7} and all its differentiated series
converge on $\{\|u\|_s<2R\}$. Thus coefficient evaluation is continuous
into scalar analytic functions on smaller balls.

Ordered differentiation on finitely supported vectors, followed by
density using \eqref{max-init-proof-5}--\eqref{max-init-proof-7}, gives
\[
 dF_q(u)[X_{\chi_T}(u)]=\{F_q,\chi_T\}(u),\qquad
 dH_2(u)[X_{\chi_T}(u)]=-B_T(u).
\]
Taylor expansion along the analytic flow now gives
\eqref{max-init-A}--\eqref{max-init-B}, since
$H_4-B_T=Z_4+G_T$.

\medskip\noindent\textit{The jet after opening.}

The hypotheses of Proposition~\ref{prop:max3-opening} hold for small
$r$, since
\begin{equation}\label{eq:initial-opening-ratios}
 \begin{aligned}
 t_*&=2e^{2\eta}\frac{\varrho^2}{a^2}
       =2e^{2\eta}r^{0.02\delta}\longrightarrow0,\\
 \vartheta&=\frac{\varrho}{\rho}
       +\frac{2|T|a}{\rho\sqrt{1-t_*}}
       =r^{0.02\delta}+\frac{2|T|r^{0.01\delta}}{\sqrt{1-t_*}}
       \longrightarrow0.
 \end{aligned}
\end{equation}
Taking $r$ smaller so that $\rho\le R_*$ and using parameter
independence of $\mathcal R_T$, we obtain
\begin{equation}
\|\mathcal O_T\mathcal R_T\|^{\mathrm{tot}}_{T,\eta,\varrho}
\leq C(1+L)\|\mathcal R_T\|_{\varnothing,\eta,\rho}
\leq C\rho^6. \label{max-init-proof-8}
\end{equation}
Monotonicity of the total norm in the radius gives the first bound in \eqref{max-init-D}.

Substituting
\[
 \rho_t=\xi_t+y_t\quad(t\in T),\qquad
 M=A_\xi+M_T+\sum_{t\in T}y_t
\]
in $H_2+Z_4$ gives $c_\xi,N_\xi,Q_T$. The bounded part
$N_\xi-\mathcal O_TH_2$ consists of constants, tangential actions and
the normal mass, hence belongs to the coefficient class.

Opening leaves \(G_T\) unchanged, with zero action limits and norm at most \(2e^{4\eta}r^4\). Corollary~\ref{cor:max3-opening-support} gives \(\|Q_T\|^{\mathrm{tot}}_{T,\eta,r}\leq Cr^4\), with zero derivative with respect to the new parameters, proving \eqref{max-init-D}. The coefficients are Lipschitz by \eqref{max-init-proof-8} or polynomial in \(\xi\); the symmetries follow from Proposition~\ref{prop:max3-opening}.

Theorem~\ref{thm:coefficient-initial} gives
\[
 \max\left\{
 \|I\Pi_{\le2}E\|^{\mathrm{tot}},\,
 \max_{0\le d\le2}\|\Pi_d(1-I)E\|^{\mathrm{tot}},\,
 \|\Pi_3E\|^{\mathrm{tot}}\right\}\le\|E\|^{\mathrm{tot}}.
\]

Summing these five bounds proves \eqref{max-init-E}. Integrable
terms of degree at most two are diagonal, while $2|m|_1+q\le3$ implies
$|m|_1\le1$ and $q\le3$. The required initial bounds follow from
\[
 \rho^6=r^{6-6.5\delta}r^{0.44\delta}\le r^4.
\]

It remains to verify that both sides are defined on the same Sobolev domain. Put \(W_T=\sum_{t\in T}\langle t\rangle^{2\sigma}\) and \(w_T=\max_{t\in T}\langle t\rangle^{2\sigma}\). For \(u\in U_{\sigma,\xi}(3r)\),
\[
\|u\|_\sigma^2
\leq W_T a^2+9(w_T+1)r^2.
\]
Consequently
\[
 \rho^{-2}\|u\|_\sigma^2
 \le W_Tr^{0.02\delta}+9(w_T+1)r^{2.02\delta}\longrightarrow0.
\]
This proves uniform inclusion in the original ball. Moreover,
$|y_t|/\xi_t\le18r^{2\delta}\to0$, so all actions remain positive and
the square-root expansions converge.

The opened factors and normal $\ell^1$ norm have majorant
$C_Ta+3A_sr=o(\rho)$. Thus \eqref{max-init-proof-7} gives absolute
convergence after substitution, and Proposition~\ref{prop:max3-opening}
justifies the action limits. Equation~\eqref{max-init-C} follows from
\eqref{max-init-A}.\qedhere
\end{proof}

\subsection{Near an invariant torus}
We next treat a nonempty tangential set.

\begin{theorem}[Preliminary normal form]
\label{thm:max-old-preparation}
Let \(s\ge2\), \(\sigma=s+4\), and let \(S,T\subset\mathbb N\) be finite,
nonempty and disjoint. Put \(S'=S\cup T\), \(L=\max S'\), and let \(Q_S\)
be the projection onto the normal modes outside \(S\), including both signs. Let
\[
 \mathcal B=\prod_{j\in S}[\ell_j,u_j],\qquad
 0<\ell_j<u_j,\qquad \ell_*=\min_{j\in S}(u_j-\ell_j),
\]
and let \(E\subset\mathcal B\) be Borel with positive Lebesgue measure.
Fix \(0<R_0\le1\) with \(10R_0\le\min_j\sqrt{\ell_j}\), and
\[
 C_s=4\left(1+\sum_{n\ge1}(1+n^2)^{-s}\right)^{1/2},
 \qquad \eta\ge\eta_*:=\sigma+3+\log(8C_s),\quad 0<h\le1.
\]
These data are fixed before choosing the new radius.

Suppose \(A_0,c:E\to\mathbb R\) are bounded Lipschitz functions, and
\(f:E\to\mathbb R^S\), \(g_n:E\to\mathbb R\), \(n>0\), \(n\notin S\), satisfy
\[
 \operatorname{Lip}_{\ell^1}f\le\epsilon,\qquad
 \sup_n\operatorname{Lip}g_n\le\epsilon,\qquad
 \sup_{\alpha,n}|g_n(\alpha)|\le B_0,\qquad
 \sup_\alpha(\|\alpha\|_1+\|f(\alpha)\|_1)\le X_0,
\]
where \(0\le\epsilon\le1/12\), \(0\le B_0\le1/5\), \(X_0<\infty\).
Define
\[
 \omega_j=j^3+jA_0-6j\alpha_j+jf_j,\qquad
 \Omega_n=n^3+nA_0+g_n,\qquad
 N_\alpha=c+\sum_{j\in S}\frac{\omega_j}{j}y_j
              +\sum_{n>0,n\notin S}\frac{\Omega_n}{n}u_nu_{-n}.
\]
Extend \(\Omega_n\) oddly. Let \(Q_4,R\) be real, momentum-preserving,
even coefficient families on \(E\) with finite total norms at
\((S,\eta+3h,R_0)\). Assume
\[
 IQ_4=Q_4=\Pi_4Q_4,\qquad
 \Pi_{\le3}R=0,\qquad I\Pi_4R=0,\qquad H=N+Q_4+R.
\]

Fix \(\delta=10^{-3}\), \(\tau>|S|+8\), and for \(0<r<1\) put
\[
 \rho=r^{1-1.01\delta},\quad a=r^{1-\delta},\quad
 \varrho=r^{1-0.99\delta},\quad
 \kappa=\min(1,\ell_*)r^{0.01\delta},\quad
 \mathcal B_T=[a^2/2,a^2]^T.
\]
For an old monomial of degree \(d\), let \(h_T\) count its normal factors
at sites in \(\pm T\), including multiplicities, and set
\[
 B_{\rm pre}=\mathbf1_{\{h_T\ge d-3\}}\Pi_{4,5}R,\quad
 G_{\rm pre}=\mathbf1_{\{h_T<d-3\}}\Pi_{4,5}R,\quad R_6=\Pi_{\ge6}R.
\]
Let \(E_r\subset E\) be the set defined by the inequalities
\[
 |k\cdot\omega+\textstyle\sum_b\Omega_{n_b}|
       \ge\kappa(1+|k|_1)^{-\tau}
\]
for all nonintegrable, even, momentum-preserving tuples with at most five
normal entries and fourth normal index at most \(L\).
There exist $r_*>0$ and $C<\infty$, depending only on these fixed data
and bounds, such that the following assertions hold for $0<r\le r_*$.

The set \(E_r\) is Borel, has positive measure, and
\[
 |E\setminus E_r|\le Cr^{0.01\delta}|\mathcal B|.
\]
The coefficient formula
\[
 (\chi_{\rm pre})_{k,m,q}(\mathbf n)
 =\frac{(B_{\rm pre})_{k,m,q}(\mathbf n)}
        {i(k\cdot\omega+\sum_b\Omega_{n_b})}
\]
on the support of $B_{\rm pre}$, and zero elsewhere, defines a real
Hamiltonian $\chi_{\rm pre}\in\mathcal C_G^{\mathrm{aff}}(S;L)$ of
weighted degrees four and five, satisfying
\begin{equation}
 \{N,\chi_{\rm pre}\}=-B_{\rm pre},\qquad
 \|B_{\rm pre}\|^{\rm tot}_{\eta+3h,\rho}
 +\|G_{\rm pre}\|^{\rm tot}_{\eta+3h,\rho}\le C\rho^4,\qquad
 \|\chi_{\rm pre}\|^{\rm tot}_{G,\eta+5h/2,\rho}
       \le C\rho^4/\kappa.\label{eq:preparation-p1}
\end{equation}
The first two norms are on \(E\); the generator norm is on \(E_r\).

For each \(\alpha\in E_r\), define
\[
 U_{v,\alpha}(b\rho)=
 \left\{u:\sum_{j\in S}\big||u_j|^2-\alpha_j\big|<b^2\rho^2,\
                 \|Q_Su\|_v<b\rho\right\},\qquad0<b\le3.
\]
The Hamiltonian flow \(\Phi_\alpha^t\), \(|t|\le1\), exists on
\(U_{\sigma,\alpha}(\rho)\), stays in \(U_{\sigma,\alpha}(2\rho)\), and
has the regularity and weak symplectic variational operator of Proposition~\ref{prop:flows-finite}. It satisfies
\begin{equation}
 \|\Phi_\alpha^t(u)-u\|_s\le C|t|\rho^2/\kappa,\qquad
 e^{-C|t|\rho/\kappa}\|u-w\|_v
 \le\|\Phi_\alpha^t(u)-\Phi_\alpha^t(w)\|_v
 \le e^{C|t|\rho/\kappa}\|u-w\|_v,\label{eq:preparation-p2}
\end{equation}
for \(u,w\in H^\sigma\) in that domain and \(s\le v\le\sigma-1\).
In the decomposition
\(X_{\chi_{\rm pre}}=V_\alpha+b_\alpha\partial_xQ_Su\), the trajectories satisfy
\begin{equation}
 \int_{-1}^1\sum_{j\in S}\langle j\rangle^{2\sigma}
 |2\operatorname{Re}(\overline{u_j}(V_\alpha)_j)|\,dt
 \le C\rho^4/\kappa,\qquad
 \int_{-1}^1\|Q_SV_\alpha\|_\sigma\,dt\le C\rho^3/\kappa.\label{eq:preparation-p3}
\end{equation}

There is a real, even, momentum-preserving remainder $R^\sharp$ with
$\Pi_{\le5}R^\sharp=0$ such that, coefficientwise and on the old
$H^\sigma$ annulus,
\begin{equation}
 H_\alpha\circ\Phi_\alpha^1
   =N_\alpha+(Q_4)_\alpha+(G_{\rm pre})_\alpha+(R^\sharp)_\alpha,
 \qquad
 \|R^\sharp\|^{\rm tot}_{S,\eta+h,\rho;E_r}\le C\rho^6/\kappa.\label{eq:preparation-p4}
\end{equation}
Writing \(F=Q_4+R\) and \(\operatorname{ad}_\chi U=\{U,\chi_{\rm pre}\}\),
the remainder is the absolutely convergent series
\begin{equation}
 R^\sharp=R_6+\sum_{n\ge1}
 \left(\frac{\operatorname{ad}_\chi^nF}{n!}
       -\frac{\operatorname{ad}_\chi^nB_{\rm pre}}{(n+1)!}\right).\label{eq:preparation-p5}
\end{equation}

Open the sites \(T\) at \(\zeta\in\mathcal B_T\), and put
\[
 \bar N=\Pi_{\le2}\mathcal O_T(N+Q_4),\quad
 \bar Q=\Pi_4\mathcal O_TQ_4,\quad
 \bar G=\mathcal O_TG_{\rm pre},\quad
 \bar F=\mathcal O_TR^\sharp.
\]
Then
\begin{equation}
 \mathcal O_T(H\circ\Phi^1)=\bar N+\bar Q+\bar G+\bar F,\label{eq:preparation-p6}
\end{equation}
where \(\bar N\) is diagonal, \(\bar Q\) is an integrable homogeneous quartic,
and \(\Pi_{\le3}\bar G=I\bar G=0\). On the whole product \(E_r\times\mathcal B_T\),
\begin{equation}
 \|\bar Q\|^{\rm tot}_{S',\eta,r}\le Cr^4,\qquad
 \|\bar G\|^{\rm tot}_{S',\eta,r}\le Cr^{4-0.08\delta},\qquad
 \|\bar F\|^{\rm tot}_{S',\eta,r}
 \le\|\bar F\|^{\rm tot}_{S',\eta,\varrho}
 \le C\rho^6/\kappa\le Cr^{6-6.07\delta}.\label{eq:preparation-p7}
\end{equation}
The unbounded part of \(\bar N\) is obtained by exact substitution.
If
\[
 \widetilde N=\bar N+I\Pi_{\le2}\bar F,\qquad
 J=(1-I)\Pi_{\le2}\bar F,\qquad
 P=\bar Q+\bar G+\Pi_{\ge3}\bar F,
\]
then the decomposition \(\widetilde N+J+P\) is exact, \(\widetilde N\)
is diagonal, and at width \(\eta\), radius \(r\),
\begin{equation}
 \|\widetilde N-\bar N\|^{\rm tot}
 +\sum_{d=0}^2\|\Pi_dJ\|^{\rm tot}+\|\Pi_3P\|^{\rm tot}
 \le Cr^{6-6.07\delta},\qquad
 \|P\|^{\rm tot}\le Cr^{4-0.08\delta}.\label{eq:preparation-p8}
\end{equation}
For small $r$, the jet and cubic norms are bounded by
$r^{6-6.5\delta}$, and $\|P\|^{\rm tot}\le r^{4-\delta/4}$.
The jet and cubic part satisfy $|m|_1\le1$ and $q\le3$.

These equalities also hold as scalar identities on the new annuli
\[
 \left\{u:
 \sum_{j\in S}\big||u_j|^2-\alpha_j\big|
 +\sum_{t\in T}\big||u_t|^2-\zeta_t\big|<9r^2,\
 \|Q_{S'}u\|_\sigma<3r\right\}.
\]
They lie uniformly inside \(U_{\sigma,\alpha}(\rho)\).
All resulting families of finite coefficient norm are Borel, real,
momentum preserving and even.
The width loss is \(3h\). Proposition~\ref{prop:coherent-production} verifies the frequency hypotheses needed to apply the inner iteration to this decomposition.
\end{theorem}

\begin{proof}
\medskip\noindent\textit{The generating function.}
For $\Pi_{\ge d}U=U$ and $0<R\le R_0$, homogeneity gives
\begin{equation}
 \|U\|^{\rm tot}_{S,\zeta,R}
       \le (R/R_0)^d\|U\|^{\rm tot}_{S,\zeta,R_0}.\label{eq:preparation-p9}
\end{equation}
The Lipschitz term gains the additional factor $R^2/R_0^2\le1$.
By Corollary~\ref{cor:max3-opening-support},
\[
 \|Q_4\|^{\rm tot}_{\eta+3h,\rho}
 +\|B_{\rm pre}\|^{\rm tot}_{\eta+3h,\rho}
 +\|G_{\rm pre}\|^{\rm tot}_{\eta+3h,\rho}\le C\rho^4,\qquad
 \|R_6\|^{\rm tot}_{\eta+3h,\rho}\le C\rho^6.
\]
The same corollary gives $\mu_4(B_{\rm pre})\le L$ and
$|m|_1\le1$, with all symmetries preserved.

Propositions~\ref{prop:max4-exclusion} and~\ref{prop:max4-inverse},
applied with $\mu=L$, give
\[
 |E\setminus E_r|\le C\kappa|\mathcal B|/\ell_*
                  \le Cr^{0.01\delta}|\mathcal B|.
\]
Since $|E|>0$ is fixed, $|E_r|>0$ for small $r$.
With width loss $h/2$, the inverse estimate reads
\[
 \|\chi_{\rm pre}\|^{\rm tot}_{G,\eta+5h/2,\rho}
 \le C(\kappa^{-1}+\rho^2\kappa^{-2})\rho^4
 \le C\rho^4/\kappa.
\]
Here \(\rho^2/\kappa=\min(1,\ell_*)^{-1}r^{2-2.03\delta}\to0\).
The fixed powers of \(h^{-1}\) are absorbed in \(C\). The common term
\(A_0\) cancels by zero momentum, and division preserves degree, proving
\eqref{eq:preparation-p1}.

\medskip\noindent\textit{The unbounded quadratic Hamiltonian.}
In the old coordinates, the quadratic Hamiltonian is
\[
 H_2^S=\sum_{j\in S}j^2(\alpha_j+y_j)
              +\sum_{n>0,n\notin S}n^2u_nu_{-n}.
\]
Writing \(N-H_2^S=D_0+D_2\), one has
\[
 D_0=c-\sum_{j\in S}j^2\alpha_j,\qquad
 D_2=\sum_{j\in S}(A_0-6\alpha_j+f_j)y_j+
          \sum_{n>0,n\notin S}(A_0+g_n/n)u_nu_{-n}.
\]
The normal coefficient has the decomposition
\[
 \tfrac12(A_0+g_n/n)=A_0/2+g_n/(2n),\qquad
 (1+n)|g_n|/(2n)\le B_0,
\]
with the analogous two-point estimate. Since the tangential coefficients
form a fixed finite set of bounded Lipschitz functions,
\[
 \|D_2\|^{\rm tot}_{\eta+3h,\rho}\le C\rho^2,
\]
while \(D_0\) is a bounded scalar with zero bracket. Applying the bracket
estimate at width \(\eta+3h/2\), with both losses \(h/2\), proves
\begin{equation}
 Q_0:=\{H_2^S,\chi_{\rm pre}\}
 =-B_{\rm pre}-\{D_2,\chi_{\rm pre}\}
 \in\mathcal H_{\eta+3h/2,\rho}(S),\qquad
 \|Q_0\|^{\rm tot}\le C\rho^4/\kappa.\label{eq:preparation-p10}
\end{equation}
Hence Proposition~\ref{prop:flows-lie} applies to $H_2^S$.

\medskip\noindent\textit{The transformed Hamiltonian.}
Apply Proposition~\ref{prop:flows-affine-realization} with generator
size $C\rho^4/\kappa$. Since $\rho/\kappa\to0$,
Proposition~\ref{prop:flows-finite} gives the flow, weak symplecticity
and \eqref{eq:preparation-p2}. The tangential action derivative and
bounded normal field are $O(\rho^4/\kappa)$ and
$O(\rho^3/\kappa)$, respectively; integration gives
\eqref{eq:preparation-p3}.

For the coefficient Lie estimate at base width \(\eta+h\), increment
\(h/2\), the constant in the Lie estimate satisfies \(K\le C\rho^2/\kappa\).
For a Hamiltonian \(U\) at width \(\eta+3h/2\), Theorem~\ref{thm:max-coefficient-algebra} gives
\[
 \sum_{n\ge1}\frac{\|\operatorname{ad}_\chi^nU\|^{\rm tot}_{\eta+h,\rho}}{n!}
 \le\sum_{n\ge1}(eK)^n\|U\|^{\rm tot}_{\eta+3h/2,\rho}
 \le2eK\|U\|^{\rm tot}_{\eta+3h/2,\rho}
\]
when \(eK\le1/2\). Apply this to \(F=Q_4+R\) and \(B_{\rm pre}\),
whose norms are \(O(\rho^4)\). Formula \eqref{eq:preparation-p5} then has norm
\(C\rho^6+C\rho^6/\kappa\le C\rho^6/\kappa\). Every nonzero
bracket raises weighted degree by at least two; hence it has degree at
least six. Poisson brackets and limits in the coefficient norm preserve all stated symmetries.

Since all widths exceed $\eta_*$,
Propositions~\ref{prop:flows-evaluation} and~\ref{prop:flows-lie},
together with \eqref{eq:preparation-p10}, identify this Lie series
with the Hamiltonian pullback. Now
$F-B_{\rm pre}=Q_4+G_{\rm pre}+R_6$ gives
\eqref{eq:preparation-p4}--\eqref{eq:preparation-p5}.

\medskip\noindent\textit{The jet after opening.}
The ratios in \eqref{eq:initial-opening-ratios} again tend to zero.
Proposition~\ref{prop:max3-opening} therefore applies from
\(\eta+h,\rho\) to \(\eta,\varrho\), in the total norm on the product parameter set. It gives \(\|\bar F\|^{\rm tot}_{\eta,\varrho}\le C\rho^6/\kappa\)
and \(\|\bar G\|^{\rm tot}_{\eta,\varrho}\le C\rho^4\).
Corollary~\ref{cor:max3-opening-support} gives \(\Pi_{\le3}\bar G=I\bar G=0\).
Restricting \(\bar G\) to radius \(r\) and using \eqref{eq:preparation-p9} gives
\[
 \|\bar G\|^{\rm tot}_{\eta,r}
       \le C(r/\varrho)^4\rho^4=Cr^{4-0.08\delta}.
\]
For \(\bar F\), radius monotonicity and
\(\rho^6/\kappa=\min(1,\ell_*)^{-1}r^{6-6.07\delta}\) give the
remaining bound. Corollary~\ref{cor:max3-opening-support} yields
\(\|\bar Q\|^{\rm tot}_{\eta,r}\le Cr^4\): its hypothesis at radius one
follows from the fixed bounds at radius \(R_0\), with factors
\(R_0^{-4}\) and \(R_0^{-6}\) in the supremum and weighted Lipschitz
parts, respectively.

The degree formula in Corollary~\ref{cor:max3-opening-support} gives
the diagonal splitting \eqref{eq:preparation-p6}. The projection bounds
used in \eqref{max-init-E} give \eqref{eq:preparation-p8}, with total
bound $5\|\bar F\|^{\rm tot}$. Finally,
\[
 Cr^{6-6.07\delta}=r^{6-6.5\delta}Cr^{0.43\delta},\qquad
 Cr^{4-0.08\delta}=r^{4-\delta/4}Cr^{0.17\delta},
\]
and $2|m|_1+q\le3$ implies $|m|_1\le1$, $q\le3$.
These are the required initial estimates and support conditions.

\medskip\noindent\textit{Domains and parameter estimates.}
The sum of the absolute deviations of the old tangential actions on the new
annulus is less than \(9r^2<\rho^2\).
Put \(W_T=\sum_{t\in T}\langle t\rangle^{2\sigma}\) and
\(w_T=\max_{t\in T}\langle t\rangle^{2\sigma}\). The old normal
component includes the newly opened sites and satisfies
\[
 \|Q_Su\|_\sigma^2
 \le W_Ta^2+9(w_T+1)r^2,\qquad
 \frac{W_Ta^2+9(w_T+1)r^2}{\rho^2}
 =W_Tr^{0.02\delta}+9(w_T+1)r^{2.02\delta}\longrightarrow0.
\]
This proves uniform inclusion in the old domain. Since \(9r^2/(a^2/2)\to0\),
the new actions are positive and their square-root expansions converge.

The scalar identity follows as in Theorem~\ref{thm:max-initial-cutoff},
using $a/\rho,r/\rho,r^2/\rho^2\to0$ in the substituted majorant.
We omit the identical convergence argument.

Borel dependence follows from the global Lipschitz estimates and the
countable divisor conditions. Adjoining the full new box preserves the
normalized measure loss. The new diagonal frequencies are estimated in
Proposition~\ref{prop:max4-retained-bootstrap}.
\end{proof}

\section{Successive openings}
\label{sec:loop}
We construct normal forms inductively on increasing tangential sets. One step consists of the preliminary
normal form of Section~\ref{sec:preparatory-normal-forms}, the opening
estimates of Section~\ref{sec:verified-opening}, and the KAM and cubic
normalizations of Section~\ref{sec:verified-newton}. The resulting
normal form satisfies the hypotheses needed for the next step.
Following \cite{BGR25,BP26}, we distinguish this outer iteration from
the KAM iteration on a fixed tangential set.

\subsection{One step of the loop}
\label{sec:repeatable-opening}

Fix \(s\geq2\), \(\sigma=s+4\), and \(\delta=10^{-3}\). Let
\(S_*\subset\mathbb N\) be infinite with
$L_*:=\sum_{j\in S_*}j^{-1}<\infty$.
We retain the Fourier and Sobolev conventions of
Section~\ref{sec:verified-flows}.
Put \(\rho_n=u_nu_{-n}\), \(M=\sum_{n>0}\rho_n\), and
\begin{equation}\label{eq:max3-prod-fixed-hamiltonian}
 \begin{gathered}
 H_{\mathrm{fix}}=H_2+H_4+H_6,\qquad
 H_2=\sum_{n>0}n^2\rho_n,\qquad
 H_4=\tfrac12\langle u^4\rangle,\qquad
 H_6=\tfrac16\langle u^6\rangle,\\
 Z_4^{\mathrm{base}}=6M^2-3\sum_{n>0}\rho_n^2.
 \end{gathered}
\end{equation}
Here \(\langle u^k\rangle\) is the spatial average. The Gardner
vector field is
\[
 F(u)=-u_{xxx}+6u^2u_x+\partial_x(u^5).
\]
We use the norms of Proposition~\ref{prop:max3-concrete-inner} and
Theorem~\ref{thm:max-old-preparation}, with parameter weight $R^2$
at radius $R$.

For a finite set \(A\subset\mathbb N\), let \(Q_A\) be the Fourier
projection onto the modes whose absolute indices lie outside \(A\). Define
\begin{equation}\label{eq:max3-prod-domain}
 D_{A,\alpha}(b)=
 \left\{u\in H^\sigma_0:
 \sum_{j\in A}\langle j\rangle^{2\sigma}
       \big||u_j|^2-\alpha_j\big|\leq b^2,
 \quad \|Q_Au\|_\sigma\leq b\right\}.
\end{equation}
Set
\[
 C_s=4\left(1+\sum_{n\geq1}\langle n\rangle^{-2s}\right)^{1/2},
 \qquad \eta_*=\sigma+3+\log(8C_s).
\]

\begin{proposition}[Induction step]
\label{prop:repeatable-opening}
Fix the following data before choosing the new radius.
\begin{enumerate}
\item
Let \(S,T\subset S_*\) be disjoint finite nonempty sets and
\(A=S\cup T\). Fix a box
\(B_S=\prod_{j\in S}[b_j/2,b_j]\), with all \(b_j>0\), and a
Borel set \(E\subset B_S\) such that \(|E|/|B_S|\geq m_*>0\).
Fix \(0<R\leq1\) with
\(10R\leq\min_{j\in S}\sqrt{b_j/2}\).

\item
Fix nested nonempty sets \(A_1\subset\cdots\subset A_m=S\),
radii \(r_q\in(0,1]\), widths \(\eta_q\geq0\), Borel domains
\(E_q\), and integrable families \(Z^{(q)}\) of weighted degree at
most four on \(E_q\). Assume the coefficient hypotheses of
Proposition~\ref{prop:frequency-interface}, with each family depending
only on its original parameters. Set
\[
 K_q=200(1+\max A_q)(10+3L_*+2|A_q|),\qquad
 d_q=K_q\frac{\|Z^{(q)}\|^{\mathrm{tot}}_{A_q,\eta_q,r_q;E_q}}{r_q^4}.
\]
Assume \(\sum_{q=1}^m d_q\leq1/48\), and
\(\alpha|_{A_q}\in E_q\) for every \(\alpha\in E\). Fix positive
tail margins \(t_q\) with
\begin{equation}\label{eq:max3-prod-old-tail}
 \sum_{j\in S\setminus A_q}b_j\leq r_q^2-t_q.
\end{equation}
The resulting integrable Hamiltonian in the original actions is
\begin{equation}\label{eq:max3-prod-inherited}
 \begin{split}
 H_{\mathrm{ret}}(\rho;\alpha)
 &=H_2+Z_4^{\mathrm{base}}
   +\sum_{q=1}^m Z^{(q)}_{\alpha|_{A_q}}(Y^{(q)}),\\
 Y_j^{(q)}&=
 \begin{cases}
 \rho_j-\alpha_j,&j\in A_q,\\
 \rho_j,&j\notin A_q.
 \end{cases}
 \end{split}
\end{equation}
\item
In the current coordinates $y_j=\rho_j-\alpha_j$, $j\in S$, set
$N_S=\Pi_{\le2}H_{\mathrm{ret}}$ and $W_S=\Pi_4H_{\mathrm{ret}}$.
Suppose the preceding normal form is
\begin{equation}\label{eq:max3-prod-old-normal-form}
 H_{\mathrm{old}}=N_S+W_S+R_{\mathrm{old}},\qquad
 \Pi_{\leq3}R_{\mathrm{old}}=0,\qquad I\Pi_4R_{\mathrm{old}}=0.
\end{equation}
Fix \(0<h\leq1/32\) and \(\eta_f\geq\eta_*\), with
\(\eta_{\mathrm{old}}=\eta_f+27h\leq\min_q\eta_q\).
The total norms of \(W_S\) and \(R_{\mathrm{old}}\) at
\((S,\eta_{\mathrm{old}},R;E)\) have fixed finite upper bounds.
The constant coefficient of $N_S$ has fixed supremum and global
Lipschitz bounds, inherited from the preceding finite sum of corrections.

\item\label{hyp:max3-prod-old-realization}
There is a change of variables \(\Psi_\alpha\) in the original Fourier coordinates,
locally defined at every point of \(D_{S,\alpha}(R)\), jointly
Borel in parameter and state, and continuous in the strong
\(H^\sigma\) topology. On an \(H^\sigma\)-open neighborhood of
\(D_{S,\alpha}(R/8)\),
\begin{equation}\label{eq:max3-prod-old-scalar}
 H_{\mathrm{fix}}\circ\Psi_\alpha=H_{\mathrm{old},\alpha}.
\end{equation}
It has the regularity of finite compositions of the actual maps
in Propositions~\ref{prop:flows-finite} and~\ref{prop:flows-inner}:
on each local domain it is locally Lipschitz in \(H^v\) for
\(s\leq v\leq\sigma-1\); in directions belonging to \(H^\sigma\), its directional
derivative agrees with an invertible bounded variational operator
\(E_\Psi\) on each such \(H^v\), preserving
\(\Omega(h,k)=\langle\partial_x^{-1}h,k\rangle\).
This operator is strongly continuous on each fixed vector when
its base point varies in \(H^\sigma\). In particular, for every path in its domain,
\[
 w\in C(H^\sigma)\cap C^1(H^{s+1})
 \quad\Longrightarrow\quad
 \Psi_\alpha\circ w\in C(H^\sigma)\cap C^1(H^{s+1}),\qquad
 (\Psi_\alpha\circ w)'=E_\Psi(w)w'.
\] A prescribed range for the preceding transformation is also included
in the data.
\end{enumerate}
There exist $C\ge1$ and $r_*>0$, depending only on these data and
bounds, such that the following holds for $0<r<r_*$. They are uniform
in the parameters and in Borel sets with the same lower measure bound.
Set
\begin{equation}\label{eq:max3-prod-scales}
 \begin{gathered}
 \rho=r^{1-1.01\delta},\qquad a=r^{1-\delta},\qquad
 \varrho=r^{1-0.99\delta},\qquad B_T(r)=[a^2/2,a^2]^T,\\
 \epsilon=r^{6-6.5\delta},\qquad K=r^{4-\delta/4}.
 \end{gathered}
\end{equation}
There are a Borel set \(E^+\subset E\times B_T(r)\), a single new
integrable correction \(Z^+\) of weighted degree at most four,
a remainder \(R^+\), and maps \(\tau_\beta\) in the original Fourier coordinates, for
\(\beta=(\alpha,\zeta)\in E^+\), with the following properties.

\begin{enumerate}
\item[(a)] The excluded parameter set satisfies
\begin{equation}\label{eq:max3-prod-loss}
 |(E\times B_T(r))\setminus E^+|
 \leq Cr^{0.01\delta}|B_S\times B_T(r)|,\qquad
 \frac{|E^+|}{|B_S\times B_T(r)|}>0.
\end{equation}
\item[(b)] In the new partial action--angle coordinates, the Hamiltonian has the normal form
\begin{equation}\label{eq:max3-prod-new-normal-form}
 H^+=H_2+Z_4^{\mathrm{base}}
       +\sum_{q=1}^m Z^{(q)}_{\alpha|_{A_q}}(Y^{(q)})+Z^+_\beta+R^+_\beta.
\end{equation}
The previously constructed coefficient families are unchanged. At width \(\eta_f\) and radius \(r\),
\begin{equation}\label{eq:max3-prod-new-norms}
 \begin{gathered}
 \|Z^+\|^{\mathrm{tot}}\leq6\epsilon\leq r^{4+1/3},\qquad
 \|R^+\|^{\mathrm{tot}}\leq4K+7\epsilon\leq r^{4-2\delta},\\
 \Pi_{\leq3}R^+=0,\qquad I\Pi_4R^+=0.
 \end{gathered}
\end{equation}
For \(K_A=200(1+\max A)(10+3L_*+2|A|)\), the contribution
of the new correction to subsequent frequency estimates satisfies
\begin{equation}\label{eq:max3-prod-new-budget}
 d^+:=K_A\frac{\|Z^+\|^{\mathrm{tot}}}{r^4}
 \leq6K_A r^{2-6.5\delta}.
\end{equation}
This bound holds for subsequent openings with new action mass at most
$r^2$, and tends to zero with $r$.

\item[(c)] The map \(\tau_\beta\) is locally defined at every
point of \(D_{A,\beta}(r)\), is jointly Borel and strongly
\(H^\sigma\)-continuous, and has the regularity specified in
assumption~\ref{hyp:max3-prod-old-realization}. For
\(u,v\in D_{A,\beta}(r)\),
\begin{equation}\label{eq:max3-prod-metric}
 \begin{gathered}
 \|\tau_\beta u-u\|_s\leq Cr^{2-6\delta},\\
 e^{-Cr^{1-6\delta}}\|u-v\|_s
 \leq\|\tau_\beta u-\tau_\beta v\|_s
 \leq e^{Cr^{1-6\delta}}\|u-v\|_s.
 \end{gathered}
\end{equation}
The integrated estimates along the flow paths are
\begin{equation}\label{eq:max3-prod-path-budgets}
 A_{\mathrm{path}}\leq Cr^{4-6\delta},\qquad
 B_{\mathrm{path}}\leq Cr^{2-6\delta},
\end{equation}
for the absolute variation of the previous tangential actions with
weights \(\langle j\rangle^{2\sigma}\), and the \(H^\sigma\) norm of the bounded part of the
vector field on the entire previous normal subspace \(Q_S\), including
\(T\). These estimates hold for every finite approximant and every constituent
flow, and yield the limiting endpoint bounds
\begin{equation}\label{eq:max3-prod-path-endpoints}
 \begin{gathered}
 \sum_{j\in S}\langle j\rangle^{2\sigma}
       \big||\tau_\beta(u)_j|^2-|u_j|^2\big|\leq A_{\mathrm{path}},\\
 \|Q_S\tau_\beta(u)\|_\sigma\leq\|Q_Su\|_\sigma+B_{\mathrm{path}}.
 \end{gathered}
\end{equation}
By decreasing \(r_*\), we may also impose
\begin{equation}\label{eq:max3-prod-strict-reserve}
 \begin{gathered}
 r^2+A_{\mathrm{path}}<(R/8)^2,\\
 \sqrt{a^2\sum_{t\in T}\langle t\rangle^{2\sigma}+2r^2}
      +B_{\mathrm{path}}<R/8.
 \end{gathered}
\end{equation}
Consequently all intermediate points and their limits belong to
$D_{S,\alpha}(R/8)$.

\item[(d)] On an \(H^\sigma\)-open neighborhood of \(D_{A,\beta}(r/8)\),
\begin{equation}\label{eq:max3-prod-scalar-conjugacy}
 H_{\mathrm{old},\alpha}\circ\tau_\beta=H^+_\beta,\qquad
 H_{\mathrm{fix}}\circ(\Psi_\alpha\circ\tau_\beta)=H^+_\beta.
\end{equation}
The composition has the regularity properties in
assumption~\ref{hyp:max3-prod-old-realization} and the prescribed
range of the preceding transformation. Let \(\omega^+_j\) be \(j\) times the
coefficient of the linear tangential action \(y_j\) in \(H^+\).
Define \(z_\beta(\theta)_j=\sqrt{\beta_j}e^{i\theta_j}\) on
\(A\), zero on the remaining positive modes, with conjugate
negative modes. For every angle \(\theta\),
\begin{equation}\label{eq:max3-prod-finite-solution}
 t\longmapsto(\Psi_\alpha\circ\tau_\beta)
      z_\beta(\theta+t\omega^+(\beta))
\end{equation}
is a global solution of the fixed equation and belongs to
\(C(\mathbb R;H^\sigma_0)\cap C^1(\mathbb R;H^{s+1}_0)\).
\end{enumerate}
\end{proposition}

\begin{proof}
\emph{The preliminary normal form.}
\label{par:max3-prod-preparation}
By Section~\ref{sec:frequency-replacement} and the tail bounds
\eqref{eq:max3-prod-old-tail}, the current frequencies have the form
\[
 \omega_j=j^3+ja_S-6j\alpha_j+jf_j,\qquad
 \Omega_n=n^3+na_S+g_n.
\]
Writing $\widetilde a_S=a_S-12\sum_{j\in S}\alpha_j$, we have
\[
 \operatorname{Lip}(f;\ell^1)+\sup_n\operatorname{Lip}g_n
       +\operatorname{Lip}\widetilde a_S\le\sum_qd_q\le1/48,
\]
\[
 \sup\bigl(\|f\|_1+\sup_n|g_n|+|\widetilde a_S|\bigr)
       \le\sum_qd_qr_q^2\le1/48.
\]
Thus Theorem~\ref{thm:max-old-preparation} applies; the common
coefficient and constant term have the required bounded Lipschitz norms.

Use Theorem~\ref{thm:max-old-preparation} with
\(\eta=\eta_f+24h\), so that its input width is
\(\eta+3h=\eta_{\mathrm{old}}\). Set
\[
 \ell_S=\min_{j\in S}b_j/2,\qquad
 \kappa_{\mathrm{pre}}=\min(1,\ell_S)r^{0.01\delta}.
\]
It yields a Borel set \(E_r\subset E\), with
\begin{equation}\label{eq:max3-prod-preliminary-loss}
 \frac{|E\setminus E_r|}{|B_S|}
       \leq C_{\mathrm{pre}}r^{0.01\delta},
\end{equation}
a preliminary Hamiltonian time-one map
\(\mathcal T_{\mathrm{pre},\alpha}\), and, after opening, the identity
\begin{equation}\label{eq:max3-prod-opened-identity}
 H_{\mathrm{old},\alpha}\circ\mathcal T_{\mathrm{pre},\alpha}
       =N_{\mathrm{bar}}+W_{\mathrm{bar}}+G_{\mathrm{bar}}+F_{\mathrm{bar}}.
\end{equation}
The right-hand side is written in the opened coordinates. The identity
holds on the $H^\sigma$ annulus with unweighted action deviation less
than $9r^2$ and normal norm less than $3r$, which is contained in the
preliminary domain.

The constants depend on $E$ only through $m_*$. Choose $r$ so that
$C_{\mathrm{pre}}r^{0.01\delta}<m_*/2$.

On the entire product \(E_r\times B_T(r)\), preparation gives
\begin{equation}\label{eq:max3-prod-prepared-bounds}
 \begin{aligned}
 \|F_{\mathrm{bar}}\|^{\mathrm{tot}}_{A,\eta,r}
      &\leq C_{\mathrm{seed}}r^{6-6.07\delta},\\
 \|W_{\mathrm{bar}}\|^{\mathrm{tot}}_{A,\eta,r}
      &\leq C_{\mathrm{seed}}r^4,\\
 \|G_{\mathrm{bar}}\|^{\mathrm{tot}}_{A,\eta,r}
      &\leq C_{\mathrm{seed}}r^{4-0.08\delta}.
 \end{aligned}
\end{equation}
It also gives \(IG_{\mathrm{bar}}=0\) and
\(\Pi_{\leq3}G_{\mathrm{bar}}=0\). Define
\begin{equation}\label{eq:max3-prod-seed-split}
 \begin{aligned}
 U&=I\Pi_{\leq2}F_{\mathrm{bar}},&N_0&=N_{\mathrm{bar}}+U,\\
 J_0&=(1-I)\Pi_{\leq2}F_{\mathrm{bar}},&
 V_4&=I\Pi_4F_{\mathrm{bar}},\\
 P_0&=W_{\mathrm{bar}}+G_{\mathrm{bar}}+\Pi_{\geq3}F_{\mathrm{bar}}.
 \end{aligned}
\end{equation}
The projection bounds of the preparation theorem imply, after
decreasing \(r_*\),
\begin{equation}\label{eq:max3-prod-seed-bounds}
 \begin{gathered}
 \sum_{d=0}^2\|\Pi_dJ_0\|^{\mathrm{tot}}\leq\epsilon,\qquad
 \|\Pi_3P_0\|^{\mathrm{tot}}\leq\epsilon,\qquad
 \|P_0\|^{\mathrm{tot}}\leq K,\\
 \|U\|^{\mathrm{tot}}\leq\epsilon,\qquad
 \|V_4\|^{\mathrm{tot}}\leq\epsilon.
 \end{gathered}
\end{equation}
The gains in the exponents are $0.43\delta$ and $0.17\delta$.
Writing $u,v$ for the bounds on $U,V_4$, we retain
\[
 u,v\le C_{\mathrm{seed}}r^{6-6.07\delta}.
\]
The jet of degree at most three satisfies $|m|_1\le1$, $q\le3$.

In the notation of Proposition~\ref{prop:max4-retained-bootstrap},
$N_{\mathrm{bar}}=N_{\mathrm{ret}}$ and
$W_{\mathrm{bar}}=W_{\mathrm{ret}}$, so
\begin{equation}\label{eq:max3-prod-quartic-seed}
 I\Pi_4P_0=W_{\mathrm{ret}}+V_4.
\end{equation}
Impose
\(|T|a^2\leq\tfrac12\min_qt_q\), so that every old tail
inequality remains valid with a positive margin after opening.

\medskip\noindent\textit{KAM and measure estimates.}
\label{par:max3-prod-inner-constants}
Fix
\[
 \tau_A=|A|+10,\qquad
 A_{\mathrm{amp}}=\max(1,\max_{j\in S}\sqrt{b_j}),\qquad
 M_0=\sum_{j\in S}b_j+|T|,\qquad X=M_0+2.
\]
For $r<1$, these bounds determine $C_A,C_D,C_P$ in
Proposition~\ref{prop:max3-concrete-inner}, independently of $r$.
For sufficiently small $r$, the shortest side of the new box is
$\ell(r)=a^2/2=r^{2-2\delta}/2\le\ell_S$. Choose
\begin{equation}\label{eq:max3-prod-kappa}
 \kappa=\ell(r)r^{\delta/2}=\tfrac12 r^{2-1.5\delta}.
\end{equation}
The sum of the preceding correction bounds is at most \(1/48\).
The remaining initial bounds required by Proposition~\ref{prop:max4-retained-bootstrap}
are
\[
 (|A|+1)u/r^4\leq1/48,\qquad
 4(|A|+1)\epsilon/r^4\leq1/24.
\]
These follow from \eqref{eq:max3-prod-seed-bounds} and
$\epsilon/r^4=r^{2-6.5\delta}\to0$. Part A then gives the frequency
bounds throughout the required diagonal neighborhood, on every Borel
restriction.

Substitution into the constants of
Proposition~\ref{prop:max3-concrete-inner} gives
\begin{equation}\label{eq:max3-prod-inner-powers}
 \begin{aligned}
 B_0&=C_Ah^{-3}r^{-2},\\
 T_0&=C_D\bigl(\kappa^{-1}h^{-\tau_A-4}
          +r^2\kappa^{-2}h^{-2\tau_A-5}\bigr)
          \leq C_T r^{-2+1.5\delta},\\
 D_0&=B_0T_0\leq C_D' r^{-4+1.5\delta},\\
 t_0&=D_0(2K+3\epsilon)\leq C_t r^{1.25\delta},\\
 x_0&=9D_0(1+t_0)^5\,2^{6(2\tau_A+10)}\epsilon
          \leq C_x r^{2-5\delta},\\
 T_0\epsilon/r^3&\leq C_T r^{1-5\delta},\qquad
 3C_D\kappa/\ell(r)=3C_Dr^{\delta/2}.
 \end{aligned}
\end{equation}
Here $r^2\kappa^{-2}=4r^{-2+3\delta}\le4r^{-2+1.5\delta}$.
We impose $t_0\le1/64$; all constants in
\eqref{eq:max3-prod-inner-powers} are then independent of $r$.

The positive exponents in \eqref{eq:max3-prod-inner-powers} allow us to
impose all smallness conditions of
Proposition~\ref{prop:max3-concrete-inner}, together with
\begin{equation}\label{eq:max3-prod-loss-positivity}
 C_{\mathrm{pre}}r^{0.01\delta}+3C_Dr^{\delta/2}<m_*.
\end{equation}
We also require
\[
 r^{-\delta}/\sqrt2\ge40,\qquad
 40r\le\min_{j\in S}\sqrt{b_j/2},\qquad r\le1/4,
\]
which give its amplitude and radius conditions.

Apply that proposition to \(N_0+J_0+P_0\) on
\(E_0=E_r\times B_T(r)\), at input width \(\eta_f+24h\).
It gives the Newton generators, the limiting diagonal, the cubic
generator and $E^+$. Together with \eqref{eq:max3-prod-preliminary-loss},
\[
 \frac{|(E\times B_T(r))\setminus E^+|}{|B_S\times B_T(r)|}
 \le C_{\mathrm{pre}}r^{0.01\delta}+3C_Dr^{\delta/2}
 \le Cr^{0.01\delta}.
\]
This and \eqref{eq:max3-prod-loss-positivity} prove
\eqref{eq:max3-prod-loss}.

\medskip\noindent\textit{The new integrable correction.}
\label{par:max3-prod-correction}
The inner normal form is \(N_0+Z_{\mathrm{in}}+R^+\), with
\[
 Z_{\mathrm{in}}=\Delta_\infty+I\Pi_4P_\infty
                  +\tfrac12 I\{\Pi_3P_\infty,\Gamma\}.
\]
Define
\begin{equation}\label{eq:max3-prod-correction-subtraction}
 Z^+=U+Z_{\mathrm{in}}-W_{\mathrm{ret}}.
\end{equation}
Using \eqref{eq:max3-prod-quartic-seed}, this is
\begin{equation}\label{eq:max3-prod-correction-expansion}
 Z^+=U+\Delta_\infty+V_4+I\Pi_4(P_\infty-P_0)
                    +\tfrac12 I\{\Pi_3P_\infty,\Gamma\}.
\end{equation}
Part C of Proposition~\ref{prop:max4-retained-bootstrap} gives
\(\|Z^+\|^{\mathrm{tot}}\leq u+v+4\epsilon\leq6\epsilon\)
and the estimate for subsequent openings~\eqref{eq:max3-prod-new-budget}.
Proposition~\ref{prop:max3-concrete-inner} also gives
\[
 \|R^+\|^{\mathrm{tot}}\le4K+7\epsilon,\qquad
 \Pi_{\le3}R^+=I\Pi_4R^+=0.
\]
Equations~\eqref{eq:max3-prod-opened-identity},
\eqref{eq:max3-prod-seed-split}, \eqref{eq:max3-prod-quartic-seed},
and~\eqref{eq:max3-prod-correction-subtraction} prove
\eqref{eq:max3-prod-new-normal-form}. The preceding corrections keep their original parameter dependence.

Finally,
\[
 \frac{6\epsilon}{r^{4+1/3}}=6r^{5/3-6.5\delta}\longrightarrow0,
 \qquad
 \frac{4K+7\epsilon}{r^{4-2\delta}}\longrightarrow0.
\]
Hence \eqref{eq:max3-prod-new-norms} holds for sufficiently small $r$;
likewise, $d^+\le6K_A r^{2-6.5\delta}\to0$.

\medskip\noindent\textit{Parameter estimates.}
\label{par:max3-prod-physical}
At radius \(\rho\), the preparation generator obeys
\begin{equation}\label{eq:max3-prod-preliminary-generator}
 P:=\|\chi_{\mathrm{pre}}\|^{\mathrm{tot}}_{G,S,\eta+5h/2,\rho}
 \leq C\rho^4/\kappa_{\mathrm{pre}}
 \leq K_{\mathrm{phys}}\rho^4r^{-0.05\delta}.
\end{equation}
The constant $K_{\mathrm{phys}}$ absorbs $\min(1,\ell_S)^{-1}$.
The subsequent generators have degree at most three, are affine in
$y$, and have fourth-index cutoff zero and width at least $\eta_f$.
Their generating norms dominate those in
Theorem~\ref{thm:physical-opening}. Its constants, including second
state derivatives and parameter differences, are bounded by
$c_{\mathrm{phys}}=c_{\mathrm{phys}}(A,A_{\mathrm{amp}},s,\eta_f)$.
The estimate for the sum of the generating Hamiltonians in Proposition~\ref{prop:max3-concrete-inner} and
\eqref{eq:max3-prod-inner-powers} give
\begin{equation}\label{eq:max3-prod-inner-generators}
 Q:=c_{\mathrm{phys}}\left(
       \sum_{\nu\geq0}\|F_\nu\|^{\mathrm{tot}}_G
       +\|\Gamma\|^{\mathrm{tot}}_G\right)
 \leq4c_{\mathrm{phys}}T_0\epsilon
 \leq K_{\mathrm{phys}}r^{4-5.10\delta}.
\end{equation}
Enlarge $K_{\mathrm{phys}}$ to cover both estimates.

Apply the parameter-dependent construction of
Theorem~\ref{thm:physical-opening} with preceding parameter set \(E_r\), new
admissible set \(E^+\), preliminary cutoff \(\max A\), and radius \(R/8\) in the preceding coordinates. Its constants are uniform in the Borel sets $E_r,E^+$. We obtain
\begin{equation}\label{eq:max3-prod-actual-map}
 \tau_\beta=\mathcal T_{\mathrm{pre},\alpha}\circ
   \left(\lim_{N\to\infty}\mathcal T_{0,\beta}\circ\cdots
                                \circ\mathcal T_{N,\beta}\right)
   \circ\mathcal T_{\Gamma,\beta}
\end{equation}
on \(D_{A,\beta}(r)\), with pointwise convergence in \(H^\sigma\) and joint continuity
in the parameter and the \(H^\sigma\) topology. The change to action--angle variables leaves the phase point $u$ unchanged. The theorem gives
\eqref{eq:max3-prod-metric}--\eqref{eq:max3-prod-strict-reserve}
after decreasing $r_*$. The projection $Q_S$ includes $T$.

Before substituting the powers of \(r\), the relevant bounds are
\begin{equation}\label{eq:max3-prod-physical-budgets}
 \begin{aligned}
 d&=C_{\mathrm{dom}}(P/\rho^2+Q/r^2),&
 e&=C_{\mathrm{dom}}(P/\rho^3+Q/r^3),\\
 A_{\mathrm{path}}&=C_{\mathrm{dom}}(P+Q),&
 B_{\mathrm{path}}&=C_{\mathrm{dom}}(P/\rho+Q/r^2).
 \end{aligned}
\end{equation}
Here $C_{\mathrm{dom}}$ depends only on the fixed data; $d,e$ bound
the displacement and logarithmic metric distortion.
Substituting \eqref{eq:max3-prod-preliminary-generator} and
\eqref{eq:max3-prod-inner-generators} into
\eqref{eq:max3-prod-physical-budgets}, with $\rho=r^{1-1.01\delta}$, gives
\[
 \begin{aligned}
 d&\le Cr^{2-6\delta},& e&\le Cr^{1-6\delta},\\
 A_{\mathrm{path}}&\le Cr^{4-6\delta},&
 B_{\mathrm{path}}&\le Cr^{2-6\delta},
 \end{aligned}
\]
which proves \eqref{eq:max3-prod-metric}
and~\eqref{eq:max3-prod-path-budgets}.
The flows after opening are defined using radius \(4r\). For generators of
degree at most five,
\begin{equation}\label{eq:max3-prod-generator-rescaling}
 \|\chi\|_{G,\eta,4r}\leq4^5\|\chi\|_{G,\eta,r},\qquad
 \|\chi\|^{\mathrm{tot}}_{G,\eta,4r}
       \leq4^7\|\chi\|^{\mathrm{tot}}_{G,\eta,r}.
\end{equation}
The factor $4^7$ includes the parameter weight. The amplitude bounds
and $4r\le1$ have already been imposed.

The working annuli contain $H^\sigma$-open neighborhoods of the closed
input domain. Theorem~\ref{thm:physical-opening}, with
Propositions~\ref{prop:flows-finite} and~\ref{prop:flows-inner}, therefore
gives the variational, symplectic and metric estimates in
$H^v$, $s\le v\le\sigma-1$, using
\eqref{eq:max3-prod-generator-rescaling}.
Uniqueness identifies these flows with those of
Proposition~\ref{prop:max3-concrete-inner}; strong convergence
identifies their limiting products.

\medskip\noindent\textit{Conjugacy on the common domain.}
\label{par:max3-prod-chain-rule}
The Hamiltonian identity in Proposition~\ref{prop:max3-concrete-inner}
holds on
\begin{equation}\label{eq:max3-prod-inner-scalar-domain}
 \mathcal U_\sigma(r/4;\beta):=
 \left\{u:\sum_{j\in A}\big||u_j|^2-\beta_j\big|<r^2/16,
                    \quad\|Q_Au\|_\sigma<r/4\right\}.
\end{equation}
All intermediate inner points lie in the radius-$r$ annulus, hence in
the radius-$3r$ domain of \eqref{eq:max3-prod-opened-identity}.
Composition gives $H_{\mathrm{old},\alpha}\circ\tau_\beta=H^+_\beta$
on \eqref{eq:max3-prod-inner-scalar-domain}. Since the action weights
are at least one,
\[
 D_{A,\beta}(r/8)\subset\mathcal U_\sigma(r/4;\beta),\qquad
 \tau_\beta(D_{A,\beta}(r))
   \subset\operatorname{int}_{H^\sigma}D_{S,\alpha}(R/8),
\]
where the second inclusion follows from
\eqref{eq:max3-prod-strict-reserve}. Strong continuity then permits
composition with \eqref{eq:max3-prod-old-scalar} on an $H^\sigma$-open
neighborhood of $D_{A,\beta}(r/8)$, proving
\eqref{eq:max3-prod-scalar-conjugacy}.

Let $T$ be a finite flow or a limiting composition on an
$H^\sigma$-open domain. By Propositions~\ref{prop:flows-finite}
and~\ref{prop:flows-inner}, its variational extension $E_T$ is locally
uniformly bounded on $H^v$, $s\le v\le\sigma-1$. The
$C^1(H^\sigma;H^{\sigma-1})$ regularity and density give
\[
 u_n\to u\text{ in }H^\sigma
 \quad\Longrightarrow\quad
 E_T(u_n)v\to E_T(u)v\text{ in }H^{s+1}
 \quad(v\in H^{s+1}).
\]
For $w\in C(H^\sigma)\cap C^1(H^{s+1})$ in this domain, the segment
from $w(t)$ to $w(t+a)$ remains in the domain for small $a$. The
fundamental theorem of calculus gives
\begin{equation}\label{eq:max3-prod-weak-chain-rule}
 \frac{T(w(t+a))-T(w(t))}{a}
 =\int_0^1 E_T\bigl(w(t)+\lambda(w(t+a)-w(t))\bigr)
             \frac{w(t+a)-w(t)}{a}\,d\lambda.
\end{equation}
The difference quotients converge in $H^{s+1}$. Local uniform
boundedness and strong continuity of $E_T$ allow $a\to0$, giving
\[
 (T\circ w)'(t)=E_T(w(t))w'(t)\quad\text{in }H^{s+1}.
\]
The right-hand side is continuous. The same argument in $H^v$ gives
the remaining regularity in
assumption~\ref{hyp:max3-prod-old-realization}. These properties,
bounded invertibility and symplecticity persist under finite composition.

Differentiating \eqref{eq:max3-prod-scalar-conjugacy} in a direction
$h\in H^\sigma$ is justified by the preceding chain rule. As in
Proposition~\ref{prop:max3-concrete-inner}, weak symplecticity gives
\begin{equation}\label{eq:max3-prod-symplectic-pairing}
 \begin{aligned}
 &\Omega\bigl(F((\Psi_\alpha\circ\tau_\beta)(u)),
                    E_{\Psi\circ\tau}(u)h\bigr)\\
 &\qquad=\Omega(X_{H^+}(u),h)\\
 &\qquad=\Omega\bigl(E_{\Psi\circ\tau}(u)X_{H^+}(u),
                    E_{\Psi\circ\tau}(u)h\bigr).
 \end{aligned}
\end{equation}
The diagonal representation and Proposition~\ref{prop:flows-conjugacy}
give vector fields in $H^{s+1}$. Since the variational operators are
invertible on $H^{s+1}$, density and nondegeneracy of $\Omega$ yield
\begin{equation}\label{eq:max3-prod-vector-conjugacy}
 F((\Psi_\alpha\circ\tau_\beta)(u))
     =E_{\Psi\circ\tau}(u)X_{H^+}(u)
     \quad\text{in }H^{s+1}.
\end{equation}

On the reference torus, the first state derivatives of $R^+$ and of
the homogeneous integrable quartic vanish. Therefore
\[
 X_{H^+}(z_\beta(\theta))
     =D_\theta z_\beta(\theta)\omega^+(\beta).
\]
The reference torus is smooth and lies in the interior of the domain.
Equations~\eqref{eq:max3-prod-weak-chain-rule}--\eqref{eq:max3-prod-vector-conjugacy}
therefore give the global solution \eqref{eq:max3-prod-finite-solution},
with the stated regularity and range.
\end{proof}

\subsection{Iteration of the loop}
\label{subsec:iteration-of-loop}
We initialize the construction with an empty tangential set and then
apply Proposition~\ref{prop:repeatable-opening} recursively. At each
stage, all previous data and the losses of analyticity are fixed before
the new radius is chosen.

\begin{proposition}[Finite-dimensional tori]
\label{prop:coherent-production}
For every increasing enumeration
\(S_*=\{j_1<j_2<\cdots\}\), every \(0<r_0\leq1\), and every
\(0<\vartheta_*\leq1/16\), there exists a compatible sequence of normal forms on
$S_p=\{j_1,\ldots,j_p\}$. Set
\begin{equation}\label{eq:max3-prod-width-sequence}
 \eta_0=\eta_*+2,\qquad h_p=2^{-p}/54,\qquad
 \eta_p=\eta_0-\sum_{q=1}^p27h_q\geq\eta_*+\tfrac32.
\end{equation}
There are radii \(r_p>0\), intervals
\(I_p=[r_p^{2-2\delta}/2,r_p^{2-2\delta}]\), and Borel sets
\(O_p\subset B^{(p)}:=\prod_{q\leq p}I_q\) such that
\begin{equation}\label{eq:max3-prod-recursive-measure}
 O_p\subset O_{p-1}\times I_p,\qquad
 \frac{|O_p|}{|B^{(p)}|}\geq1-\vartheta_*,\qquad
 d_p^{\mathrm{ret}}\leq\frac{2^{-p}}{96}.
\end{equation}
Here \(O_0=B^{(0)}\) is the one-point parameter space and
\(d_p^{\mathrm{ret}}\) is defined by \eqref{eq:max3-prod-new-budget}.
The first stage is supplied by
Theorem~\ref{thm:max-initial-cutoff} and
Proposition~\ref{prop:max3-concrete-inner}.

Each $Z^{(q)}$ and $\tau_q$ retains its original parameter dependence.
At step $p$, the Hamiltonian is
\[
 H^{(p)}=H_2+Z_4^{\mathrm{base}}+\sum_{q\leq p}Z^{(q)}+R^{(p)},
\]
expressed in the current action coordinates. The new correction and
remainder satisfy \eqref{eq:max3-prod-new-norms} at
\((\eta_p,r_p)\). The composition
\[
 \Psi_{p,\alpha}=
 \tau_{1,\alpha|_{S_1}}\circ\cdots\circ\tau_{p,\alpha}
\]
is locally defined on \(D_{S_p,\alpha}(r_p)\), has the
regularity properties of
assumption~\ref{hyp:max3-prod-old-realization}, maps this domain
into the \(H^\sigma\) ball of radius \(r_0\), and obeys
\begin{equation}\label{eq:max3-prod-recursive-scalar}
 H_{\mathrm{fix}}\circ\Psi_{p,\alpha}=H^{(p)}_\alpha
       \quad\text{near }D_{S_p,\alpha}(r_p/8).
\end{equation}
The images of its finite-dimensional reference rotations are the
global solutions in
\eqref{eq:max3-prod-finite-solution}.

At each step, constants $C_p$ and $R_p>0$ are fixed independently of
$r_p$ and the admissible parameter. The radii can be chosen below $R_p$
to satisfy the hypotheses of the opening step and
\begin{equation}\label{eq:max3-prod-radius-targets}
 \begin{gathered}
 C_pr_p^{0.01\delta}\leq\vartheta_*2^{-p},\qquad
 C_p(r_p^{2-6\delta}+r_p^{1-6\delta})\leq2^{-(p+6)},\\
 \langle j_p\rangle^{2\sigma}r_p^{2-2\delta}
       \leq2^{-(p+6)}\min_{0\leq q<p}r_q^2,\qquad
 r_p<r_{p-1}/2,\\
 \langle j_p\rangle^\sigma r_p^{1-\delta}
       \leq r_p^{1-1.01\delta}/100,\qquad
 r_p^{1-1.01\delta}\leq r_{p-1}/100,\qquad
 r_p^{-\delta}\geq100.
 \end{gathered}
\end{equation}
At the first step the preceding domain is the ball of radius $r_0$.
The resulting families satisfy the hypotheses of
Proposition~\ref{prop:xhigh3-assembly}, with $M=\kappa=1/96$. In
particular, both the size and the global Lipschitz seminorm of the
common and $\ell^1$ tangential corrections are at most $1/96$.
\end{proposition}

\begin{proof}
\emph{Inductive hypotheses.}
\label{par:max3-prod-invariant}
By \eqref{eq:max3-prod-width-sequence},
\[
 h_p\le1/108<1/32,\qquad
 \eta_{p-1}-3h_p-24h_p=\eta_p\ge\eta_*+3/2.
\]
Thus the widths remain uniformly above $\eta_*$ after preparation and
inner normalization. We maintain the following inductive hypotheses.
\begin{enumerate}
\item[(i)] The Borel set \(O_p\subset B^{(p)}\) lies in
\(O_{p-1}\times I_p\) and satisfies
\(|O_p|/|B^{(p)}|\geq1-\vartheta_*\sum_{q=1}^p2^{-q}\).
Every \(Z^{(q)}\) is a fixed coefficient
family on its original domain \(O_q\), and \(\tau_q\) is the
fixed map constructed there.

\item[(ii)] The coefficient identity is
\(H^{(p)}=H_2+Z_4^{\mathrm{base}}+\sum_{q\leq p}Z^{(q)}+R^{(p)}\)
in current coordinates, with
\eqref{eq:max3-prod-new-norms} at \((r_p,\eta_p)\), and
\begin{equation}\label{eq:max3-prod-retained-sequence}
 d_q^{\mathrm{ret}}\leq2^{-q}/96.
\end{equation}
The current inherited quartic and remainder have finite uniform
total norms at \((S_p,\eta_p,r_p)\). The constant term of
the current diagonal Hamiltonian is bounded and globally Lipschitz.

\item[(iii)] The finite composition $\Psi_{p,\alpha}$
is locally defined on \(D_{S_p,\alpha}(r_p)\), has all the
properties of assumption~\ref{hyp:max3-prod-old-realization},
and maps this domain into the radius-\(r_0\) ball in
\(H^\sigma\). It satisfies
\eqref{eq:max3-prod-recursive-scalar}, and its finite solution
curves are those of \eqref{eq:max3-prod-finite-solution}.
\end{enumerate}

For a sufficiently small next radius,
Proposition~\ref{prop:repeatable-opening} preserves (i)--(iii).
Indeed, the new integrable quartic is
$W_{\mathrm{bar}}+\Pi_4Z^+$, whose total norm is finite by
\eqref{eq:max3-prod-prepared-bounds} and
\eqref{eq:max3-prod-new-norms}. The identity
\eqref{eq:max3-prod-new-normal-form} preserves the preceding action
polynomials and their parameter dependence. Finite sums and evaluation
on a finite tail preserve the boundedness and Lipschitz bounds of the
constant term. The frequency estimates follow from
Section~\ref{sec:frequency-replacement} and
\eqref{eq:max3-prod-retained-sequence}; the regularity and domain
assertions follow from parts (c)--(d) of
Proposition~\ref{prop:repeatable-opening}.

\medskip\noindent\textit{The first opening.}
\label{par:max3-prod-initialization}
Apply Theorem~\ref{thm:max-initial-cutoff} at
\(\eta_{\mathrm{in}}=\eta_1+24h_1\). For this fixed width, it gives a cutoff transformation with generator
norm $O(\rho^4)$ and no parameter exclusion. The Hamiltonian identity
holds on the new radius-$3r_1$ annulus, and the opened remainder
satisfies $\|E_\xi\|^{\mathrm{tot}}=O(\rho^6)$. The pure normal quartic
has zero integrable projection, while the integrable quartic is the
opening of $Z_4^{\mathrm{base}}$. Hence
\begin{equation}\label{eq:max3-prod-base-split}
 U=I\Pi_{\leq2}E_\xi,\qquad V_4=I\Pi_4E_\xi,\qquad
 I\Pi_4P_0=W_{\mathrm{base}}+V_4.
\end{equation}
There are no earlier integrable corrections, so take $P=0$ in
Proposition~\ref{prop:max4-retained-bootstrap}. For small \(r_1\),
\(u,v\leq\epsilon\) and
\(K_{\mathrm{initial}}\leq Cr_1^4\leq r_1^{4-\delta/4}\).
Proposition~\ref{prop:max3-concrete-inner} then gives the inner normal
form, with new correction
\eqref{eq:max3-prod-correction-subtraction}--\eqref{eq:max3-prod-correction-expansion}.

The initial cutoff excludes no parameters. With
\(\ell_1=|I_1|=r_1^{2-2\delta}/2\) and
\(\kappa_1=\ell_1r_1^{\delta/2}\), as in
\eqref{eq:max3-prod-kappa}, the measure estimate in
Proposition~\ref{prop:max3-concrete-inner} gives
\[
 \frac{|I_1\setminus O_1|}{|I_1|}
 \leq\frac{3C_D\kappa_1}{\ell_1}
 =3C_Dr_1^{\delta/2}\longrightarrow0\qquad(r_1\longrightarrow0).
\]
Here \(C_D\) depends only on the fixed first-stage data and the uniform
frequency bounds, and is independent of \(r_1\). Thus the first loss
is at most \(C_1r_1^{0.01\delta}\) for \(C_1\geq3C_D\) and
\(r_1\leq1\). The first measure condition in
\eqref{eq:max3-prod-radius-choice} therefore makes it at most
\(\vartheta_*/2\), together with all the preceding smallness conditions.

The preliminary generator of
Theorem~\ref{thm:max-initial-cutoff} satisfies the hypotheses for the case of an empty preceding
tangential set in Theorem~\ref{thm:physical-opening}, with
\(P\leq C\rho^4\), the same post-opening bound \(Q\), and
preceding radius \(r_0\). That theorem gives a transformation
on \(D_{S_1,\alpha}(r_1)\) with image in the ball of radius \(r_0\). The initial and inner Hamiltonian identities give
\eqref{eq:max3-prod-recursive-scalar} for $p=1$, with the regularity in
(iii). The coefficient estimates hold at every required width.

\medskip\noindent\textit{Choice of the radii.}
\label{par:max3-prod-recursive-radii}
Suppose all data through stage \(p-1\) have been fixed.
Fix \(j_p,h_p,\eta_p\), the relative measure bound
\(\vartheta_*2^{-p}\), and the correction bound \(2^{-p}/96\).
The preceding bounds, the set $S_p$, and the prescribed width and
cutoff determine the preparation, KAM and flow constants. Choose $C_p$
to dominate the measure and displacement constants. It is independent
of $r_p$ and the admissible parameter.

For \(p>1\), impose \eqref{eq:max3-prod-strict-reserve} with
\(R=r_{p-1}\). For \(p=1\), require the image to lie in the radius-\(r_0\) ball, as above. Choose \(r_p\) below the fixed
finite-stage thresholds and impose
\begin{equation}\label{eq:max3-prod-radius-choice}
 \begin{aligned}
 C_pr_p^{0.01\delta}&\leq\vartheta_*2^{-p},\\
 C_p(r_p^{2-6\delta}+r_p^{1-6\delta})&\leq2^{-(p+6)},\\
 m_p:=\langle j_p\rangle^{2\sigma}r_p^{2-2\delta}
      &\leq2^{-(p+6)}\min_{0\leq q<p}r_q^2,\\
 6K_{S_p}r_p^{2-6.5\delta}&\leq2^{-p}/96,\qquad r_p<r_{p-1}/2,
 \end{aligned}
\end{equation}
where \(K_{S_p}=200(1+\max S_p)(10+3L_*+2|S_p|)\).
All constants and all data from earlier stages are fixed before
\(r_p\) is chosen. In particular, the measure condition follows from
\[
 r_p\leq\left(\frac{\vartheta_*2^{-p}}{C_p}\right)^{1/(0.01\delta)}.
\]
This is a positive upper bound, which can be imposed simultaneously
with the other conditions by decreasing \(r_p\).
We also require
\begin{equation}\label{eq:max3-prod-amplitude-choice}
 \langle j_p\rangle^\sigma r_p^{1-\delta}
     \leq r_p^{1-1.01\delta}/100,\qquad
 r_p^{1-1.01\delta}\leq r_{p-1}/100,\qquad
 r_p^{-\delta}\geq100.
\end{equation}
The first condition follows from
$\langle j_p\rangle^\sigma r_p^{0.01\delta}\to0$; the last ensures
that the new amplitude exceeds $40r_p$. All these conditions, including
the smallness conditions and the strict inclusions of domains at this stage, remain valid
under a further decrease of \(r_p\).

\medskip\noindent\textit{Tail and measure estimates.}
\label{par:max3-prod-tails-measure}
For each fixed \(q\) and each later step \(P\),
the third inequality in \eqref{eq:max3-prod-radius-choice} yields
\begin{equation}\label{eq:max3-prod-tail-sum}
 \sum_{q<l\leq P}\sup I_l
 \leq\sum_{q<l\leq P}m_l
 \leq r_q^2\sum_{l>q}2^{-(l+6)}<r_q^2/64.
\end{equation}
Thus \eqref{eq:max3-prod-old-tail} holds on the full boxes with
$t_q=r_q^2/2$. The same bound controls the infinite tail. Choosing the
next interval shorter than half each remaining margin preserves these
inequalities.

The loss estimate~\eqref{eq:max3-prod-loss} and the first
inequality in \eqref{eq:max3-prod-radius-choice}, together with the
first-stage estimate above, give
\begin{equation}\label{eq:max3-prod-measure-induction}
 \begin{aligned}
 \frac{|O_p|}{|B^{(p)}|}
 &\geq\frac{|O_{p-1}|}{|B^{(p-1)}|}-\vartheta_*2^{-p}\\
 &\geq1-\vartheta_*\sum_{q=1}^p2^{-q}
 =1-\vartheta_*(1-2^{-p})\geq1-\vartheta_*.
 \end{aligned}
\end{equation}
Since \(O_0\) has normalized measure one and
\(1-\vartheta_*\geq15/16>1/2\), the preparation at every step can use the
same lower bound \(m_*=1/2\). Moreover,
\[
 \sum_{q\ge1}d_q^{\mathrm{ret}}\le1/96<1/48,
\]
so the frequency hypotheses hold at every step.

Set $\Psi_{p,\beta}=\Psi_{p-1,\beta|_{S_{p-1}}}\circ\tau_{p,\beta}$.
The coefficient formulas preserve measurability, and
Theorem~\ref{thm:physical-opening} gives joint Borel dependence of the maps.
This completes the induction.

\medskip\noindent\textit{Convergence hypotheses.}
\label{par:max3-prod-outer-interface}
Take the domains $D_p$ and constants $C_p$ in
Proposition~\ref{prop:xhigh3-assembly} as above.
Hypothesis~\ref{hyp:outer-parameters} follows from
\eqref{eq:max3-prod-loss}, the first-stage measure estimate, and
\eqref{eq:max3-prod-radius-choice};
\ref{hyp:outer-maps} follows from \eqref{eq:max3-prod-metric} and
\eqref{eq:max3-prod-path-budgets}--\eqref{eq:max3-prod-path-endpoints};
and \ref{hyp:outer-pde} follows from \eqref{eq:max3-prod-finite-solution}.
Here the previous normal projection $Q_{S_{p-1}}$ includes $j_p$,
and the finite solutions have the stronger regularity
$C(H^\sigma)\cap C^1(H^{s+1})$.

For the frequency hypothesis, apply the summation in
Section~\ref{sec:frequency-replacement} to the retained coefficient families and
the tails in \eqref{eq:max3-prod-tail-sum}. For the common correction $c_p$ and tangential correction $f^{(p)}$,
with $\|(c,f)\|_{\mathbb R\oplus\ell^1}=|c|+\|f\|_{\ell^1}$, it gives
\[
 \begin{aligned}
 \|(c_p,f^{(p)})\|_{\mathbb R\oplus\ell^1}^{\mathrm{lip}}
   &\le\sum_{q\le p}d_q^{\mathrm{ret}}\le1/96,\\
 \|(c_p,f^{(p)})\|_{\mathbb R\oplus\ell^1}^{\sup}
   &\le\sum_{q\le p}d_q^{\mathrm{ret}}r_q^2\le1/96.
 \end{aligned}
\]
The frequencies in the original equation are exactly those of the linear part
of \eqref{eq:max3-prod-new-normal-form}:
\begin{equation}\label{eq:max3-prod-final-frequencies}
 \omega_j^{(p)}=j^3+12j\sum_{i\in S_p}\alpha_i-6j\alpha_j
                    +jc_p+jf_j^{(p)}.
\end{equation}
Thus \(M=1/96\) and \(\kappa=1/96\) meet that hypothesis.

Finally, \eqref{eq:max3-prod-radius-choice} gives the required
smallness, \eqref{eq:max3-prod-strict-reserve} the nesting of domains,
and \eqref{eq:max3-prod-amplitude-choice} the amplitude bounds. At the
first step the preceding domain is the $H^\sigma$ ball of radius $r_0$.
Thus all hypotheses of Proposition~\ref{prop:xhigh3-assembly} hold.
\end{proof}

\section{The limiting torus}
\label{sec:xhigh3-outer}
We now pass to the limit in the finite-dimensional tori of
Proposition~\ref{prop:coherent-production}.

Fix \(s\geq2\), put \(\sigma=s+4\), and enumerate the prescribed infinite set
\(S=\{j_1<j_2<\cdots\}\), where \(\sum_{j\in S}j^{-1}<\infty\).
Set \(S_p=\{j_1,\ldots,j_p\}\), \(S_0=\varnothing\), with the Fourier conventions of Section~\ref{sec:hamiltonian}.
All Fourier projections include the conjugate negative modes.
Complements of tangential sets are taken in $\mathbb N$.

\subsection{Convergence criterion}

Put \(\delta=10^{-3}\), fix \(0<\vartheta_*\leq1/16\) and
\(0<r_0\leq1\), take \(0<r_p\leq1\), and write
\[
 b_{j_p}=r_p^{2-2\delta},\qquad
 m_p=\langle j_p\rangle^{2\sigma}b_{j_p},\qquad
 I_p=[b_{j_p}/2,b_{j_p}],\qquad
 \mathcal B_p=\prod_{q=1}^pI_q.
\]
For finite constants \(C_p\geq1\), define
\[
 d_p=C_pr_p^{2-6\delta},\quad e_p=C_pr_p^{1-6\delta},\quad
 A_p=C_pr_p^{4-6\delta},\quad B_p=C_pr_p^{2-6\delta},\quad
 \beta_p=C_pr_p^{0.01\delta}.
\]
The quantities $A_p,B_p$ bound the action and normal displacements in
$H^\sigma$. Assume
\begin{align}
 d_p+e_p&\leq2^{-p-6},&
 \beta_p&\leq\vartheta_*2^{-p},&
 m_p&\leq2^{-p-6}\min_{0\leq q<p}r_q^2,
 \label{eq:xhigh3-scales}\\
 r_p^2+A_p&\leq r_{p-1}^2,&
 \sqrt{m_p+2r_p^2}+B_p&\leq r_{p-1}.
 \label{eq:xhigh3-nesting}
\end{align}
For \(\alpha\in\mathcal B_p\), define
\[
 D_p(\alpha)=
 \left\{u\in H^\sigma_0:
 \sum_{j\in S_p}\langle j\rangle^{2\sigma}
       \big||u_j|^2-\alpha_j\big|\leq r_p^2,\quad
 \|\Pi_{S_p^c}u\|_\sigma\leq r_p\right\},
 \qquad D_0=\{\|u\|_\sigma\leq r_0\}.
\]

\begin{proposition}[The limiting torus]
\label{prop:xhigh3-assembly}
Assume \eqref{eq:xhigh3-scales}--\eqref{eq:xhigh3-nesting}. Suppose that,
for every $p\ge1$, the following properties hold.

\begin{enumerate}[label=(H\arabic*),ref=H\arabic*]
\item\label{hyp:outer-parameters}
\textit{Parameter sets and measure estimates.}
There are Borel sets \(O_p\subset\mathcal B_p\), with \(O_0\) the
one-point empty parameter space, such that
\[
 O_p\subset O_{p-1}\times I_p,\qquad
 \operatorname{Leb}((O_{p-1}\times I_p)\setminus O_p)
       \leq\beta_p\operatorname{Leb}(\mathcal B_p).
\]

\item\label{hyp:outer-maps}
\textit{Transformations and domain estimates in $H^\sigma$.}
For every \(\alpha\in O_p\) there is a map
\(\tau_{p,\alpha}:D_p(\alpha)\to H^\sigma_0\) satisfying, for every pair
\(u,v\in D_p(\alpha)\),
\[
 \|\tau_{p,\alpha}u-u\|_s\leq d_p,\qquad
 e^{-e_p}\|u-v\|_s
 \leq\|\tau_{p,\alpha}u-\tau_{p,\alpha}v\|_s
 \leq e^{e_p}\|u-v\|_s.
\]
If \(w=\tau_{p,\alpha}u\), it also satisfies
\[
 \sum_{j\in S_{p-1}}\langle j\rangle^{2\sigma}
       \big||w_j|^2-|u_j|^2\big|\leq A_p,\qquad
 \|\Pi_{S_{p-1}^c}w\|_\sigma
       \leq\|\Pi_{S_{p-1}^c}u\|_\sigma+B_p.
\]
The maps from earlier steps are extended to the new parameter sets by
restriction to their original coordinates.

\item\label{hyp:outer-pde}
\textit{Finite-dimensional invariant tori.}
There are Borel frequency maps
\(\omega^{(p)}:O_p\to\mathbb R^{S_p}\). Let
\[
 z_{\alpha,p}(\theta)_j=
 \begin{cases}\sqrt{\alpha_j}e^{i\theta_j},&j\in S_p,\\
 0,&j>0,\ j\notin S_p.\end{cases}
\]
For every \(\alpha\in O_p\) and every \(\theta\in\mathbb T^{S_p}\), the curve
\[
 \big(\tau_{1,\alpha|_{S_1}}\circ\cdots\circ\tau_{p,\alpha}\big)
 z_{\alpha,p}(\theta+t\omega^{(p)}(\alpha))
\]
is a global solution, in
\(C(\mathbb R;H^{s+3}_0)\cap C^1(\mathbb R;H^s_0)\), of the same fixed equation
\[
 u_t=F(u):=-u_{xxx}+\partial_x(2u^3+u^5).
\]

\item\label{hyp:outer-frequencies}
\textit{Frequency estimates.}
There are constants \(M<\infty\), \(0\leq\kappa<3/2\), and Borel
functions \(c_p:O_p\to\mathbb R\), \(f^{(p)}:O_p\to\mathbb R^{S_p}\), such that
\[
 \omega^{(p)}_j(\alpha)
 =j^3+12j\sum_{i\in S_p}\alpha_i-6j\alpha_j
      +jc_p(\alpha)+jf^{(p)}_j(\alpha),
 \qquad
 |c_p(\alpha)|+\|f^{(p)}(\alpha)\|_{\ell^1(S_p)}\leq M,
\]
and for all \(\alpha,\alpha'\in O_p\),
\[
 |c_p(\alpha)-c_p(\alpha')|
 +\|f^{(p)}(\alpha)-f^{(p)}(\alpha')\|_{\ell^1(S_p)}
 \leq\kappa\|\alpha-\alpha'\|_{\ell^1(S_p)}.
\]
\end{enumerate}

Then \(\sum_{j\in S}\langle j\rangle^{2\sigma}b_j<\infty\). For the product
probability
\[
 \mu=\bigotimes_{p\geq1}\operatorname{Unif}(I_p)
 \quad\hbox{on}\quad \mathcal B=\prod_{p\geq1}I_p,
\]
there is a Borel set $G$ satisfying
\[
 \mu(G)\geq1-\sum_{p\geq1}\beta_p\geq1-\vartheta_*
\]
such that, for each $\xi\in G$, the following holds.
Define $z_\xi$ on all of $S$ by the
same Fourier formula and put $T_\xi=z_\xi(\mathbb T^S)$.
There is an embedding $\Psi_\xi:T_\xi\to H^s_0$ with
\[
 C^{-1}\|u-v\|_s\leq\|\Psi_\xi u-\Psi_\xi v\|_s\leq C\|u-v\|_s,
 \qquad C=\exp\Big(\sum_{p\geq1}e_p\Big)\leq e^{1/64},
\]
uniformly in \(\xi\). Its image is compact and homeomorphic to
\(\mathbb T^{\mathbb N}\) in \(H^s\) and in \(H^{s+3}\), its parametrization
is continuous into \(H^{s+3}_0\), and every image point has \(H^{s+4}\) norm
at most \(r_0\).

There is a real sequence $\omega(\xi)$, with
$\omega_j(\xi)=\lim_{p\to\infty}\omega_j^{(p)}(\xi|_{S_p})$ for each $j\in S$,
such that
\[
 U_\theta(t)=\Psi_\xi z_\xi(\theta+t\omega(\xi))
\]
is a global solution of the fixed equation in
\(C(\mathbb R;H^{s+3}_0)\cap C^1(\mathbb R;H^s_0)\), for every angle
\(\theta\). It is Bohr almost periodic in \(H^{s+3}_0\), satisfies
\(\sup_t\|U_\theta(t)\|_{s+4}\leq r_0\), and every forward orbit is dense in
the image torus. Moreover,
\[
 k\cdot\omega(\xi)\ne0
 \qquad\text{for every }0\ne k\in\mathbb Z^{(S)},
\]
where $\mathbb Z^{(S)}$ denotes the finitely supported integer sequences.
\end{proposition}

\medskip\noindent\textit{Choice of the radii.}
Proposition~\ref{prop:coherent-production} fixes $C_p$ and an upper
bound for $r_p$ before choosing $r_p$. All powers of $r_p$ are positive, so
\eqref{eq:xhigh3-scales}--\eqref{eq:xhigh3-nesting} can be imposed
recursively for every prescribed $0<\vartheta_*\leq1/16$.
It also gives \ref{hyp:outer-frequencies} with
$M=\kappa=1/96$.

\subsection{Convergence and invariance}

For \(u\in D_p(\alpha)\),
\[
 \|\Pi_{S_{p-1}^c}u\|_\sigma^2\leq m_p+2r_p^2.
\]
The two endpoint estimates in $H^\sigma$ and
\eqref{eq:xhigh3-nesting} therefore give
\[
 \tau_{p,\alpha}(D_p(\alpha))
 \subset D_{p-1}(\alpha|_{S_{p-1}}).
\]
Thus all finite compositions are defined, by
Theorem~\ref{thm:physical-opening} and
Proposition~\ref{prop:repeatable-opening}.

\begin{proof}[Proof of Proposition \ref{prop:xhigh3-assembly}]
\emph{The limiting embedding.}
For every \(p\geq0\), the tail condition in
\eqref{eq:xhigh3-scales} gives
\[
 \sum_{q>p}m_q\leq2^{-p-6}r_p^2.
\]
In particular, $\sum_{j\in S}\langle j\rangle^{2\sigma}b_j<\infty$.
Let $\widetilde O_p$ be the cylinder over $O_p$ in $\mathcal B$.
By \ref{hyp:outer-parameters}, these sets are decreasing. Thus, by
continuity of the probability measure from above,
\[
 G_0=\bigcap_{p\geq1}\widetilde O_p,\qquad
 \mu(G_0)
 =\lim_{p\to\infty}\frac{\operatorname{Leb}(O_p)}
                          {\operatorname{Leb}(\mathcal B_p)}
 \geq1-\sum_{p\geq1}\beta_p\geq1-\vartheta_*.
\]

Fix \(\xi\in G_0\), suppress parameter subscripts, and put
\(\Psi_p=\tau_1\circ\cdots\circ\tau_p\).
The domain inclusions give
\[
 \Psi_p(D_p)\subset D_0,\qquad
 C^{-1}\|u-v\|_s\leq\|\Psi_pu-\Psi_pv\|_s\leq C\|u-v\|_s
 \quad(u,v\in D_p).
\]
The tail bound gives continuity of $z_\xi:\mathbb T^S\to H^\sigma$.
Positivity of the amplitudes gives injectivity, hence a homeomorphism
onto $T_\xi$. Both $T_\xi$ and its $S_p$-truncation lie in $D_p$.
For $u\in T_\xi$, the points $u,\tau_{p+1}u$ belong to $D_p$, so
\[
 \sup_{u\in T_\xi}\|\Psi_{p+1}u-\Psi_pu\|_s\le Cd_{p+1},\qquad
 \sup_{u\in T_\xi}\|\Psi_q u-\Psi_pu\|_s
       \le C\sum_{l=p+1}^q d_l\longrightarrow0.
\]
Thus $\Psi_p|_{T_\xi}\to\Psi_\xi$ in $C(T_\xi;H^s)$, with the
same bi-Lipschitz bounds. Lower semicontinuity gives
$\|\Psi_\xi u\|_\sigma\le r_0$.
The Fourier interpolation inequality
\[
 \|h\|_v\leq
 \|h\|_s^{(\sigma-v)/(\sigma-s)}
 \|h\|_\sigma^{(v-s)/(\sigma-s)},\qquad s\leq v<\sigma,
\]
gives convergence in $C(T_\xi;H^{s+3})$. Compactness and injectivity
give the asserted homeomorphisms in both topologies.

\smallskip\noindent
\emph{The limiting equation.}
The Sobolev algebra estimate
$\|uv\|_{H^t}\le C_t\|u\|_{H^t}\|v\|_{H^t}$, $t>1/2$, gives
\[
 \|F(u)-F(v)\|_s
 \le C_{s,R}\|u-v\|_{s+3}
 \qquad(\|u\|_{s+3},\|v\|_{s+3}\le R).
\]
Let $K=\sup_{u\in D_0}\|F(u)\|_s<\infty$.
The exact finite solutions and the inverse metric inequality imply
\[
 \|z_{\xi,p}(\theta+t\omega^{(p)})-z_{\xi,p}(\theta)\|_s
       \leq CK|t|.
\]
Differentiating this finite Fourier curve at zero yields
\[
 \sum_{j\in S_p}\langle j\rangle^{2s}\xi_j
                  |\omega_j^{(p)}|^2\leq C^2K^2.
\]
Each frequency coordinate is bounded. Choose a diagonal subsequence
converging coordinatewise to $\omega$. For every $T>0$ and $h\le p$,
\[
 \begin{split}
 &\sup_{|t|\le T,\,\theta}
 \|z_{\xi,p}(\theta+t\omega^{(p)})-z_\xi(\theta+t\omega)\|_s^2\\
 &\quad\le T^2\sum_{j\in S_h}\langle j\rangle^{2s}\xi_j
                    |\omega_j^{(p)}-\omega_j|^2
             +4\sum_{q>h}m_q.
 \end{split}
\]
First let $p\to\infty$ along the subsequence and then $h\to\infty$.
The metric bounds and uniform convergence of $\Psi_p|_{T_\xi}$ give
convergence of the finite solutions to
$U_\theta(t)=\Psi_\xi z_\xi(\theta+t\omega)$, locally uniformly in
time in $H^s$. Interpolation gives convergence in $H^{s+3}$, so
\[
 U_\theta(t)=U_\theta(0)+\int_0^t F(U_\theta(\tau))\,d\tau
 \quad\text{in }H^s.
\]
Consequently $U_\theta\in C^1(\mathbb R;H^s)$ and
$\partial_tU_\theta=F(U_\theta)$.

For uniqueness, put $w=u-v$ and
\[
 A=2(u^2+uv+v^2)+u^4+u^3v+u^2v^2+uv^3+v^4.
\]
Then \(w_t=-w_{xxx}+\partial_x(Aw)\), and periodic integration by parts gives
\[
 \frac12\frac{d}{dt}\|w\|_{L^2}^2
       =\frac12\int_{\mathbb T}A_xw^2\,\frac{dx}{2\pi}.
\]
Since $A_x\in L^\infty$ locally in time, Gronwall's inequality gives
uniqueness. Two frequency accumulation points therefore generate the
same solution of the perturbed mKdV equation. Injectivity of $\Psi_\xi$ and positivity of
the amplitudes identify their rotations and hence their frequencies.
Thus the full sequence converges coordinatewise, with a Borel limit on
$G_0$.

\smallskip\noindent
\emph{Limiting frequencies.}
By \ref{hyp:outer-frequencies}, $c_p$ is bounded. Along a subsequence
with $c_p\to c$, the frequency formula gives
\[
 f_j^{(p)}\longrightarrow
 f_j=\omega_j/j-j^2-12\sum_{i\in S}\xi_i+6\xi_j-c .
\]
Fatou's lemma gives $|c|+\|f\|_{\ell^1(S)}\le M$. Two distinct limits
of $c_p$ would produce reduced limits whose difference is a nonzero
constant sequence, contrary to $f\in\ell^1(S)$. Hence $c_p\to c$ and
$f_j^{(p)}\to f_j$ for every $j$. Applying Fatou's lemma to the
parameter-difference estimate gives
\[
 |c(\xi)-c(\xi')|+\sum_{j\in S}|f_j(\xi)-f_j(\xi')|
          \leq\kappa\|\xi-\xi'\|_{\ell^1(S)} .
\]
In particular
\[
 \omega_j=j^3+12j\sum_{i\in S}\xi_i-6j\xi_j+jc(\xi)+jf_j(\xi).
\]

\smallskip\noindent
\emph{Rational independence.}
Fix a nonzero finitely supported integer sequence \(k\), and set
\[
 w_j=jk_j,\quad P=\sum_jw_j,\quad W=\max_j|w_j|,\quad
 A_k=\max_{m\in\operatorname{supp}k}|2P-w_m|.
\]
If \(d=|\operatorname{supp}k|\), then
\(\sum_m(2P-w_m)=(2d-1)P\), so \(|P|\leq A_k\) and \(W\leq3A_k\).
Also \(A_k>0\). Choose a maximizing coordinate \(m\). Two points of \(G_0\)
which differ only in that coordinate satisfy
\[
 |\Delta(k\cdot\omega)-6(2P-w_m)\Delta\xi_m|
       \leq4\kappa A_k|\Delta\xi_m|.
\]
Since $4\kappa<6$, each $\xi_m$-fiber contains at most one zero of
$k\cdot\omega$. The resonance set is Borel, and the coordinate measure
is nonatomic; Fubini's theorem therefore gives
\[
 \mu\{\xi\in G_0:k\cdot\omega(\xi)=0\}=0.
\]
Removing the countable union over nonzero finitely supported $k$ yields
$G\subset G_0$ with $\mu(G)=\mu(G_0)$.

\smallskip\noindent
\emph{Density and almost periodicity.}
By the finite-dimensional Kronecker theorem, every finite projection of
each forward orbit is dense. Since cylinder sets form a basis of the
product topology, the orbit is dense in $\mathbb T^S$ and hence in
the image torus.

By \cite[Lemma~10.2]{BGR25}, applied in the real normed space
$H^{s+3}_0$ to the homeomorphism
$\varphi\mapsto\Psi_\xi z_\xi(\theta+\varphi)$ and the rationally
independent frequency sequence, $U_\theta$ is Bohr almost periodic.
\end{proof}

\subsection{Proof of the main theorem}
\begin{proof}[Proof of Theorem~\ref{thm:main}]
Fix $0<\vartheta<1$, take $r_0=\varepsilon\le1$, set
$\vartheta_*=\min\{\vartheta,1/16\}$, and write
$S_\infty=\{j_1<j_2<\cdots\}$.
Proposition~\ref{prop:coherent-production}, applied with this
$\vartheta_*$, gives parameter sets, radii
and transformations satisfying
\eqref{eq:xhigh3-scales}--\eqref{eq:xhigh3-nesting} and
\ref{hyp:outer-parameters}--\ref{hyp:outer-frequencies}, with
$M=\kappa=1/96$. Set $b_{j_p}=r_p^{2-2\delta}$; these numbers may
depend on both $\varepsilon$ and $\vartheta$.
Proposition~\ref{prop:xhigh3-assembly} gives
\[
 \mathbb P(\mathcal G)\geq1-\vartheta_*\geq1-\vartheta,
 \qquad
 \sup_{w\in\mathcal T_\xi}\|w\|_{s+4}\leq r_0=\varepsilon,
\]
together with all the remaining conclusions of the theorem.
\end{proof}

\section*{Use of AI tools}
Generative AI tools were used to assist with proof development, checking,
and editorial preparation of this manuscript.
\bibliographystyle{amsplain}
\providecommand{\bysame}{\leavevmode\hbox to3em{\hrulefill}\thinspace}
\providecommand{\MR}{\relax\ifhmode\unskip\space\fi MR }
\providecommand{\MRhref}[2]{%
  \href{http://www.ams.org/mathscinet-getitem?mr=#1}{#2}
}
\providecommand{\href}[2]{#2}

\end{document}